\documentclass{article}
\newtheorem{theorem}{Theorem}
\newtheorem{lemma}{Lemma}

\usepackage{graphicx}
\usepackage{amsfonts}
\DeclareSymbolFont{AMSa}{U}{msa}{m}{n}
\DeclareMathDelimiter\ulcorner{\mathopen} {AMSa}{"70}{AMSa}{"70}
\DeclareMathDelimiter\urcorner{\mathclose}{AMSa}{"71}{AMSa}{"71}
\makeatletter
\def\uufill{$\m@th\mathopen\ulcorner\mkern-7mu%
  \cleaders\hbox{\rule[6pt]{1dd}{1dd}}\hfill
  \mkern-7mu\mathclose\urcorner$}
\def\overbrack#1{\vbox{\m@th\ialign{##\crcr
      \uufill\crcr\noalign{\kern-\p@\nointerlineskip}%
      $\hfil\displaystyle{#1}\hfil$\crcr}}}
\makeatother

\title{Quantum Knots and  Riemann Hypothesis}

\author{Sze Kui Ng
\\ {\small Department of Mathematics and Aplied Mathematics, Hanshan Normal University, China}\\{\small szekuing@yahoo.com.hk}
}

\begin{document}
\date{}
\maketitle
\begin{abstract}
In this paper
we propose a quantum gauge system from which we construct generalized Wilson loops which will be as quantum knots. From quantum knots we give a classification table of knots where knots are one-to-one assigned with an integer such that prime knots are bijectively assigned with prime numbers and the prime number $2$ corresponds to the trefoil knot.

Then by considering the quantum knots as periodic orbits of the quantum system and by the identity of knots with integers and an approach which is similar to the quantum chaos approach of Berry and Keating we derive a trace formula which may be called the von Mangoldt-Selberg-Gutzwiller trace formula. From this trace formula
we then give a proof of the Riemann Hypothesis. For our proof of the Riemann Hypothesis we show that the Hilbert-Polya conjecture holds that there is a self-adjoint operator for the nontrivial zeros of the Riemann zeta function and this operator is the Virasoro energy operator with central charge $c=\frac12$. Our approach for proving the Riemann Hypothesis can also be extended to prove the Extended Riemann Hypothesis. We also investigate the relation of our approach for proving the Riemann Hypothesis with the Random Matrix Theory for $L$-functions.


\end{abstract}

\section{Introduction}\label{sec00}

It is well known that the Jones polynomial as a knot invariant can be derived from a quantum Chern-Simon gauge field theory
\cite{Jon}\cite{Witten}. Inspired by this work
in this paper we shall also propose a quantum gauge model. In this quantum model we generalize the way of defining Wilson loops to construct generalized Wilson loops which will be as quantum knots. From quantum knots we  give a classification table of knots where knots are one-to-one assigned with an integer such that prime knots are bijectively assigned with prime numbers and the prime number $2$ corresponds to the trefoil knot.

Then by considering the quantum knots as periodic orbits of the quantum model and by the identity of knots with integers and an approach which is similar to the quantum chaos approach of Berry and Keating we derive a trace formula which may be called the von Mangoldt-Selberg-Gutzwiller trace formula. From this trace formula
we then give a proof of the Riemann Hypothesis \cite{Sel}-\cite{Kea}.

From the quantum gauge model we first define the classical Wilson loop and Wilson line. Then from the quantum gauge model we derive a definition for the generator of the Wilson line.

Then we derive two quantum Knizhnik-Zamolodchikov (KZ) equations which are dual to each other for the product of quantum Wilson lines.
This quantum KZ equation in dual form may be regarded as a quantum Yang-Mill equation as analogous to the classical Yang-Mill equation derived from the classical Yang-Mill theory since this quantum KZ equation is as the basic quantum equation derived from the quantum gauge model. Solutions of this quantum Yang-Mill equation are then used to construct generalized Wilson loops which are as quantum knots (These quantum knots may be regarded as solitons as similar to the instantons of the classical Yang-Mill equation).
In deriving this quantum KZ equation we first derive a conformal field theory consisting of the Kac-Moody algebra and the Virasoro energy operator and Virasoro algebra.

Then from the quantum knots we derive a knot invariant. From this knot invariant
we give the classification table of knots.

Then the quantum knots as the periodic orbits of the quantum gauge system and the identity of prime knots with prime numbers are as the two basic ingredients for proving the Riemann Hypothesis.

For our proof of the Riemann Hypothesis we show that the Hilbert-Polya conjecture holds that there is a self-adjoint operator for the nontrivial zeros of the Riemann zeta function and this operator is the Virasoro energy operator with central charge $c=\frac12$ \cite{Fra}-\cite{Fuc}. Our approach for proving the Riemann Hypothesis can also be extended to prove the Extended Riemann Hypothesis. We also investigate the relation of our approach for proving the Riemann Hypothesis with the Random Matrix Theory for $L$-functions
\cite{Dys}-\cite{Kea2}.

This paper is organized as follows. In section 2 we give a brief
description of a quantum gauge model of electrodynamics  and its
nonabelian generalization. In this paper we shall consider a
nonabelian generalization with a $SU(2)$ gauge symmetry. With this
quantum  model in section 3 we introduce the definition of classical Wilson loop and Wilson line.

In section 4
we derive the defintion of the generator of the Wilson line. From this definition in section 4 and 5 we derive a conformal field theory which includes the Virasolo enegry operator and Virasolo algebra, the affine
Kac-Moody algebra and the quantum KZ equation in dual form.
In section 6 we compute the solutions of the quantum KZ equation in dual form.
In section 7 we compute the quantum Wilson lines.
In section 8
we represent the braiding
of two pieces of curves by defining the braiding of
 two quantum Wilson lines.
 By this representation in section 10 we define the generalized Wilson loop which will be as a quantum knot.
In section 9
we compute the quantum Wilson loop.
In section 10
we define generalized Wilson loops which will be shown to have properties of the corresponding knot diagram and will be regarded as quantum knots.
In section 11
we give some examples of generalized Wilson loops and show that they have the properties of the corresponding knot diagram and thus may be regarded as quantum knots.
In  section 12
we show
that this generalized Wilson loop is a complete copy of the corresponding knot diagram and thus we may call a generalized Wilson loop as a quantum knot. From quantum knots we have a knot invariant of the form $Tr R^{-m} W(z,z)$ where $W(z,z)$ denotes
a quantum Wilson loop  and $R$ is the
braiding
matrix and is the monodromy of the quantum KZ equation and $m$ is an integer. We show that this knot invariant classifies knots and that knots can be one-to-one assigned with the integer $m$.
In section 13
we give more computations of quantum knots and their knot invariant.
Then in section 14  and 15
with the integer $m$ we give a
classification table
of knots  where we show that
prime knots (and only prime knots) are assigned with prime integer $m$.

Then in section 16
by using the classification table of knots and by considering the quantum knots as the periodic orbits of the quantum gauge model we then, by using the quantum chaos approach of Gutzwiller, Berry and Keating and the von Mangoldt-Selberg approach of proving the Riemann Hypothesis, derive the von Mangoldt-Selberg-Gutzwiller trace formula. From this trace formula we then prove the Riemann Hypothesis. We also generalize this approach to prove the Extended Riemann Hypothesis.
In section 17
we determine that the central charge $c$ for the Riemann zeta function and the Dirichlet $L$-functions is equal to $\frac12$.
Then in section 18
by our approach for proving Riemann Hypothesis we give a commutative diagram relating
Virasoro energy operator and Virasoro algebra, automorphic (or modular) forms and $L$-functions.
In section 19
we show a connection of our approach with the Random Matrix Theory approach for $L$-functions by showing that this connection is from the conformal field theory.

\section{A Quantum Gauge Model}\label{sec2}

We shall first establish a quantum gauge model. This quantum gauge model will be as a physical motivation for introducing operators which will be called Wilson loop and Wilson line as analogous to the Wilson loops in the existing quantum field theories. Then the definition of classical Wilson loop and Wilson line and the definition of a generator $J$ of the Wilson line will be as the basis of the mathematical foundation of this paper (In order to simplify the mathematics of this paper we treat this quantum gauge model as a physical motivation instead of as the mathematical foundation of this paper).

We shall show that the generator $J$ gives an affine Kac-Moody algebra and a Virasoro energy operator $T$ with central charge $c$.
From  $J$ and $T$ we shall derive the quantum KZ equation in dual form which will be regarded as the quantum Yang-Mills equation.
From this quantum KZ equation we then construct generalized Wilson loops which will be as quantum  knots.

Let us
construct a quantum gauge model, as follows. In
probability theory we have the Wiener measure $\nu$ which is a
measure on the space $C[t_0,t_1]$ of continuous functions
\cite{Jaf}. This measure is a well defined mathematical theory for
the Brownian motion and it may be symbolically written in the
following form:
\begin{equation}
d\nu =e^{-L_0}dx
\label{wiener}
\end{equation}
where $L_0 :=
\frac12\int_{t_0}^{t_1}\left(\frac{dx}{dt}\right)^2dt$ is the
energy integral of the Brownian particle and $dx =
\frac{1}{N}\prod_{t}dx(t)$ is symbolically a product of Lebesgue
measures $dx(t)$ and $N$ is a normalized constant.

Once the Wiener measure is defined we may then define other
measures on $C[t_0,t_1]$ as follows\cite{Jaf}. Let a potential
term $\frac12\int_{t_0}^{t_1}Vdt$ be added to $L_0$. Then we have
a measure $\nu_1$ on $C[t_0,t_1]$ defined by:
\begin{equation}
d\nu_1 =e^{-\frac12\int_{t_0}^{t_1}Vdt}d \nu
\label{wiener2}
\end{equation}
Under some condition on $V$ we have that $\nu_1$ is well defined
on $C[t_0,t_1]$. Let us call (\ref{wiener2}) as the Feymann-Kac
formula \cite{Jaf}.

Let us then follow this formula to construct a quantum  model 
electrodynamics, as follows.
Then similar to the formula (\ref{wiener2}) we construct a quantum
model of electrodynamics from the following energy
integral:
\begin{equation}
\begin{array}{rl}
&-\int_{s_0}^{s_1}D ds:= -\int_{s_0}^{s_1}[\frac12\left(\frac{\partial A_1}{\partial
x^2}-\frac{\partial A_2}{\partial x^1}\right)^*
\left(\frac{\partial A_1}{\partial x^2}-\frac{\partial A_2}{\partial x^1}\right)+ 
 \left(
\frac{dZ^*}{ds} +ie_0 A_z Z^*\right)
\left( \frac{dZ}{ds}
-ie_0A_z  Z\right)]ds \\
\end{array}
\label{1.1}
\end{equation}
where the complex variable $Z=Z(z(s))$ and the real variables
$A_1=A_1(z(s))$ and $A_2=A_2(z(s))$ are continuous functions in a
form that they are in terms of a (continuously differentiable)
curve $z(s)=(x^1(s),x^2(s)), s_0\leq s\leq s_1,
z(s_0)=z(s_1)$ in the complex plane where $s$ is a parameter
representing the proper time in relativity.
The complex
variable $Z=Z(z(s))$ represents a field of matter( such as the
electron) ($Z^*$ denotes its complex conjugate) and the real
variables $A_1=A_1(z(s))$ and $A_2=A_2(z(s))$ represent a
connection (or the gauge field of the photon), and $A_z =\sum_{j=1}^2A_j\frac{dx^j}{ds}$,  and $e_0$ denotes
the bare electric charge.

The integral (\ref{1.1}) has the following gauge symmetry:
\begin{equation}
\begin{array}{rl}
Z^{\prime}(z(s))  := Z(z(s))e^{ie_0a(z(s))},  \quad  A'_j (z(s))  := A_j(z(s))+
\frac{\partial a}{\partial x^j} \quad j=1,2
\end{array}
\label{1.2}
\end{equation}
where $a=a(z)$ is a continuously differentiable real-valued
function of $z$.

As the Wiener measure is based on the Banach space $C\lbrack s_0, s_1 \rbrack$, 
the above gauge model 
is based on the Banach space $C^3\lbrack s_0, s_1 \rbrack$ of continuous functions $Z(z(s)), A_j(z(s)), j=1,2, s_0\leq s\leq s_1$
on the one dimensional interval $\lbrack s_0, s_1 \rbrack$.

Similar to the differential $dx$ of the Wiener measure,  the differential of this gauge model is of the form $dA_zdZ^*dZ$. 

Since $D$ is positive and the model is one dimensional
we have that this gauge model is similar to the Wiener measure except that this gauge model has a gauge symmetry. 

An advantage of this  gauge model is that it is free of the difficulty of the degenerate degree of freedom of the gauge symmetry where the two degrees of  $A_1, A_2$ is combined into one degree of $A_z$. In the physics literature the usual way to treat the degenerate degree of freedom of gauge symmetry is to introduce a gauge fixing condition to eliminate the degenerate degree of freedom \cite{Fad}. In this gauge model because it is free of the difficulty of the degenerate degree of freedom of the gauge symmetry that the mathematical difficulty of introducing gauge fixing condition can be avoid.

We have the  following theorem:
\begin{theorem}
The above gauge model is a measure 
defined on the Banach space $C^3\lbrack s_0, s_1 \rbrack$.
\end{theorem}

{\bf Proof}. 
For this gauge model let us choose gauges such that the above gauge model is a measure
defined on the Banach space $C^3\lbrack s_0, s_1 \rbrack$. As an example let us choose a gauge $ (A_1, A_2)$  such that the following form holds:
\begin{equation}
\frac{\partial A_1}{\partial
x^2}-\frac{\partial A_2}{\partial x^1}=\rho\frac{d A_z(z(s))}{ds}
\label{measure3}
\end{equation}
for some constant $\rho>0$.
Then the energy integral is of the following form:
\begin{equation}
\begin{array}{rl}
&-\int_{s_0}^{s_1}D ds= -\int_{s_0}^{s_1}[\frac{\rho^2}2\left(\frac{d A_z}{ds}\right)^2
+ 
 \left(
\frac{dZ^*}{ds} +ie_0 A_z Z^*\right)
\left( \frac{dZ}{ds}
-ie_0A_z  Z\right)]ds \\
\end{array}
\label{me}
\end{equation}
Then the singular pure gauge $ (A_1, A_2)$  part  of  the gauge model is of the form of  the Wiener measure and thus is a well defined measure.

On the other hand  the $Z$ part does not have the singularity of the gauge model and it gives a measure of the form of the above Feynman-Kac type. Thus from the energy integral (\ref{me})  we have a well defined measure on  the Banach space $C^3\lbrack s_0, s_1 \rbrack$. This proves the theorem.  $\diamond$

We remark that this model
 is not
formulated with the four-dimensional space-time but is formulated
with the one dimensional proper time in relativity theory. This gets rid of the difficulty of ultraviolet divergences of the four-dimensional space-time quantum field theory


Similar to the usual Yang-Mills gauge theory we can generalize this gauge model with $U(1)$ gauge symmetry to
nonabelian gauge models. As an illustration let us consider
$SU(2)$ gauge symmetry. Similar to (\ref{1.1}) we consider the
following energy integral:
\begin{equation}
L := 
\int_{s_0}^{s_1} [\frac12 Tr (D_1A_2-D_2A_1)^{*}(D_1A_2-D_2A_1) +
(D_0^*Z^*)(D_0Z)]ds \label{n1}
\end{equation}
where $Z= (z_1, z_2)^{T}$ is a two dimensional complex vector;
$A_j =\sum_{k=0}^{3}A_j^k e_0 t^k $ $(j=1,2)$ where $A_j^k$ denotes a
component of a gauge field $A^k$; $t^k=i T^k$ denotes a
generator of $SU(2)\otimes U(1)$ where $T^k$ denotes a
self-adjoint generator of $SU(2)\otimes U(1)$ 
 (here for simplicity we choose a
convention that the complex $i$ is absorbed by $t^k$ and $t^k$ is
absorbed by $A_j$; and  the notation $A_j$ is with a little
confusion with the notation $A_j$ in the above formulation of
(\ref{1.1}) where $A_j, j=1,2$ are real valued); and
$D_j=\frac{\partial}{\partial x_j}-A_j $ for $j=1,2$;  and
$D_0=\frac{d}{ds}- A_z$ where $A_z=\sum_{j=1}^2A_j\frac{dx^j}{ds}$.

From (\ref{n1}) we can develop a nonabelian gauge model as similar
to that for the above abelian gauge model.
We have that (\ref{n1}) is invariant under the following
gauge transformation:
\begin{equation}
\begin{array}{rl}
Z^{\prime}(z(s)) & :=U(a(z(s)))Z(z(s)) \\
A_j^{\prime}(z(s)) & := U(a(z(s)))A_j(z(s))U^{-1}(a(z(s)))+
 U(a(z(s)))\frac{\partial U^{-1}}{\partial x^j}(a(z(s))),
\quad  j =1,2
\end{array}
\label{n2}
\end{equation}
where $U(a(z(s)))=e^{a(z(s))}$; $a(z(s))=\sum_k  a^k (z(s))e_0 t^k$  for some functions $a^k$.
We shall mainly consider the case that $a$ is a function of the form $a(z(s))
=\sum_k \mbox{Re}\, \omega^k(z(s))e_0 t^k$ where $\omega^k$ are
analytic functions of $z$ (We let
$\omega(z(s)):=\sum_k\omega^k(z(s))e_0 t^k$ and we write
$a(z)=\mbox{Re}\,\omega(z)$).

Then as the above abelian case, this gauge model is a well defined measure on the Banach space 
of continuous functions $Z(z(s)), A_j(z(s)), j=1,2, s_0\leq s\leq s_1$
on the one dimensional interval $\lbrack s_0, s_1 \rbrack$..



From this gauge model we shall derive 
 the quantum KZ equation in dual form which will be regarded as a quantum Yang-Mills equation since its role will be similar to the classical Yang-Mills equation derived from the classical Yang-Mills gauge model.


{\bf Remark}. In this paper the main aim of introducing this quantum gauge model is to derive the quantum KZ equation in dual form which will be regarded as a quantum Yang-Mills equation. From this quantum KZ equation in dual form we then construct quantum knots. From quantum knots we then prove the Riemann Hypothesis.

\section{Dirac-Wilson loop } \label{sec4}

Similar to the Wilson loop in quantum field theory \cite{Witten} from
our quantum theory we introduce an analogue of Wilson loop, as
follows (We shall also call a Wilson loop as a Dirac-Wilson loop).

{\bf Definition}.
 A classical Wilson loop 
$W_{DC}(\gamma)$ 
is defined by :
\begin{equation}
W_{DC}(\gamma):= W_C(z_0, z_1):= Pe^{\int_{\gamma} A_jdx^j}
 \label{n4}
\end{equation}
where $D$ denotes a representation of $SU(2)$ and $C$ means classical (For simplicity we shall omit the notation $D$: $W_{C}(\gamma):=W_{DC}(\gamma)$); $\gamma=z(\cdot)$ is a  closed curve;
 $z_0=z_1$ is a starting and end point of $\gamma=z(\cdot)$;
 and the quantum gauge model is based
on $z(\cdot)$ as specific in the above section (We keep in mind that the definition of the classical Wilson loop $ W_C(z_0, z_1)$ is defined and depended on the whole curve $z(\cdot)$). As usual the notation $P$
in the definition of $W_C(\gamma)$ denotes a path-ordered product
\cite{Witten}\cite{Kau}\cite{Baez}:
\begin{equation}
 Pe^{\int_{\gamma} A_jdx^j}:=\lim_{\max_k \{r_{k+1}-r_k\}\to 0}\prod_{k=0}^n e^{ \sum_{j=1}^2A_j(z(r_k))\Delta x^j(r_k))}
 \label{n4p}
\end{equation}
where $s_0=r_1<r_2<\cdot\cdot\cdot <r_n=s_1$ is a partition of $\lbrack s_0, s_1\rbrack$.  

{\bf Remarks}.

1). We  extend the definition of $W_C(z(\cdot))$ to the case that
$z(\cdot)$ is not a closed curve with $z_0\neq z_1$. When
$z(\cdot)$ is not a closed curve we shall call $W_C(z_0, z_1)$ as a
Wilson line.


2) We shall  use the notation $W(z_0, z_1)$ to denote the quatum version of the classical Wilsion line  $W_C(z_0, z_1)$ and we call $W(z_0, z_1)$ as the quantum Wilson line.

We first have the following theorem on $W_C(z_0,z_1)$:
\begin{theorem}
For a given continuous path $A_i, i=1,2$  on $[s_0, s_1]$
the Wilson line $W_C(z_0,z_1)$ exists on this path and has the
following transition property:
\begin{equation}
W_C(z_0,z_1)=W_C(z_0,z)W_C(z,z_1)
 \label{df2}
\end{equation}
where $W_C(z_0,z_1)$ denotes the Wilson line of a
curve $z(\cdot)$ which is with $z_0$ as the starting
point and $z_1$ as the ending point and $z$ is a
point on $z(\cdot)$ between $z_0$ and $z_1$.
\end{theorem}

{\bf Proof}. We have that $W_C(z_0,z_1)$ is a limit
(whenever exists)
of ordered product of $e^{A_i\triangle x^i}$ and thus can be
written in the following form:
\begin{equation}
\begin{array}{rl}
W_C(z_0,z_1)= & I +
\int_{s^{\prime}}^{s^{\prime\prime}}
A_i(z(s))\frac{dx^i(s)}{ds}ds \\
 & + \int_{s^{\prime}}^{s^{\prime\prime}}
[\int_{s^{\prime}}^{s_1} A_i(z(s_1))\frac{dx^i(s_1)}{ds}ds_1]
A_i(z(s_2))\frac{dx^i(s_2)}{ds}ds_2 +\cdot\cdot\cdot
\end{array}
\label{df3}
\end{equation}
where $z(s^{\prime})=z_0$ and $z(s^{\prime\prime})=z_1$. Then
since $A_i$ are continuous on $[s^{\prime}, s^{\prime\prime}]$ and
$x^i(z(\cdot))$ are continuously differentiable on $[s^{\prime},
s^{\prime\prime}]$ we have that the series in (\ref{df3}) is
absolutely convergent. Thus the Wilson line $W_C(z_0,z_1)$ exists.
Then since $W_C(z_0,z_1)$ is the limit of ordered
product  we can write $W_C(z_0,z_1)$ in the form $W_C(z_0,z)W_C(z,z_1)$
by dividing $z(\cdot)$ into two parts at $z$. This proves the
theorem. $\diamond$



By following the usual approach of deriving a chiral symmetry from a gauge transformation of a gauge field we have the following chiral symmetry which is derived by applying an analytic gauge transformation with an analytic function $\omega$ for the transformation:



\begin{theorem}
Under an analytic gauge
transformation with an analytic function $\omega$
we have the following symmetry:
\begin{equation}
W_C(z_0, z_1) \mapsto W^{\prime}_C(z_0, z_1)=U(\omega(z_1))
W_C(z_0, z_1)U^{-1}(\omega(z_0))
\label{n5}
\end{equation}
where $W_C^{\prime}(z_0, z_1)$ is a Wilson line with gauge field:
\begin{equation}
A_{\mu}^{\prime} =  U(z)  \frac{\partial U^{-1}(z)}{\partial x^{\mu}} + U(z)A_{\mu}U^{-1}(z)
 \label{an5}
\end{equation}
\end{theorem}

{\bf Proof}. 
Let a gauge transformation $U(\omega(z))$ be with an analytic variation 
$\omega$ of the form $\omega=\sum_k\omega^k e_0 t^k$ with the  generators $e_0t^k$.

Then let us prove the symmetry (\ref{n5}) as follows. Let $U(z):=
U(\omega(z(s)))$ and $U^{-1}(z+dz)\approx U^{-1}(z)+\frac{\partial
U^{-1}(z)}{\partial x^{\mu}}dx^{\mu}$ where $dz=(dx^1,dx^2)$. Following
 \cite{Kau} we have
\begin{equation}
\begin{array}{rl}
& U(z)(1+ dx^{\mu}A_{\mu})U^{-1}(z+ dz)\\
= & U(z)U^{-1}(z+ dz)
+ dx^{\mu}U(z)A_{\mu}U^{-1}(z+ dz) \\
\approx & 1+U(z) \frac{\partial U^{-1}(z)}{\partial
x^{\mu}}dx^{\mu}
  + dx^{\mu}U(z)A_{\mu}U^{-1}(z+ dz) \\
\approx & 1+U(z) \frac{\partial U^{-1}(z)}{\partial
x^{\mu}}dx^{\mu}
+ dx^{\mu}U(z)A_{\mu}U^{-1}(z) \\
=:& 1 + dx^{\mu}A_{\mu}^{\prime}
\end{array}
\label{n5b}
\end{equation}

From (\ref{n5b}) we have that (\ref{n5}) holds since (\ref{n5}) is
the limit of ordered product in which the right-side factor $U^{-1}(z_i+ dz_i)$ in (\ref{n5b}) 
is canceled by the left-side
factor $U(z_{i+1})$ of (\ref{n5b}) where $ z_{i+1}=z_i+ dz_i$. This proves the theorem. $\diamond$

As analogous to
the WZW model in
conformal field theory \cite{Fra}\cite{Fuc},
from the above symmetry  we have the following formulas for the
variations $\delta_{\omega}W_C$  (where the notation $\delta_{\omega}$  means the variation with respect to the variations $\omega^k$), and $\delta_{\omega^{\prime}}W_C$  with
respect to this symmetry (\cite{Fra} p.621):
\begin{equation}
\delta_{\omega}W_C(z,z')=W_C(z,z')\omega(z), \quad 
\delta_{\omega^{\prime}}W_C(z,z')=-\omega^{\prime}(z')W_C(z,z')
\label{k2a}
\end{equation}
where $z$ and $z'$ are independent variables and
$\omega^{\prime}(z')=\omega(z)$ when $z'=z$. In (\ref{k2a})
 the
variation $\delta_{\omega}$ is with respect to the $z$ variable while
the variation $\delta_{\omega^{\prime}}$ is with respect to the $z'$ variable. This
two-side-variations when $z\neq z'$ can be derived as follows. For
the left variation we may let $\omega$ be analytic in a
neighborhood of $z$ and  extended as a continuously differentiable
function to a neighborhood of $z'$ such that $\omega(z')=0$ in
this neighborhood of $z'$. Then from (\ref{n5}) we have that the first formula of
(\ref{k2a}) holds. Similarly we may let $\omega^{\prime}$ be
analytic in a neighborhood of $z'$ and extended as a continuously
differentiable function to a neighborhood of $z$ such that
$\omega^{\prime}(z)=0$ in this neighborhood of $z$. Then we have
that the second formula of (\ref{k2a}) holds.

In terms of 
an analytic variation 
$\omega$ of the form $\omega=\sum_k\omega^k t^k$ where the  generators $t^k$ are without the charge $e_0$, 
 from  
 (\ref{k2a}) we have the following   formulas for the
variations $\delta_{\omega}W_C$ and $\delta_{\omega^{\prime}}W_C$:
\begin{equation}
\delta_{\omega}W_C(z,z')=W_C(z,z')\omega(z)e_0, \quad 
\delta_{\omega^{\prime}}W_C(z,z')=-e_0\omega^{\prime}(z')W_C(z,z')
\label{k1c}
\end{equation}

\section{Affine Kac-Moody algebra} \label{sec6}


 We shall derive a quantum loop algebra (or the
affine Kac-Moody algebra) structure from the Wilson line $W_C(z,z')$
for the generator $J$ of $W(z,z')$. To this end let us first
consider the classical case. Since $W_C(z,z')$ is constructed from $
SU(2)$ we have that the mapping $z \to W_C(z,z')$ (We consider
$W_C(z,z')$ as a function of $z$ with $z'$ being fixed) has a loop
group structure \cite{Lus}\cite{Seg}. For a loop group we have the
following generators:
\begin{equation}
J_n^a = t^a z^n \qquad n=0, \pm 1, \pm 2, ...
\label{km1}
\end{equation}
These generators satisfy the following algebra:
\begin{equation}
[J_m^a, J_n^b] =
if_{abc}J_{m+n}^c
\label{km2}
\end{equation}
This is  the so called loop algebra \cite{Lus}\cite{Seg}. Let us
then introduce the following generating function $J$:
\begin{equation}
J(w) = \sum_a J^a(w)=\sum_a j^a(w) t^a
\label{km3}
\end{equation}
where we define
\begin{equation}
J^a(w)= j^a(w) t^a :=
\sum_{n=-\infty}^{\infty}J_n^a(z) (w-z)^{-n-1}
\label{km3a}
\end{equation}

From $J$ we have
\begin{equation}
J_n^a=  \frac{1}{2\pi i}\oint_z dw (w-z)^{n}J^a(w)
\label{km4}
\end{equation}
where $\oint_z$ denotes a closed contour integral  with center $z$. This formula
can be interpreted as that
$J$ is the generator of the loop group and that
$J_n^a$ is the directional generator in the direction
$\omega^a(w)= (w-z)^n$. We may generalize $(\ref{km4})$
to the following  directional generator:
\begin{equation}
  \frac{1}{2\pi i}\oint_z dw \omega(w)J(w)
\label{km5}
\end{equation}
where the analytic function
\begin{equation}
\omega(w)=\sum_a \omega^a(w) t^a
\label{variation}
\end{equation}
 is regarded
as a variation direction and we define
\begin{equation}
 \omega(w)J(w):= \sum_a \omega^a(w)J^a
\label{km5a}
\end{equation}

Then since $W_C(z,z')\in SU(2)$, from the variational formula
(\ref{km5}) for the loop algebra of the loop group of $SU(2)$ we
have that the variation of $W_C(z,z')$ in the direction $\omega(w)$
is given by
\begin{equation}
W_C(z,z')
  \frac{1}{2\pi i}\oint_z dw \omega(w)J(w)
\label{km6}
\end{equation}

Now let us consider the quantum case which is based on the quantum
gauge theory in section 2. For this quantum case we shall define a
quantum generator $J$ which is analogous to the $J$ in
(\ref{km3}). We shall choose the equations (\ref{n8b}) and
(\ref{n6}) as  the equations for defining the quantum generator
$J$. Let us first give a formal derivation of the equation
(\ref{n8b}), as follows.
 Let us consider the
following 
functional integration based on the measure of the gauge model:
\begin{equation}
\langle W(z,z')B(z) \rangle := \int dA_zdZ^{*}dZ  e^{-L}
W_C(z,z')B_C(z) \label{n8a}
\end{equation}
where  
$A_z(z(s))=\sum_{j=1}^2
A_j\frac{d x^j}{ds}$ and
$B_C(z)$ denotes a field from the quantum gauge theory (We
first let $z'$ be fixed as a parameter) and $B(z)$ is the corresponding quantum operator.

Let us  do a calculus of variation on this integral to derive a variational
equation by applying a gauge transformation on (\ref{n8a}) as follows
(We remark that such variational equations are usually called the
Ward identity in the  physics literature).

Let 
$(A_1,A_2,Z)$ be regarded as a coordinate system of the integral
(\ref{n8a}).
Under a gauge transformation (regarded as
a change of coordinate) with gauge function $a(z(s))$ this coordinate
is changed to another coordinate denoted by
$(A_1^{\prime}, A_2^{\prime}, Z^{\prime})$.
As similar to the usual change of variable for integration we have that
the integral  (\ref{n8a}) is unchanged
under a change of variable and we have the following
equality:
\begin{equation}
\begin{array}{rl}
& \int dA_z^{\prime}
dZ^{\prime *}dZ^{\prime}
 e^{-L^{\prime}} W_C^{\prime}(z,z')B_C^{\prime}(z)
= \int dA_zdZ^{*}dZ  e^{-L} W_C(z,z')B_C(z)
\end{array}
\label{int}
\end{equation}
where $W_C^{\prime}(z,z')$ denotes the Wilson line based on
$A_1^{\prime}$ and $A_2^{\prime}$ and similarly $B_C^{\prime}(z)$
denotes  the field obtained from $B_C(z)$ with  $(A_1, A_2,Z)$
replaced by $(A_1^{\prime}, A_2^{\prime},Z^{\prime})$.


Then it can be shown that the differential is unchanged under a
gauge transformation \cite{Fad}:
\begin{equation}
dA_z^{\prime}dZ^{\prime *}dZ^{\prime}
= dA_zdZ^{*}dZ
\label{int2}
\end{equation}
Also by the gauge invariance property the factor $e^{-L}$ is
unchanged under a gauge transformation. Thus from (\ref{int}) we
have
\begin{equation}
0 = \langle W^{\prime}(z,z')B^{\prime}(z)\rangle -
  \langle W(z,z')B(z)\rangle
\label{w1}
\end{equation}
where the correlation notation $\langle \cdot\rangle$ denotes the integral with respect to the differential (or the measure of the gauge model):
 \begin{equation} e^{-L}dA_zdZ^{*}dZ \label{w1a} \end{equation}

We can now carry out the calculus of variation. From the gauge
transformation we have the formula:
\begin{equation}
W_C^{\prime}(z,z')=U(a(z))W_C(z,z')U^{-1}(a(z'))\,, \label{aw1a}
\end{equation}
where $a(z)=\mbox{Re}\,\omega(z)$. This
gauge transformation gives a variation of $W_C(z,z')$ with the
gauge function $a(z)$
as the variational direction $a$
in the variational formulas (\ref{km5}) and  (\ref{km6}). Thus analogous
to the variational formula (\ref{km6}) we have that the variation
of $W(z,z')$ under this gauge transformation is given by
\begin{equation}
W(z,z')
  \frac{1}{2\pi i}\oint_z dw a(w)J(w)
\label{int3}
\end{equation}
where the generator $J$ for this variation is to
be specific. This $J$ will be a quantum generator
which generalizes the classical generator $J$ in
(\ref{km6}).
Thus under a gauge transformation with gauge function $a(z)$ from (\ref{w1}) we have the
following variational equation:
\begin{equation}
0= \langle W(z,z')[\delta_{a}B(z)+\frac{1}{2\pi i}\oint_z
dw a(w)J(w)B(z)]\rangle
\label{w2}
\end{equation}
where $\delta_{a}B(z) $
denotes the variation of the field
$B(z)$ in the direction $a(z)$.
From this equation an ansatz of
$J$ is that $J$ satisfies the following equation:
\begin{equation}
W(z,z')[\delta_{a}B(z)+\frac{1}{2\pi i}\oint_z
dw a(w)J(w)B(z)] =0 \label{n8bb}
\end{equation}
From this equation we have the following variational equation:
\begin{equation}
\delta_{a}B(z)=\frac{-1}{2\pi i}\oint_z dw a(w)J(w)B(z)
\label{n8bre}
\end{equation}
where $a(w)J(w)$ is defined by  (\ref{km5a}). 
This completes the 
calculus of variation. 

From (\ref{n8bre}) we have the following  
equation for defining the generator $J$:
\begin{equation}
\delta_{\omega}B(z)=\frac{-1}{2\pi i}\oint_z dw\omega(w)J(w)B(z)
\label{n8b}
\end{equation}
where we generalize the direction $a(z)=\mbox{Re}\,\omega(z)$ to
the analytic variation direction $\omega(z)$ in (\ref{variation}).

Let us 
then add one more condition for determining the quantum generator
$J$ in (\ref{n8b}). 
As analogous to the WZW model in conformal field
theory \cite{Fra}\cite{Fuc} \cite{Kni}  let us consider a classical generator $J_C$
given by
\begin{equation}
J_C(z) := -k_0 W_C^{-1}(z, z')\partial_z W_C(z, z') \label{n6}
\end{equation}
where we define $\partial_z=\partial_{ x^1} +i\partial_{ x^2} $ and we set $z'=z$ after
the differentiation with respect to $z$; $ k_0>0 $ is a constant which is fixed when
the $J_C$ is determined to be of the form (\ref{n6}) and the minus sign is chosen by
convention. In the WZW model \cite{Fra}\cite{Kni}
 the $J_C$ of the form (\ref{n6})
is the  generator  of the chiral symmetry of the WZW model. We can
write the $J_C$ in (\ref{n6}) in the following form:
\begin{equation}
 J_C(w) = \sum_a J_C^a(w) =
\sum_a j_C^a(w)e_0 t^a
\label{km82}
\end{equation}
Thus the quantum version $J$ of $ J_C$ is of the following form:
\begin{equation}
 J(w) = \sum_a J^a(w) =
\sum_a j^a(w)e_0 t^a
\label{km8}
\end{equation}
where  $ j^a(w)$ are quantum generators.

We see that the generators $t^a$ of $SU(2)$ appear in this form of
$J_C$ and $J$ and  this form is analogous to the classical generator $J$ in
(\ref{km3}). This shows that
 this $J$ is a possible candidate for the quantum generator
$J$ in (\ref{n8b}).


Then from (\ref{n5}) and (\ref{n6}) we
have that the variation $\delta_{\omega}J_C$ of the generator $J_C$ in
(\ref{n6}) is given by \cite{Fra}(p.622) \cite{Kni}:
\begin{equation}
\delta_{\omega}J_C= \lbrack J_C, \omega\rbrack -k_0\partial_z \omega
\label{n8c2}
\end{equation}

Then  the quantum generator $J$ is required also to have this property  (\ref{n8c2}) of  the classical generator $J_C$:
\begin{equation}
\delta_{\omega}J= \lbrack J, \omega\rbrack -k_0\partial_z \omega
\label{n8c}
\end{equation}




Now we 
show that this quantum  generator $J$ in (\ref{n8b})   and  (\ref{n8c}) can be uniquely solved.
From (\ref{n8b}) and (\ref{n8c}) we have that $J$ satisfies the
following relation of current algebra
\cite{Fra}\cite{Fuc}\cite{Kni}:
\begin{equation}
J^a(w)J^b(z)=\frac{k_0\delta_{ab}}{(w-z)^2}
+\sum_{c}if_{abc}\frac{J^c(z)}{(w-z)} \label{n8d}
\end{equation}
where as a convention the regular term of the product
$J^a(w)J^b(z)$ is omitted. Then by following
\cite{Fra}\cite{Fuc}\cite{Kni} from (\ref{n8d}) and (\ref{km8})
we can show that the $J_n^a$ in (\ref{km3})  for the corresponding Laurent series of
the quantum generator $J$
satisfy the
following  Kac-Moody algebra:
\begin{equation}
[J_m^a, J_n^b] = if_{abc}J_{m+n}^c + k_0 m\delta_{ab}\delta_{m+n, 0}
\label{n8}
\end{equation}

We remark that the constant $k_0$ is  called the central extension or the level of
the Kac-Moody algebra.

{\bf Remark}. Let us also consider the other side of the above chiral symmetry. Similar to the $J$ in (\ref{n6}) we define a generator $J_C^{\prime}$ by:
\begin{equation}
J_C^{\prime}(z')= k_0 \partial_{z'}W_C(z, z')W_C^{-1}(z, z') \label{d1}
\end{equation}
where after differentiation with respect to $z'$ we set $z=z'$. Let us then consider the following formal correlation:
\begin{equation}
\langle B(z')W(z,z') \rangle := \int dAdZ^{*}dZ  B(z')W_C(z,z')e^{-L}
\label{n8aa}
\end{equation}
where $z$ is fixed. By an approach similar to the above derivation of (\ref{n8b}) we have the following  variational equation:
\begin{equation}
\delta_{\omega^{\prime}}B(z') =\frac{-1}{2\pi i}\oint_{z^{\prime}}
dwB(z')J^{\prime}(w) \omega^{\prime}(w) \label{n8b1}
\end{equation}
Then similar to (\ref{n8c}) we also have
\begin{equation}
\delta_{\omega^{\prime}}J^{\prime}= \lbrack  J^{\prime},
\omega^{\prime}\rbrack -k_0 \partial_{z'} \omega^{\prime}
\label{n8c1}
\end{equation}
Then from (\ref{n8b1}) and (\ref{n8c1}) we can derive the current algebra and the Kac-Moody algebra for $J^{\prime}$ which are of the same form of (\ref{n8d}) and (\ref{n8}). From this we  have $J^{\prime}=J$. $\diamond$

\section{Quantum Knizhnik-Zamolodchikov Equation In Dual Form} \label{sec7}

With the above current algebra $J$ and the formula (\ref{n8b}) we can now
follow the usual approach
in conformal field theory to derive a quantum
Knizhnik-Zamolodchikov (KZ) equation for the product of
primary fields in a conformal field theory \cite{Fra}\cite{Fuc}\cite{Kni}.
We shall derive the KZ equation for the product of $n$ Wilson lines $W(z, z')$.
Here an important point is that from the two sides of
$W(z, z')$  we can derive two quantum KZ equations which are
dual to each other. These two quantum KZ equations are different from the usual KZ equation in that they are equations for the quantum operators $W(z, z')$ while the usual KZ equation is for the correlations of quantum operators.

With this difference the following derivation of  KZ equation for deriving these two quantum KZ equations is well known  in conformal field theory  \cite{Fra}\cite{Fuc}. The reader may skip this derivation of  KZ equation and just look at the form of the Virasoro energy operator $T(z)$ and the Virasoro algebra and the form of these two quantum KZ equations.

 Let us first consider the following quantum version of (\ref{k1c}) with $W_C(z, z')$ replaced by the quantum $W(z, z')$:
\begin{equation}
\delta_{\omega}W(z,z')=W(z,z')\omega(z)e_0, \quad 
\delta_{\omega^{\prime}}W(z,z')=-e_0\omega^{\prime}(z')W(z,z')
\label{k1}
\end{equation}
From (\ref{n8b}) and (\ref{k1}) 
we have:
\begin{equation}
J^a(z)W(w, w') = \frac{-e_0 t^aW(w,w')}{z-w}
\label{k3}
\end{equation}
where as a convention the regular term of the product $J^a(z)W(w, w')$
is omitted.

Following \cite{Fra} and \cite{Fuc}
let us define a Virasoro energy operator (also called an energy-momentum tensor) $T(z)$ by:
\begin{equation}
T(z) := \frac{1}{2(k_0+g_0)}\sum_a :J^a(z)J^a(z):
\label{k4}
\end{equation}
where $g_0$ is the dual Coxeter number of the gauge group \cite{Fra}. In (\ref{k4})
the symbol $:J^a(z)J^a(z):$ denotes the normal ordering  of the operator $J^a(z)J^a(z)$ which can be
defined as follows \cite{Fra}\cite{Fuc}. Let a product of operators
$A(z)B(w)$ be written in the following Laurent series
form:
\begin{equation}
A(z)B(w)= \sum_{n=-n_0}^{\infty}
 a_n(w)(z-w)^n
\label{normal}
\end{equation}
The singular part
of (\ref{normal}) is called the contraction
of $A(z)B(w)$ and will be denoted by $\overbrack{A(z)B}(w)$.
Then the term $a_0(w)$ is called the normal ordering
of $A(z)B(w)$ and we denote $a_0(w)$ by $:A(w)B(w):$.

The above definition of the energy operator $T(z)$ is called the Sugawara construction \cite{Fra}. We first have the following well known theorem on $T(z)$ in conformal field theory \cite{Fra}:
\begin{theorem}

The operator product $T(z)T(w)$ is given by the following
formula:
\begin{equation}
T(z)T(w)=\frac{c}{2(z-w)^4}+
         \frac{2T(w)}{(z-w)^2}+\frac{\partial T(w)}{(z-w)}
\label{k5}
\end{equation}
for some constant $c=\frac{4k_0}{k_0+g_0}$  
and as a convention we omit the regular term of this product.
\end{theorem}
{\bf Proof}. In \cite{Fra} there is a detail proof of this theorem where the constant $ \frac{1}{2(k_0+g_0)}$ is shown to be a normalization constant such that this theorem holds. Here we want to remark that the formula (\ref{n8d}) for the product $J^a(z)J^b(x)$ is used for the proof of this theorem. $\diamond$

From this theorem we then have the following Virasoro algebra
of the mode expansion of $T(z)$ \cite{Fra}\cite{Fuc}:
\begin{theorem}

Let us write $T(z)$ in the following Laurent series form:
\begin{equation}
T(z)= \sum_{n=-\infty}^{\infty}(z-w)^{-n-2}L_n(w)
\label{conform9}
\end{equation}
This means that the modes $L_n(w)$ are defined by
\begin{equation}
L_n(w):= \frac{1}{2\pi i}\oint_w dz (z-w)^{n+1}T(z)
\label{conform10}
\end{equation}
Then we have that $L_n$ form a Virasoro algebra:
 \begin{equation}
[L_n, L_m] = (n-m)L_{n+m} + \frac{c}{12}n(n^2-1)\delta_{n+m,0}
\label{conform11}
\end{equation}
\end{theorem}

From  the formula (\ref{n8d}) for the product $J^a(z)J^b(w)$ we have the following operator product expansion \cite{Fra}:
\begin{equation}
\overbrack{T(z) J^a}(w)
= \frac{J^a(w)}{(z-w)^2}+ \frac{\partial J^a(w)}{(z-w)}
\label{conform7}
\end{equation}

Then we have the following operator product of $T(z)$ with an operator $A(w)$:
 \begin{equation}
T(z)A(w) = \sum_{n=-\infty}^{\infty}(z-w)^{-n-2}L_nA(w)
\label{conform13}
\end{equation}

From (\ref{conform7}) and (\ref{conform13}) we have that $L_{-1}J^a(w)= \partial J^a(w)$ and $L_{-1} =  \frac{\partial}{\partial z}$.
Thus we have
\begin{equation}
L_{-1}W(w, w')= \frac{\partial W(w, w')}{\partial w}
\label{conform14}
\end{equation}

On the other hand as shown in \cite{Fra},  by using the Laurent series expansion of
$J^a(z)$ in the section on Kac-Moody algebra we can compute
the normal ordering $:J^a(z)J^a(z):$ from which we have the
Laurent series expansion of $T(z)$ with $L_{-1}$ given by \cite{Fra}:
\begin{equation}
L_{-1}=\frac{1}{2(k_0+g_0)}\sum_a \lbrack
\sum_{m\leq -1}J_m^aJ_{-1-m}^a +
\sum_{m\geq 0}J_{-1-m}^a J_m^a \rbrack
\label{conform15}
\end{equation}
where since $J_m^a$ and $J_{-1-m}^a$ commute each other the ordering of them
is irrelevant.

From (\ref{conform15}) we then have
\begin{equation}
\begin{array}{rl}
 & L_{-1}W(w,w') \\
 =& \frac{1}{2(k_0+g_0)}\sum_a \lbrack
\sum_{m\leq -1}J_m^a(w)J_{-1-m}^a(w) +
\sum_{m\geq 0}J_{-1-m}^a(w) J_m^a(w) \rbrack W(w,w') \\
 = & \frac{1}{(k_0+g_0)}J_{-1}^a(w) J_{0}^a(w) W(w,w')
\end{array}
\label{conform16}
\end{equation}
since $J_m^a W(w,w')=0$ for $m>0$ by  (\ref{k3}) and the quantum version of $J_m^a$ in (\ref{km4})  (where the $\sum_a$ is omitted).

It follows from (\ref{conform14}) and (\ref{conform16}) that we have:
the following equality:
\begin{equation}
\partial_w W(w, w')  = \frac{1}{(k_0+g_0)}J_{-1}^a(w) J_{0}^a(w) W(w,w')
\label{k7}
\end{equation}
Then from (\ref{k3}) and the quantum version of $J_0^a$ in (\ref{km4}) we have:
\begin{equation}
J_{0}^a(w)W(w,w')=-e_0 t^aW(w,w')
\label{k8}
\end{equation}
From (\ref{k7}) and (\ref{k8}) we then have
\begin{equation}
\partial_z W(z, z')=\frac{-e_0}{k_0+g_0}J_{-1}^a(z)t^aW(z,z')
\label{k9}
\end{equation}

Now let us
consider a product of $n$ Wilson lines:
$ W(z_1, z_1^{\prime})\cdot\cdot\cdot
 W(z_n, z_n^{\prime})$.
Let this product be represented as a tensor product when
$z_i$ and $z_j^{\prime}$, $i,j=1,...,n$ are all independent variables.
 Then from (\ref{k9}) we have
\begin{equation}
\begin{array}{rl}
& \partial_{z_{i}} W(z_1, z_1^{\prime})\cdot\cdot\cdot
W(z_i, z_i^{\prime})\cdot\cdot\cdot
 W(z_n, z_n^{\prime})\\
=& \frac{-e_0}{k_0+g_0}W(z_1, z_1^{\prime})\cdot\cdot\cdot
J_{-1}^a(z_i)t^aW(z_i,z_i^{\prime})\cdot\cdot\cdot
W(z_n, z_n^{\prime})\\
=& \frac{-e_0}{k_0+g_0}J_{-1}^a(z_i)t^a W(z_1, z_1^{\prime})\cdot\cdot\cdot
W(z_i,z_i^{\prime})\cdot\cdot\cdot
W(z_n, z_n^{\prime})
\end{array}
\label{k9a}
\end{equation}
where the second equality is from the definition
of tensor product for which we define
\begin{equation}
t^a W(z_1, z_1^{\prime})\cdot\cdot\cdot
W(z_i,z_i^{\prime})\cdot\cdot\cdot
W(z_n, z_n^{\prime})
:=W(z_1, z_1^{\prime})\cdot\cdot\cdot
[t^a W(z_i,z_i^{\prime})]\cdot\cdot\cdot
W(z_n, z_n^{\prime})
\label{k9aa}
\end{equation}

With this formula (\ref{k9a}) we can now follow \cite{Fra} and \cite{Fuc} to
derive the KZ equation.
For a easy reference let us present this derivation in \cite{Fra} and \cite{Fuc} as follows.
From the Laurent series of $J^a$ we have
\begin{equation}
J_{-1}^a(z_i) = \frac{1}{2\pi i}
\oint_{z_i} \frac{dz}{z-z_i}J^a(z)
\label{norm2}
\end{equation}
where the line integral is on a contour encircling  $z_i$.  We let this contour be a circle with center $z_i$ and with small radious such that all other $z_j $ of $W(z_j, z_j^{\prime})$ for $j=1,...,n$, $ j\neq i$,  are excluded.
Then we have:
\begin{equation}
\begin{array}{rl}
 & J_{-1}^a(z_i)W(z_1, z_1^{\prime})\cdot\cdot\cdot
 W(z_n, z_n^{\prime}) \\
& \\
= & \frac{1}{2\pi i}\oint_{z_i} \frac{dz}{z-z_i}
 J^a(z)W(z_1, z_1^{\prime})\cdot\cdot\cdot
W(z_n, z_n^{\prime})
 \\
& \\
= & \frac{1}{2\pi i}\oint_{z_i} \frac{dz}{z-z_i}
 \sum_{j=1}^n W(z_1, z_1^{\prime})\cdot\cdot\cdot
 [\frac{-e_0 t^a}{z-z_j}W(z_j,z_j^{\prime})]\cdot\cdot\cdot
W(z_n, z_n^{\prime})
\end{array}
\label{norm3aa}
\end{equation}
where the second equality is from the $JW$ product formula (\ref{k3}). Then we have that 
 (\ref{norm3aa}) is equal to:.
\begin{equation}
\begin{array}{rl}
& \sum_{j=1,j\neq i}^n
\frac{-e_0 t_j^a}{z_i-z_j}W(z_1, z_1^{\prime})\cdot\cdot\cdot
 W(z_n, z_n^{\prime})
 \end{array}
\label{norm3}
\end{equation}
where 
we have used the definition
of tensor product.

From (\ref{norm3})
 and by applying (\ref{k9a})
to $z_i$ for $i=1,...,n$ we have the following
Knizhnik-Zamolodchikov equation \cite{Fra}\cite{Fuc}\cite{Kni}:
\begin{equation}
\partial_{z_i}
 W(z_1, z_1^{\prime})\cdot\cdot\cdot
W(z_n, z_n^{\prime})
=\frac{e_0^2}{k_0+g_0}
\sum_{j\neq i}^{n}
\frac{\sum_a t_i^a \otimes t_j^a}{z_i-z_j}
 W(z_1, z_1^{\prime})\cdot\cdot\cdot
W(z_n, z_n^{\prime})
\label{n9}
\end{equation}
for $i=1, ..., n$. We remark that in (\ref{n9}) we have
defined $t_i^a:= t^a$ and
\begin{equation}
\begin{array}{rl}
 & t_i^a \otimes t_j^a W(z_1, z_1^{\prime})\cdot\cdot\cdot
W(z_n, z_n^{\prime}) \\
:=& W(z_1, z_1^{\prime})\cdot\cdot\cdot
 [t^aW(z_i, z_i^{\prime})]\cdot\cdot\cdot
[t^aW(z_j, z_j^{\prime})]\cdot\cdot\cdot
W(z_n, z_n^{\prime})
\end{array}
\label{n9a}
\end{equation}

 It is interesting and important that we also have
another KZ equation with respect to the
$z_i^{\prime}$ variables. The derivation of this  KZ equation is
dual to the above derivation in that the operator products and
their corresponding variables are with reverse order to that in
the above derivation.

From (\ref{n8b1}) and (\ref{k1}) we have a $WJ^{\prime}$ operator product given
by
\begin{equation}
W(w, w')J^{\prime a}(z') = \frac{-W(w, w')e_0 t^a}{w'-z'}
\label{d2}
\end{equation}
where we have omitted the regular term of the product.
Then similar to the above derivation of the KZ equation from (\ref{d2})
we can then derive the following Knizhnik-Zamolodchikov
equation which is dual to (\ref{n9}):
\begin{equation}
\partial_{z_i^{\prime}}
 W(z_1,z_1^{\prime})\cdot\cdot\cdot W(z_n,z_n^{\prime})
= \frac{e_0^2}{k_0+g_0}\sum_{j\neq i}^{n}
 W(z_1, z_1^{\prime})\cdot\cdot\cdot
W(z_n, z_n^{\prime})
\frac{\sum_a t_i^a\otimes t_j^a}{z_j^{\prime}-z_i^{\prime}}
\label{d8}
\end{equation}
for $i=1, ..., n$ where we have defined:
\begin{equation}
\begin{array}{rl}
 &  W(z_1, z_1^{\prime})\cdot\cdot\cdot
W(z_n, z_n^{\prime})t_i^a \otimes t_j^a \\
:=& W(z_1, z_1^{\prime})\cdot\cdot\cdot
 [W(z_i, z_i^{\prime})t^a]\cdot\cdot\cdot
[W(z_j, z_j^{\prime})t^a]\cdot\cdot\cdot
W(z_n, z_n^{\prime})
\end{array}
\label{d8a}
\end{equation}

{\bf Remark}. From the generator $J$ and the Kac-Moody algebra
we have derived a quantum KZ equation in dual form.
This quantum KZ equation in dual form may be considered as a quantum Yang-Mills equation since it is analogous to the classical Yang-Mills equation which is derived from the classical Yang-Mills gauge model. This quantum KZ equation in dual form will be as the starting point for the construction of quantum knots. $\diamond$

\section{Solving Quantum KZ Equation In Dual Form}\label{sec8a}

Let us consider the following product of two quantum
Wilson lines:
\begin{equation}
G(z_1,z_2, z_3, z_4):=
 W(z_1, z_2)W(z_3, z_4)
\label{m1}
\end{equation}
where the two quantum Wilson lines $W(z_1, z_2)$ and
$W(z_3, z_4)$ represent two pieces
of curves starting at $z_1$ and $z_3$ and ending at
$z_2$ and $z_4$ respectively.

We have that this product $G$ satisfies the KZ equation for the
variables $z_1$, $z_3$ and satisfies the dual KZ equation
for the variables $z_2$ and $z_4$.
Then
by solving the two-variables-KZ equation in (\ref{n9}) we have that a form of $G$ is
given by \cite{Koh}\cite{Dri}\cite{Chari}:
\begin{equation}
e^{-\hat{t}\log [\pm (z_1-z_3)]}C_1
\label{m2}
\end{equation}
where $\hat{t}:=\frac{-e_0^2}{k_0+g_0}\sum_a t^a \otimes t^a$
and $C_1$ denotes a constant matrix which is independent
of the variable $z_1-z_3$.

We see that $G$ is a multivalued analytic function
where the determination of the $\pm$ sign depended on the choice of the
branch.

Similarly by solving the dual two-variable-KZ equation
 in (\ref{d8}) we have that
$G$ is of the form
\begin{equation}
C_2e^{\hat{t}\log [\pm (z_4-z_2)]}
\label{m3}
\end{equation}
where $C_2$ denotes a constant matrix which is independent
of the variable $z_4-z_2$.

From (\ref{m2}), (\ref{m3}) and we let
$C_1={\bf A}e^{\hat{t}\log[\pm (z_4-z_2)]}$,
$C_2= e^{-\hat{t}\log[\pm (z_1-z_3)]}{\bf A}$ where ${\bf A}$ is a constant matrix  we have that
$G$ is given by
\begin{equation}
G(z_1, z_2, z_3, z_4)=
e^{-\hat{t}\log [\pm (z_1-z_3)]}{\bf A}e^{\hat{t}\log [\pm (z_4-z_2)]}
\label{m4}
\end{equation}
where at the singular case that $z_1=z_3$ we simply define $\log [\pm (z_1-z_3)]=0$. Similarly for $z_2=z_4$.

Let us find a form of the initial operator ${\bf A}$. We notice that there are two operators $\Phi_{\pm}(z_1-z_3):=e^{-\hat{t}\log [\pm (z_1-z_3)]}$
and  $\Psi_{\pm}(z_4-z_2):=e^{-\hat{t}\log [\pm (z_4-z_2)]}$ 
acting on the two sides of ${\bf A}$ respectively where  the two independent variables  $z_1, z_3$ of $\Phi_{\pm}$ are mixedly from the two quantum Wilson lines $W(z_1, z_2)$ and
$W(z_3, z_4)$ respectively and the two independent variables  $z_2, z_4$ of $\Psi_{\pm}$ are mixedly from the two quantum Wilson lines $W(z_1, z_2)$ and $W(z_3, z_4)$ respectively. From this we determine the form of $A$ as follows.

Let $D$ denote a representation of $SU(2)$. Let $D(g)$ represent an element $g$ of  $SU(2)$
and let $D(g)\otimes D(g)$ denote the tensor product representation of $SU(2)$. Then
in the KZ equation we define
\begin{equation}
[t^a\otimes t^a] [D(g_1)\otimes D(g_1)]\otimes
[D(g_2)\otimes D(g_2)]
:=[t^aD(g_1)\otimes D(g_1)]\otimes
[t^aD(g_2)\otimes D(g_2)]
\label{tensorproduct}
\end{equation}
and
\begin{equation}
[D(g_1)\otimes D(g_1)]\otimes
[D(g_2)\otimes D(g_2)][t^a\otimes t^a]
:=[D(g_1)\otimes D(g_1)t^a]\otimes
[D(g_2)\otimes D(g_2)t^a]
\label{tensorproduct2}
\end{equation}

Then we let $U({\bf a})$
denote the universal
enveloping algebra
where ${\bf a}$ denotes an algebra which is formed by the Lie
algebra $su(2)$ and the identity matrix.

Now let the initial operator ${\bf A}$ be of the form ${\bf A}_1\otimes {\bf A}_2\otimes
{\bf A}_3\otimes{\bf A}_4$ where ${\bf A}_i,i=1,...,4$ denote operators
taking values in $U({\bf a})$.
In
this case we have that in (\ref{m4}) the operator
$\Phi_{\pm}(z_1-z_3):=e^{-\hat{t}\log [\pm (z_1-z_3)]}$ acts on ${\bf A}$ from
the left via the following formula:
\begin{equation}
t^a\otimes t^a{\bf A}=
[t^a {\bf A}_1]\otimes{\bf A}_2\otimes [t^a{\bf A}_3]\otimes{\bf A}_4
\label{ini2}
\end{equation}

Similarly the operator
$\Psi_{\pm}(z_4-z_2):=e^{\hat{t}\log [\pm (z_4-z_2)]}$
in (\ref{m4}) acts on ${\bf A}$ from the right via the following formula:
\begin{equation}
{\bf A} t^a\otimes t^a =
{\bf A}_1\otimes [{\bf A}_2 t^a]\otimes{\bf A}_3\otimes[{\bf A}_4 t^a]
\label{ini3}
\end{equation}

We may generalize the above tensor product of
two quantum Wilson lines as follows.
Let us consider a tensor product of $n$ quantum Wilson lines:
$W(z_1, z_1^{\prime})\cdot\cdot\cdot W(z_n, z_n^{\prime})$
where the variables $z_i$, $z_i^{\prime}$
are all independent. By solving the two KZ equations
we have that this tensor product is given by:
\begin{equation}
W(z_1, z_1^{\prime})\cdot\cdot\cdot W(z_n, z_n^{\prime})
=\prod_{ij} \Phi_{\pm}(z_i-z_j)
{\bf A}\prod_{ij}
\Psi_{\pm}(z_i^{\prime}-z_j^{\prime})
\label{tensor}
\end{equation}
where $\prod_{ij}$ denotes a product of
$\Phi_{\pm}(z_i-z_j)$ or
$\Psi_{\pm}(z_i^{\prime}-z_j^{\prime})$
for $i,j=1,...,n$ where $i\neq j$.
In (\ref{tensor}) the initial operator
${\bf A}$ is represented as a tensor product of operators ${\bf A}_{iji^{\prime}j^{\prime}}, i,j,i^{\prime}, j^{\prime}=1,...,n$ where each ${\bf A}_{iji^{\prime}j^{\prime}}$ is of the form of the initial operator ${\bf A}$ in the above tensor product of two-Wilson-lines case and is acted by $\Phi_{\pm}(z_i-z_j)$ or
$\Psi_{\pm}(z_i^{\prime}-z_j^{\prime})$ on its two sides respectively.

\section{Computation of Quantum Wilson Lines}\label{sec 8aa}

Let us consider the following product of two quantum
Wilson lines:
\begin{equation}
G(z_1,z_2, z_3, z_4):=
W(z_1, z_2)W(z_3, z_4)
\label{h1}
\end{equation}
where the two quantum Wilson lines $W(z_1, z_2)$ and
$W(z_3, z_4)$ represent two pieces
of curves starting at $z_1$ and $z_3$ and ending at
$z_2$ and $z_4$ respectively.
As shown in the above section we have that $G$
is given by the following formula:
\begin{equation}
G(z_1, z_2, z_3, z_4)=
e^{-\hat{t}\log [\pm (z_1-z_3)]}{\bf A}e^{\hat{t}\log [\pm (z_4-z_2)]}
\label{m4a}
\end{equation}
where the product is
a 4-tensor.

Let us set $z_2=z_3$. Then
the 4-tensor $W(z_1, z_2)W(z_3, z_4)$ is reduced to the 2-tensor
$W(z_1, z_2)W(z_2, z_4)$.
By using (\ref{m4a}) the 2-tensor
$W(z_1, z_2)W(z_2, z_4)$
is given by:
\begin{equation}
W(z_1, z_2)W(z_2, z_4)
=e^{-\hat{t}\log [\pm (z_1-z_2)]}{\bf A}_{14}e^{\hat{t}\log [\pm (z_4-z_2)]}
\label{closed1}
\end{equation}
where ${\bf A}_{14}={\bf A}_1\otimes{\bf A}_4$ is a 2-tensor reduced from the 4-tensor
${\bf A}={\bf A}_1\otimes{\bf A}_2\otimes{\bf A}_3\otimes{\bf A}_4$ in (\ref{m4a}). In this reduction the $\hat{t}$ operator of $\Phi=e^{-\hat{t}\log [\pm (z_1-z_2)]}$ acting on the left side of ${\bf A}_1$ and ${\bf A}_3$ in ${\bf A}$ is reduced to acting on the left side of ${\bf A}_1$ and ${\bf A}_4$ in ${\bf A}_{14}$. Similarly  the $\hat{t}$ operator of $\Psi=e^{-\hat{t}\log [\pm (z_4-z_2)]}$ acting on the right side of ${\bf A}_2$ and ${\bf A}_4$ in ${\bf A}$ is reduced to acting on the right side of ${\bf A}_1$ and ${\bf A}_4$ in ${\bf A}_{14}$.

Then since $\hat{t}$ is a 2-tensor operator we have that $\hat{t}$ is as a matrix acting on the two sides of the 2-tensor ${\bf A}_{14}$ which is also as a matrix with the same dimension as $\hat{t}$.
Thus $\Phi$ and $\Psi$ are as matrices of the same dimension as the matrix
${\bf A}_{14}$  acting on ${\bf A}_{14}$ by the usual matrix operation.
Then since $\hat{t}$ is a Casimir operator for the 2-tensor group representation of $SU(2)$ we have that
$\Phi $  and $\Psi $ commute
with ${\bf A}_{14}$ since  $\Phi $  and $\Psi$ are exponentials
of $\hat{t}$ (We remark that $\Phi $  and $\Psi $ are in general not commute with the 4-tensor initial operator ${\bf A}$).
Thus we have
\begin{equation}
e^{-\hat{t}\log [\pm (z_1-z_2)]}{\bf A}_{14}e^{\hat{t}\log[\pm (z_4-z_2)]}
=e^{-\hat{t}\log [\pm (z_1-z_2)]}e^{\hat{t}\log[\pm (z_4-z_2)]}{\bf A}_{14}
\label{closed1a}
\end{equation}

We let $W(z_1, z_2)W(z_2, z_4)$ be as a representation of the quantum Wilson line $W(z_1,z_4)$ and we write $W(z_1,z_4)=W(z_1, z_2)W(z_2, z_4)$.
Then we have the following representation of
$W(z_1,z_4)$:
\begin{equation}
W(z_1,z_4)=W(z_1,w_1)W(w_1,z_4)=e^{-\hat{t}\log [\pm (z_1-w_1)]}e^{\hat{t}\log[\pm (z_4-w_1)]}{\bf A}_{14}
\label{closed1a1}
\end{equation}
This representation of the quantum Wilson line $W(z_1,z_4)$ means that the line (or path) with end points $z_1$ and $z_4$ is specified that it passes the intermediate point $w_1=z_2$. This representation shows the quantum nature that the path is  not specified at other intermediate points except the intermediate point $w_1=z_2$. This unspecification of the path is of the same quantum nature of the Feymann path description of quantum mechanics.

Then let us consider another representation of the quantum Wilson line $W(z_1,z_4)$. We consider $W(z_1,w_1)W(w_1,w_2)W(w_2,z_4)$ which is obtained from the tensor $W(z_1,w_1)W(u_1,w_2)W(u_2,z_4)$ by two reductions where $z_j$, $w_j$, $u_j$, $j=1,2$ are independent variables. For this representation we have:
\begin{equation}
W(z_1,w_1)W(w_1,w_2)W(w_2,z_4)
=e^{-\hat{t}\log [\pm (z_1-w_1)]}e^{-\hat{t}\log [\pm (z_1-w_2)]}
e^{\hat{t}\log[\pm (z_4-w_1)]}e^{\hat{t}\log[\pm (z_4-w_2)]}{\bf A}_{14}
\label{closed1a2}
\end{equation}
This representation of the quantum Wilson line $W(z_1,z_4)$ means that the line (or path) with end points $z_1$ and $z_4$ is specified that it passes the intermediate points $w_1$ and $w_2$. This representation shows the quantum nature that the path is  not specified at other intermediate points except the intermediate points $w_1$ and $w_2$. This unspecification of the path is of the same quantum nature of the Feymann path description of quantum mechanics.

Similarly we may represent the quantum Wilson line $W(z_1,z_4)$ by path with end points $z_1$ and $z_4$ and is specified only to pass at finitely many intermediate points. Then we let the quantum Wilson line $W(z_1,z_4)$ as an equivalent class of all these representations. Thus we may write $W(z_1,z_4)=W(z_1,w_1)W(w_1,z_4)=W(z_1,w_1)W(w_1,w_2)W(w_2,z_4)=\cdot\cdot\cdot$.

{\bf Remark}. Since ${\bf A}_{14}$ is a 2-tensor
we have that a natural group representation for the Wilson line $W(z_1,z_4)$ is the 2-tensor group representation of the group $SU(2)$.

\section{Representing Braiding of Curves by Quantum Wilson Lines}\label{sec 9aa}

Consider again the product $G(z_1, z_2, z_3, z_4)=W(z_1,z_2)W(z_3,z_4)$.
We have that $G$ is a multivalued analytic function
where the determination of the $\pm$ sign depended on the choice of the
branch.

Let the two pieces of curves be crossing at $w$. Then we have $W(z_1,z_2)=W(z_1,w)W(w,z_2)$ and
 $W(z_3,z_4)=W(z_3,w)W(w,z_4)$. Thus we have
\begin{equation}
W(z_1,z_2)W(z_3,z_4)=
W(z_1,w)W(w,z_2)W(z_3,w)W(w,z_4)
\label{h2}
\end{equation}

If we interchange $z_1$ and $z_3$, then from
(\ref{h2}) we have the following ordering:
\begin{equation}
 W(z_3,w)W(w, z_2)W(z_1,w)W(w,z_4)
\label{h3}
\end{equation}

Now let us choose a  branch. Suppose that
these two curves are cut from a knot and that
following the orientation of a knot the
curve represented by  $W(z_1,z_2)$ is before the
curve represented by  $W(z_3,z_4)$. Then we fix a branch such that the  product in (\ref{m4a}) is
with two positive signs :
\begin{equation}
W(z_1,z_2)W(z_3,z_4)=
e^{-\hat{t}\log(z_1-z_3)}{\bf A}e^{\hat{t}\log(z_4-z_2)}
\label{h4}
\end{equation}

Then if we interchange $z_1$ and $z_3$ we have
\begin{equation}
W(z_3,w)W(w, z_2)W(z_1,w)W(w,z_4) =
e^{-\hat{t}\log[-(z_1-z_3)]}{\bf A}e^{\hat{t}\log(z_4-z_2)}
\label{h5}
\end{equation}
From (\ref{h4}) and (\ref{h5}) as a choice of branch we have
\begin{equation}
W(z_3,w)W(w, z_2)W(z_1,w)W(w,z_4) =
R W(z_1,w)W(w,z_2)W(z_3,w)W(w,z_4)
\label{m7a}
\end{equation}
where $R=e^{-i\pi \hat{t}}$ is the monodromy of the KZ equation.
In (\ref{m7a}) $z_1$ and $z_3$ denote two points on a closed curve
such that along the direction of the curve the point
$z_1$ is before the point $z_3$ and in this case we choose
a branch such that the angle of $z_3-z_1$ minus the angle
of $z_1-z_3$ is equal to $\pi$.

{\bf Remark}. We may use other representations of $W(z_1,z_2)W(z_3,z_4)$. For example we may use the following representation:
\begin{equation}
\begin{array}{rl}
 &W(z_1,w)W(w, z_2)W(z_3,w)W(w,z_4)\\
= &e^{-\hat{t}\log(z_1-z_3)}e^{-2\hat{t}\log(z_1-w)}e^{-\hat{t}2\log(z_3-w)}{\bf A}e^{\hat{t}\log(z_4-z_2)}e^{2\hat{t}\log(z_4-w)}e^{2\hat{t}\log(z_2-w)}
\end{array}
\label{h4a}
\end{equation}
Then the interchange of $z_1$ and $z_3$ changes only $z_1-z_3$ to $z_3-z_1$. Thus the formula (\ref{m7a}) holds. Similarly
all other representations of $W(z_1,z_2)W(z_3,z_4)$ will give the same result. $\diamond$

Now from (\ref{m7a}) we can take a convention that the ordering (\ref{h3}) represents that
the curve represented by  $W(z_1,z_2)$ is upcrossing
the curve represented by  $W(z_3,z_4)$ while
(\ref{h2}) represents zero crossing of these two
curves.

Similarly from the dual KZ equation as a choice of branch which
is consistent with the above formula we have
\begin{equation}
W(z_1,w)W(w,z_4)W(z_3,w)W(w,z_2)=
W(z_1,w)W(w,z_2)W(z_3,w)W(w,z_4)R^{-1}
\label{m8a}
\end{equation}
where $z_2$ is before $z_4$. We take a convention that the ordering in (\ref{m8a}) represents that
the curve represented by $W(z_1,z_2)$ is undercrossing the curve represented by $W(z_3,z_4)$.
Here along the orientation of a closed curve the piece of curve
represented by $W(z_1,z_2)$ is before the piece of curve represented by
$W(z_3,z_4)$. In this case since the angle of $z_3-z_1$ minus the angle
of $z_1-z_3$ is equal to $\pi$ we have that the
angle of $z_4-z_2$ minus the angle of $z_2-z_4$ is
also equal to $\pi$ and this gives the $R^{-1}$ in this formula
(\ref{m8a}).

From (\ref{m7a}) and (\ref{m8a}) we have
\begin{equation}
 W(z_3,z_4)W(z_1,z_2)=
RW(z_1,z_2)W(z_3,z_4)R^{-1} \label{m9}
\end{equation}
where $z_1$ and $z_2$ denote the end points of a curve which is before a curve with end points $z_3$ and $z_4$.
From (\ref{m9}) we see that the algebraic structure of these
quantum Wilson lines $W(z,z')$
is analogous to the quasi-triangular quantum
group \cite{Fuc}\cite{Chari}.

\section{Computation of Quantum Wilson Loop}\label{sec10a}

Let us consider again the quantum Wilson line $W(z_1,z_4)=W(z_1, z_2)W(z_2, z_4)$.
Let us set $z_1=z_4$. In this case the quantum Wilson line forms a closed loop.
Now in (\ref{closed1a}) with $z_1=z_4$ we have that $e^{-\hat{t}\log  \pm (z_1-z_2)}$
and $e^{\hat{t}\log \pm (z_1-z_2)}$ which come from the two-side KZ
equations cancel each other and from the multivalued property of
the $\log$ function we have
\begin{equation}
W(z_1, z_1) =R^{n}{\bf A}_{14} \quad\quad n=0, \pm 1, \pm 2, ...
\label{closed2}
\end{equation}
where $R=e^{-i\pi \hat{t}}$ is the monodromy of the KZ equation \cite{Chari}.

{\bf Remark}. It is clear that if we use other representation of the quantum Wilson loop $W(z_1,z_1)$ (such as the representation $W(z_1,z_1)=W(z_1,w_1)W(w_1,w_2)W(w_2,z_1)$) then we will get the same result as (\ref{closed2}).

{\bf Remark}. For simplicity we shall drop the subscript of ${\bf A}_{14}$ in (\ref{closed2}) and simply write ${\bf A}_{14}={\bf A}$.

 \section{Defining Quantum Knots and Knot Invariant}\label{sec10}

Now
we have that the quantum Wilson loop $W(z_1, z_1)$ corresponds to a closed
curve in the complex plane with starting and ending
point $z_1$.
Let this quantum Wilson loop $W(z_1, z_1)$ represents the unknot. We shall call $W(z_1, z_1)$ as the quantum unknot. Then from
(\ref{closed2}) we have the following invariant
for the unknot:
\begin{equation}
Tr W(z_1, z_1)= Tr R^{n}{\bf A} \quad\quad n=0, \pm 1, \pm 2, ...
\label{m6}
\end{equation}
where ${\bf A}={\bf A}_{14}$ is a 2-tensor constant matrix operator.

In the following let us extend the definition (\ref{m6})
to a knot invariant for nontrivial knots.
Let $W(z_i,z_j)$ represent a piece of curve
with starting point $z_i$ and ending point $z_j$.
Then we let
\begin{equation}
W(z_1,z_2)W(z_3,z_4)
\label{m11}
\end{equation}
represent two pieces of uncrossing curve.
Then by interchanging $z_1$ and $z_3$ we have
\begin{equation}
W(z_3,w)W(w,z_2)W(z_1,w)W(w,z_4)
\label{m12}
\end{equation}
represent the curve specified by $W(z_1,z_2)$ upcrossing the
curve specified by $W(z_3,z_4)$.

Now for a given knot diagram we may cut it into a sum of
parts which are formed by two pieces of curves crossing  each other.
Each of these parts is represented
by  (\ref{m12})( For a knot diagram of the unknot
with zero crossings we simply do not need to cut the
knot diagram).
Then we define the trace of a knot with a
given knot diagram by the following form:
\begin{equation}
 Tr \cdot\cdot\cdot
 W(z_3,w)W(w,z_2)W(z_1,w)W(w,z_4)
\cdot\cdot\cdot
 \label{m14}
\end{equation}
where we use (\ref{m12})  to represent the state of the
two pieces of curves specified by
 $W(z_1,z_2)$ and
$W(z_3,z_4)$. The
 $\cdot\cdot\cdot$ means the product
of a sequence of parts represented by
(\ref{m12}) according to the state of
each part. The ordering of the sequence in (\ref{m14})
 follows the ordering of the parts given by the orientation of the
knot diagram. We shall call the sequence of crossings in
the trace (\ref{m14}) as the generalized Wilson
loop of the knot diagram. For the knot diagram of the unknot with zero crossings we simply
let it be $W(z,z)$ and call it the quantum Wilson loop.

We shall
 show that the generalized Wilson loop of a knot diagram has all the properties of the knot diagram  and that
(\ref{m14}) is  a knot invariant. From this we shall call a generalized Wilson loop as a quantum knot.

\section{Examples of Quantum Knots}

Before the proof that a generalized Wilson loop of a knot diagram has all the properties of the knot diagram
in the following let us first consider
some examples to illustrate the way to define (\ref{m14}) and
the way of applying the  braiding formulas (\ref{m7a}),
 (\ref{m8a}) and (\ref{m9}) to
equivalently transform (\ref{m14}) to a simple
expression of the form  $Tr R^{-m}W(z,z)$ where $m$
is an integer.

Let us first consider the knot in Fig.1.
For this knot we have that (\ref{m14}) is given by
\begin{equation}
Tr W(z_2,w)W(w,z_2)W(z_1,w)W(w,z_1)
\label{m15a}
\end{equation}
where the product of quantum Wilson lines  is from the definition (\ref{m12})
represented a crossing at $w$.
In applying (\ref{m12}) we let $z_1$ be the
starting and the ending point.

\begin{figure}[hbt]
\centering
\includegraphics[scale=0.6]{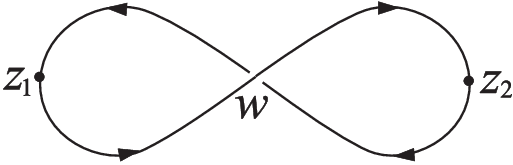}

                Fig.1
\end{figure}
Then  we have that (\ref{m15a}) is equal to
\begin{equation}
\begin{array}{rl}
&Tr W(w,z_2)W(z_1,w)W(w,z_1)W(z_2,w) \\
=&Tr RW(z_1,w)W(w,z_2)R^{-1}
RW(z_2,w)W(w,z_1)R^{-1} \\
=&Tr W(z_1,z_2)W(z_2,z_1) \\
=&Tr W(z_1,z_1)
\end{array}
\label{m16}
\end{equation}
where we have used (\ref{m9}).
We see that (\ref{m16}) is just the knot invariant (\ref{m6}) of
the unknot.
Thus the knot in Fig.1 is with the same knot invariant of the unknot and this
agrees with the fact that this knot is topologically equivalent
to the unknot.

Let us then consider a trefoil knot in Fig.2a.
By (\ref{m12}) and similar to the above examples
we have that the definition (\ref{m14})
for this knot is given by:
\begin{equation}
\begin{array}{rl}
&Tr W(z_4,w_1)W(w_1,z_2)W(z_1,w_1)W(w_1,z_5)\cdot
W(z_2,w_2)W(w_2,z_6) \\
&W(z_5,w_2)W(w_2,z_3)\cdot
W(z_6,w_3)W(w_3,z_4)W(z_3,w_3)W(w_3,z_1) \\
=&Tr W(z_4,w_1)RW(z_1,w_1)W(w_1,z_2)
R^{-1}W(w_1,z_5)\cdot
W(z_2,w_2)RW(z_5,w_2) \\
&W(w_2,z_6)R^{-1}W(w_2,z_3)\cdot
W(z_6,w_3)RW(z_3,w_3)W(w_3,z_4)R^{-1}W(w_3,z_1) \\
=&Tr
W(z_4,w_1)RW(z_1,z_2)R^{-1}W(w_1,z_5)\cdot
W(z_2,w_2)RW(z_5,z_6)R^{-1}W(w_2,z_3)\cdot \\
&W(z_6,w_3)RW(z_3,z_4)R^{-1}W(w_3,z_1) \\
=&Tr
W(z_4,w_1)RW(z_1,z_2)
W(z_2,w_2)W(w_1,z_5)W(z_5,z_6)R^{-1}W(w_2,z_3)\cdot \\
&W(z_6,w_3)RW(z_3,z_4)R^{-1}W(w_3,z_1) \\
=&Tr
W(z_4,w_1)RW(z_1,w_2)W(w_1,z_6)R^{-1}W(w_2,z_3) \\&
W(z_6,w_3)RW(z_3,z_4)R^{-1}W(w_3,z_1) \\
=&Tr
W(z_4,w_1)W(w_1,z_6)W(z_1,w_2)W(w_2,z_3) \\&
W(z_6,w_3)RW(z_3,z_4)R^{-1}W(w_3,z_1) \\
=&Tr
W(z_4,z_6)W(z_1,z_3)
W(z_6,w_3)RW(z_3,z_4)R^{-1}W(w_3,z_1) \\
=&Tr R^{-1}W(w_3,z_1)
W(z_4,z_6)W(z_1,z_3)
W(z_6,w_3)RW(z_3,z_4) \\
=&Tr
W(z_4,z_6)W(w_3,z_1)R^{-1}W(z_1,z_3)
W(z_6,w_3)RW(z_3,z_4) \\
=&Tr
RW(z_3,z_6)W(w_3,z_1)R^{-1}W(z_1,z_3)
W(z_6,w_3) \\
=&Tr
W(w_3,z_1)W(z_3,z_6)W(z_1,z_3)
W(z_6,w_3) \\
=&Tr
W(z_6,z_1)W(z_3,z_6)W(z_1,z_3)
\end{array}
\label{m21}
\end{equation}
where we have repeatly used (\ref{m9}). Then
 we have that (\ref{m21}) is equal to:
\begin{equation}
\begin{array}{rl}
&TrW(z_6,w_3)W(w_3,z_1)W(z_3,w_3)W(w_3,z_6)W(z_1,z_3)
\\
=&Tr
RW(z_3,w_3)W(w_3,z_1)W(z_6,w_3)W(w_3,z_6)W(z_1,z_3)\\
=&Tr
RW(z_3,w_3)RW(z_6,w_3)W(w_3,z_1)
R^{-1}W(w_3,z_6)W(z_1,z_3)\\
=&Tr
W(z_3,w_3)RW(z_6,z_1)
R^{-1}W(w_3,z_6)W(z_1,z_3)R\\
=&Tr
W(z_3,w_3)RW(z_6,z_3)W(w_3,z_6)\\
=&Tr W(w_3,z_6)W(z_3,w_3)RW(z_6,z_3)\\
=&Tr RW(z_3,w_3)W(w_3,z_6)W(z_6,z_3)\\
=&Tr RW(z_3,z_3)
\end{array}
\label{m22}
\end{equation}
where we have used (\ref{m7a}) and (\ref{m9}).
This is as a knot invariant for the trefoil knot in Fig.2a.
\begin{figure}[hbt]
\centering
\includegraphics[scale=0.5]{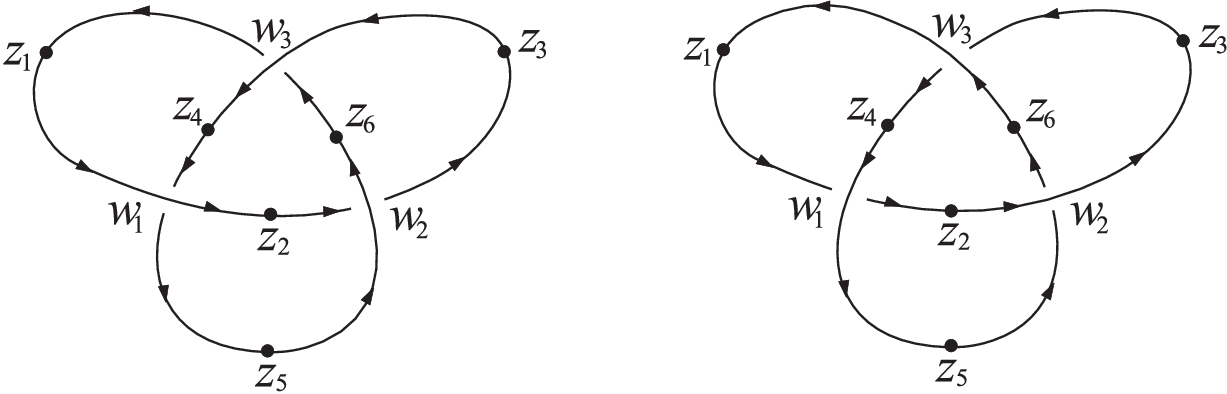}

             Fig.2a  \hspace*{5cm}  Fig.2b

\end{figure}

Then let us consider the trefoil knot in Fig. 2b
which is the mirror image of the trefoil knot in Fig.2a.
The definition (\ref{m14}) for this knot is given by:
\begin{equation}
\begin{array}{rl}
&Tr  W(z_1,w_1)W(w_1,z_5)W(z_4,w_1)W(w_1,z_2)\cdot
\\ & \qquad
W(z_5,w_2)W(w_2,z_3)W(z_2,w_2)W(w_2,z_6)\cdot \\
& \qquad W(z_3,w_3)W(w_3,z_1)W(z_6,w_3)W(w_3,z_4) \\
=&Tr W(z_5,z_1)W(z_2,z_5)W(z_1,z_2)
\end{array}
\label{m25}
\end{equation}
where similar to (\ref{m21}) we have repeatly used
(\ref{m9}).
Then we have that (\ref{m25}) is
equal to:
\begin{equation}
\begin{array}{rl}
&Tr W(z_5,z_1)W(z_2,w_1)W(w_1,z_5)W(z_1,w_1)W(w_1,z_2)\\
=&Tr
W(z_5,z_1)W(z_2,w_1)W(w_1,z_2)
W(z_1,w_1)W(w_1,z_5)R^{-1}\\
=&Tr
W(z_5,z_1)W(z_2,w_1)RW(z_1,w_1)W(w_1,z_2)
R^{-1}W(w_1,z_5)R^{-1}\\
=&Tr
R^{-1}W(z_5,z_1)W(z_2,w_1)RW(z_1,z_2)
R^{-1}W(w_1,z_5)\\
=&Tr
W(z_2,w_1)W(z_5,z_2)
R^{-1}W(w_1,z_5)\\
=&Tr
W(z_5,z_2)
R^{-1}W(w_1,z_5)W(z_2,w_1)\\
=&Tr
W(z_5,z_2)
W(z_2,w_1)W(w_1,z_5)R^{-1}\\
=&Tr
W(z_5,z_5)R^{-1}
\end{array}
\label{m26}
\end{equation}
where we have used (\ref{m8a}) and (\ref{m9}).
This is as a knot invariant for the trefoil knot in Fig.2b.
We notice that
the knot invariants for the two
trefoil knots are different. This shows that these two
trefoil knots are not topologically equivalent.

More examples of the
above quantum knots and knot invariants will be given in a following section.

\section{Generalized Wilson Loops as Quantum Knots}\label{sec11}

Let us now show that the generalized Wilson loop of a knot diagram has all the properties of the knot diagram  and that
(\ref{m14}) is  a knot invariant.
To this end let us first consider the structure of a knot. Let $K$ be
a knot. Then a knot diagram of $K$ consists of a sequence of
crossings of two pieces of curves cut from the knot $K$ where the
ordering of the crossings can be determined by the orientation of
the knot $K$.  As an example we may consider the two trefoil knots
in the above section. Each trefoil knot is represented by three
crossings of two pieces of curves. These three crossings are
ordered by the orientation of the trefoil knot starting at $z_1$.
Let us denote these three crossings by $1$, $2$ and $3$. Then the
sequence of these three crossings is given by $123$. On the other
hand if the ordering of the three crossings starts from other
$z_i$ on the knot diagram then we have sequences $231$ and $312$.
All these sequences give the same knot diagram and they can be
transformed to each other by circling as follows:
\begin{equation}
123\to 123(1) =231 \to 231(2)=312 \to 312(3)=123 \to \cdot\cdot\cdot
\label{s}
\end{equation}
where $(x)$ means that the number $x$ is to be moved to
the $(x)$ position as indicated.
Let us call (\ref{s}) as the circling property of
the trefoil knot.

As one more example let us consider the figure-eight
knot in Fig.3.
The simplest knot diagram of this knot has
four crossings.
\begin{figure}[hbt]
\centering
\includegraphics[scale=0.5]{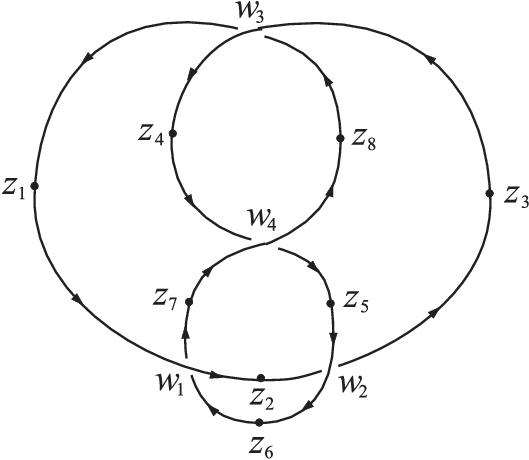}

Fig.3

\end{figure}

Starting at $z_1$ let us denote these crossings by
1, 2, 3 and 4. Then we have the following
circling property of the figure-eight
knot:
\begin{equation}
\begin{array}{rl}
& 1234\to 1234(2) =1342 \to 1342(1)=3421 \to 3421(4)=3214  \to
3214(3)=2143 \\
& \to 2143(1)=2431\to 2431(2)=4312
\to 4312(3)=4123\to 4123(4)=1234 \to \cdot\cdot\cdot
\end{array}
\label{s2}
\end{equation}
We notice that in this cirling of the figure-eight knot there are
subcirclings.

In summary we have that a knot diagram of a knot $K$
can be characterized  as a finite sequence of crossings
of curves which are cut from the knot diagram where the ordering
of the crossings is derived from the orientation of the
knot diagram and has a circling property for which
(\ref{s}) and (\ref{s2}) are examples.

Now let us represent a knot diagram of a knot $K$ by a sequence of
products of Wilson lines representing crossings as in the above
section. Let us call these products of Wilson lines by the term
W-product. Then we  call  this sequence of $W$-products as the
generalized Wilson loop of the knot diagram of a knot $K$.

Let us consider the following two $W$-products:
\begin{equation}
W(z_3,w)W(w,z_2)W(z_1,w)W(w,z_4) \quad \mbox{and}
\quad  W(z_1,z_2)W(z_3,z_4)
\label{w}
\end{equation}
 In the above section we have
shown that these two W-products faithfully represent two oriented
pieces of curves crossing or not crossing  each other where
$W(z_1,z_2)$ and $W(z_3,z_4)$ represent these two pieces of
curves.

Now there is a natural ordering of the  $W$-products
of crossings derived from the orientation of a knot as follows.
Let $W(z_1,z_2)$ and $W(z_3,z_4)$
represent two pieces of curves where the piece of curve
represented by $W(z_1,z_2)$ is before
the piece of curve
represented by $W(z_3,z_4)$  according to the
orientation of a knot. Then the ordering of these
two pieces of curves can be represented by the product
$W(z_1,z_2)W(z_3,z_4)$.
Now let 1 and 2 denote two $W$-products of crossings where we let
1  before 2 according to the orientation of a knot.
Then from the ordering of pieces of curves we have that the
product $12$ represents the ordering of the two
crossings 1 and 2.

Now let a knot diagram of a knot $K$ be given. Let the crossings
of this knot diagram be denoted by $1, 2, \cdot\cdot\cdot, n$ and
let this knot diagram be characterized by the sequence of
crossings $123\cdot\cdot\cdot n$ which is formed according to the
orientation of this knot diagram. On the other hand let us for
simplicity also denote the corresponding $W$-products of crossings
by $1, 2, \cdot\cdot\cdot, n$. Then the whole product of
$W$-products of crossings $123\cdot\cdot\cdot n$ represents the
sequence $123\cdot\cdot\cdot n$ of crossings which is identified
with the the knot diagram. This whole product $123\cdot\cdot\cdot
n$ of $W$-products of crossings is the generalized Wilson loop of
the knot diagram and we denote it by $W(K)$. In the following let us
show that this generalized Wilson loop $W(K)$ has the circling
property of the sequences of crossings of the knot diagram. It
then follows that this generalized Wilson loop represents all the
properties of the sequence $123\cdot\cdot\cdot n$ of crossings of
the knot diagram. Then since this sequence $123\cdot\cdot\cdot n$
of crossings of the knot diagram is identified with the knot
diagram we have that this generalized Wilson loop $W(K)$ can be
identified with the knot diagram. We have the following
theorem.

\begin{theorem}

Each knot $K$ can be faithfully represented by its
generalized Wilson loop $W(K)$
in the sense that if two knot diagrams have the
same generalized Wilson loop then these two knot diagrams must be topologically equivalent.
\end{theorem}
{\bf Proof.} Let us show that the generalized Wilson loop $W(K)$
of a knot diagram of $K$ has the circling property. Let us
consider a product $ W(z_1,z_2)W(z_3,z_4)$ where we first let $z_1,
z_2,z_3$ and $z_4$ be all independent. By solving the two KZ
equations as shown in the above sections we have
\begin{equation}
 W(z_1,z_2)W(z_3,z_4)
= e^{-\hat{t} \log [\pm (z_3-z_1)]}
{\bf  A}e^{\hat{t} \log [\pm (z_2-z_4)]}
\label{sa}
\end{equation}
where the initial operator ${\bf A}$ is a 4-tensor as shown in the above sections.
The sign $\pm$ in
(\ref{sa}) reflects that solutions of the KZ
equations are complex multi-valued functions.
(We remark that the 4-tensor initial operator ${\bf A}$  in general may not commute with $\Phi_{\pm}(z_1-z_2)=e^{-\hat{t} \log [\pm (z_1-z_2)]}$ and $\Psi_{\pm}(z_1-z_2)=e^{\hat{t} \log [\pm (z_1-z_2)]}$).

Then the interchange of $ W(z_1,z_2)$ and $W(z_3,z_4)$
corresponds to that
$z_1$ and $z_3$ interchange their positions and
$z_2$ and $z_4$ interchange their positions respectively.
This interchange
gives a pair of sign changes:
\begin{equation}
(z_3-z_1) \to (z_1-z_3)\qquad
\mbox{and} \qquad
(z_2-z_4) \to (z_4-z_2)
\label{sign}
\end{equation}
From this we have that $W(z_3,z_4)W(z_1,z_2)$
is given by
\begin{equation}
 W(z_3,z_4)W(z_1,z_2)
= e^{-\hat{t} \log [\pm (z_1-z_3)]}
 {\bf A}e^{\hat{t} \log [\pm (z_4-z_2)]}
\label{sb}
\end{equation}

Now let us set $z_2=z_3$ and $z_1=z_4$ such that the two products
in (\ref{sa}) and (\ref{sb}) form a closed loop. In this case we
have that the initial operator ${\bf A}$ is reduced from a 4-tensor to a 2-tensor
and that
$\Phi_{\pm}$ and $\Psi_{\pm}$ act on ${\bf A}$ by the usual matrix
operation where ${\bf A}$, $\Phi_{\pm}$ and $\Psi_{\pm}$ are matrices of
the same dimension.
In this case we have that ${\bf A}$ commutes with $\Phi_{\pm}$ and $\Psi_{\pm}$ since $\Phi_{\pm}$ and $\Psi_{\pm}$ are Casimir operators on $SU(2)$.

Let us take a definite choice of branch such that
the sign change $z_3-z_1 \to z_1-z_3$ gives
a $i\pi$ difference from the multivalued function
$\log$. Then we have that
$\Phi_{\pm}(z_3-z_1)=R \Phi_{\pm}(z_1-z_3)$. Then since $W(z_1, z_2)$
and $W(z_3, z_4)$ represent two lines with $z_1$, $z_2$ and
$z_3$, $z_4$ as starting and ending points respectively we have that the sign change
$z_2-z_4 \to z_4-z_2$ also gives the same $i\pi$ difference from the
multivalued function $\log$. Thus we have that
$\Psi_{\pm}(z_4-z_2)=R^{-1} \Psi_{\pm}(z_2-z_4)$.
It follows from this pair of sign changes and  that $A$ commutes with
$\Phi_{\pm}$ and $\Psi_{\pm}$ we have that
$W(z_1,z_2)W(z_3,z_4)=W(z_3,z_4)W(z_1,z_2)$ when $z_2=z_3$ and
$z_1=z_4$. This proves the simplest circling property of generalized
Wilson loops.

We remark that in the above proof the pair of sign changes gives
two factors $R$ and $R^{-1}$ which cancel each other and gives the
circling property. We shall later apply the same reason of pair
sign changes to get the general circling property.
We also remark that the proof of this circling property is based on the same reason as the derivation of the braiding formulas (\ref{m7a}),
(\ref{m8a}) and (\ref{m9}) as shown in the above sections.

Let us consider a product of $n$ quantum Wilson lines
$W(z_i,z_i^{\prime})$, $i=1,...,n$, with the property that the end
points $z_i$, $z_i^{\prime}$  of these quantum Wilson lines are connected
to form a closed loop. From the analysis in the above sections we
have that this product is reduced from a tensor product
to a 2-tensor. It then follows from (\ref{tensor}) that this product is
of the following form:
\begin{equation}
\prod_{ij}\Phi_{\pm}(z_i-z_j){\bf A}
\prod_{ij}
\Psi_{\pm}(z_i^{\prime}-z_j^{\prime})
\label{class}
\end{equation}
where the initial operator ${\bf A}$ is reduced to a 2-tensor
and that the $\pm$ signs of
$\Phi_{\pm}(z_i-z_j)$ and $\Psi_{\pm}(z_i-z_j)$ are to be
determined. Then since  $\Phi_{\pm}(z_i-z_j)$ and
$\Psi_{\pm}(z_i-z_j)$ commute with ${\bf A}$ we can write (\ref{class})
in the form
\begin{equation}
\prod_{ij}\Phi_{\pm}(z_i-z_j)
\prod_{ij}
\Psi_{\pm}(z_i^{\prime}-z_j^{\prime}){\bf A}
\label{basic}
\end{equation}
where $i\neq j$.
From this formula let us derive the general circling
property as follows.

Let us consider two generalized Wilson lines denoted by 1 and 2 respectively. Here by the term generalized Wilson line we mean a product of quantum Wilson lines with two open ends. As a simple example let us consider the product $W(z,z_1)W(z_2,z)$. By definition this is a generalized Wilson line with two open ends $z_1$ and $z_2$ ($z$ is not an open end). Suppose that  the two open ends of 1 and 2 are connected. Then we want to show that $12=21$. This identity is a generalization of the above interchange of $W(z_1,z_2)$ and $W(z_3,z_4)$ with $z_2=z_3$ and $z_1=z_4$.

Because $12$ and $21$ form closed loops we have that $12$ and $21$ are products of quantum Wilson lines $W(u_i, u_k)$ (where $u_i$ and $u_k$ denote some $z_p$ or $w_q$ where we use $w_q$ to denote crossing points) such that for each pair of variables $u_i$ and $u_j$ appearing at the left side of $W(u_i, u_k)$ and $W(u_j, u_l)$ there is exactly one pair of variables $u_i$ and $u_j$ appearing at the right side of $W(u_f, u_i)$ and $W(u_g, u_j)$. Thus in the formula (\ref{basic}) (with the variables $z$, $z'$ in (\ref{basic}) denoted by variables $u$) we have that the factors $\Phi_{\pm}(u_i-u_j)$ and $\Psi_{\pm}(u_i-u_j)$ appear in pairs.

As in the above case we have that the interchange of the open ends of $12$ and $21$  interchanges $12$ to $21$. This  interchange gives changes of the factors $\Phi_{\pm}(u_i-u_j)$ and $\Psi_{\pm}(u_i-u_j)$ as follows.

Let $z_1$ and $z_2$ be the open ends of $1$ and $z_3$ and $z_4$ be the open ends of $2$ such that $z_1=z_4$ and $z_2=z_3$. Consider a factor $\Phi_{\pm}(z_1-z_3)$. The interchange of $z_1$ and $z_3$ interchanges this factor to $\Phi_{\pm}(z_3-z_1)$. Then there is another factor $\Psi_{\pm}(z_2-z_4)$. The interchange of $z_2$ and $z_4$ interchanges this factor to $\Psi_{\pm}(z_4-z_2)$. Thus this is a pair of sign changes. By the same reason and the consistent choice of branch as in the above case we have that the formula (\ref{basic}) is unchanged under this pair of sign changes.

Then let us consider a factor $\Phi_{\pm}(u_i-u_j)$ of the form $\Phi_{\pm}(z_1-u_j)$ where $u_i=z_1$ and $u_j$ is not an open end. Corresponding to this factor we have the  factor $\Phi_{\pm}(z_3-u_j)$. Then under the interchange of $z_1$ and $z_3$ we have that $\Phi_{\pm}(z_1-u_j)$ and $\Phi_{\pm}(z_3-u_j)$ change to $\Phi_{\pm}(z_3-u_j)$ and $\Phi_{\pm}(z_1-u_j)$ respectively which gives no change to the formula (\ref{basic}). A similar result holds for the interchange of $z_2$ and $z_4$ for factors $\Psi_{\pm}(z_2-u_j)$ and $\Psi_{\pm}(z_4-u_j)$.

It follows that under the interchange of the open ends of $1$ and $2$ we have the pairs of sign changes from which  the formula (\ref{basic}) is unchanged. This shows that $ 12  =  21 $.

Then we consider two generalized Wilson products of crossings which are products of crossings with four open ends respectively. Let us again denote them by $1$ and $2$. Each such generalized Wilson crossing can be regarded as the crossing of two generalized Wilson lines. Then the interchange of two open ends of the two generalized lines of $1$ with the two open ends of the two generalized lines of $2$ respectively interchanges $12$ to $21$. Then let us suppose that the open ends of these two Wilson products are connected in such a way that the products $12$ and $21$ form  closed loops. In this case we want to show that $12=21$ which is a circling property of a knot diagram. The proof of this equality is again similar to the above cases. In this case we also have that the interchange of the open ends of the two generalized Wilson crossings gives  pairs of sign changes of the factors $\Phi_{\pm}(u_i-u_j)$ and $\Psi_{\pm}(u_i-u_j)$ in $ 12 $ and $ 21 $. Then by using (\ref{basic}) we have $ 12  = 21 $.

Let us then consider two generalized Wilson products
of crossings denoted by $1$ and $2$ with open ends  connected in such a way
that two open ends of $1$ (of the four open ends of $1$) are
connected to two open ends of $2$ to form a closed loop.
We want to prove that $12=21$. This will give the subcircling property.

Since a closed loop is formed we have that
each open end of $1$ or of $2$ is connected to a closed loop. In this case
as the above cases we have that
 the products $12$ is with the initial operator ${\bf A}$ being
 a 2-tensor
since the open ends of $1$ or
$2$ do not cause ${\bf A}$ to be a tensor with tensor degree more than $2$
by their connection to the closed loop. Indeed, let $z$ be an open
end of $1$ or $2$. Then it is an end point of a quantum Wilson line
$W(z,z')$ which is a part of $1$ and $2$ such that $z'$ is on the
closed loop formed by $1$ and $2$. Then we have that this quantum Wilson
line $W(z,z')$ is  connected with the closed loop at $z'$. Since
the loop is closed
from the open end $z$ we can go continuously along the closed loop to the open end of other quantum Wilson lines connected to the closed loop.
It follows that the open end $z$ gives no additional tensor degree
to the initial operator ${\bf A}$ for the product $12$ or $21$
 and that the initial operator ${\bf A}$  is
still as the initial operator for the closed loop that it is a 2-tensor.
(This is the same reason that in the above section on the computation of quantum Wilson loop we have that two quantum Wilson lines $W(z_1,z_2)$ and $W(z_2,z_4)$  connected at $z_2$ with two open ends $z_1$ and $z_4$ is with the same 2-tensor intial operator ${\bf A}$ as the case that the two quantum Wilson lines $W(z_1,z_2)$ and $W(z_2,z_4)$ form a closed loop with  $z_1=z_4$).

Now since ${\bf A}$ is a 2-tensor
we have that ${\bf A}$, $\Phi_{\pm}$ and $\Psi_{\pm}$ are as matrices of the same dimension. In this case we have that $A$ commutes
with $\Phi_{\pm}$ and $\Psi_{\pm}$. Then by interchange the open
ends of $1$ with open ends of $2$ we interchange $12$ to $21$.
This interchange again gives  pairs of sign changes. Then since
the initial operator ${\bf A}$ commutes with $\Phi_{\pm}$ and
$\Psi_{\pm}$ we have that  $12=21$, as was to be proved. Then we
let $12$ and $21$ be connected to another generalized Wilson
product of crossing denoted by $3$ to form a closed loop. Then
from $12=21$ we have $312=321$ and $123=213$. This gives the
subcircling property of generalized Wilson loops. This subcircling
property has been illustrated in the knot diagram of the
fight-eight knot. Then from a case in the above we also have the
circling property $321=213$ between $3$ and $21$.

Continuing in this way we have the circling or
subcircling properties for generalized Wilson loops
whenever the open ends of a product of generalized Wilson lines
or crossings are connected in such a way that
among the open ends a closed loop is formed.
This shows that the generalized Wilson loop of a knot
diagram has the circling property of the knot
diagram.
With this circling property it then follows that the generalized Wilson loop of
a knot diagram completely describes the structure
of the knot diagram.

Now since the generalized Wilson loop of a knot diagram is a
complete copy of this knot diagram we have that two knot diagrams
which can be equivalently moved to each other if and only if the
corresponding generalized Wilson loops can be equivalently moved
to each other. Thus we have that if two knot diagrams have the
same generalized Wilson loop then these two knot diagrams must be
equivalent. This proves the theorem. $\diamond$

{\bf Examples of generalized Wilson loops}.
As an example of generalized Wilson loops let us
consider the trefoil knots. Starting at $z_1$ let
the W-product of crossings be denoted by 1, 2 and 3.
Then we have the following circling property of
the generalized Wilson loops of the trefoil knots:
\begin{equation}
123 = 123(1) =231 = 231(2)=312 = 312(3)=123 = \cdot\cdot\cdot
\label{cr1}
\end{equation}
As one more example let us consider the figure-eight
knot.  Starting at $z_1$ let
the W-product of crossings be denoted by 1, 2, 3
and 4.
Then we have the following circling property of
the generalized Wilson loop of the figure-eight
knot:
\begin{equation}
\begin{array}{rl}
& 1234 = 1234(2) =1342 = 1342(1)=3421 = 3421(4)=3214 =
3214(3)=2143 \\
& = 2143(1)=2431 = 2431(2)=4312
= 4312(3)=4123 = 4123(4)=1234 = \cdot\cdot\cdot
\end{array}
\label{cr2}
\end{equation}
$\diamond$

{\bf Definition}.  We may call a generalized Wilson loop of a knot diagram as a quantum knot since by the above theorem this generalized Wilson loop is a complete copy of the knot diagram. $\diamond$

From the above theorem we have the following theorem.

\begin{theorem}

Let $W(K)$ denote the generalized Wilson loop of a knot $K$. Then
we can write $W(K)$ in the form $R^{-m} W(C)=R^{-m} W(z_1,z_1)$ for some integer $m$
where $C$ denotes a trivial knot and $W(C)=W(z_1,z_1)$ denotes a Wilson loop on $C$ with starting point $z_1$ and ending point $z_1$. From this form we have that the
trace $TrR^{-m}$
is a knot invariant which
classifies knots.  Thus knots can be
classified by the integer $m$.
\end{theorem}

{\bf Proof.}
Since a generalized Wilson loop $ W(K)$ is in a closed and connected form we have that a
 generalized Wilson loop $ W(K)$ can be of the form
(\ref{basic}). Thus from the multivalued property of the $\log$
function  and the two-side cancelation in (\ref{basic}) we have
that $ W(K)$ can be of the following (multivalued) form
\begin{equation}
  W(K)= R^{-k}{\bf A}
\label{BB}
\end{equation}
for some integer $k, k=0, \pm 1, \pm 2, \pm 3, ...$. Furthermore for nontrivial knot $K$ there are some factors $R^{-k_i}$ of $R^{-k}$ coming from the braidings of Wilson lines ( for which the generalized Wilson loop $W(K)$ is formed) by braiding operations such as (\ref{m7a}) and (\ref{m8a}). Thus we can write the integer $k$ in the form $k=m+n$ for some integer $m$ and for some integer $n, n=0, \pm 1, \pm 2, ...$ where $n$ is obtained by the two-side cancelations in such a way that the cancelations are obtained when the
 Wilson lines of the knot diagram for $K$ are connected together to form a
 Wilson loop $W(C)$ where $C$ is a
 closed curve which is as an unknot and is of the same form as the knot diagram for $K$ when this knot diagram of $K$ is considered only as a closed curve in the plane (such that the upcrossings and undercrossings are changed to let $K$ be the unknot $C$). From this we have $W(C)=R^{-n}{\bf A}$ for
$n=0, \pm 1, \pm 2, ...$. Thus $W(K)$ can be written in
the following form  for some $m$:
\begin{equation}
W(K) =  R^{-m} W(C)
\label{i}
\end{equation}

This number $m$ is unique since if there is another number $m_1$ such that
$W(K)=R^{-m_1} W(C)$ then we have the equality:
\begin{equation}
R^{-m} W(C)=W(K)=R^{-m_1} W(C)
\label{i1}
\end{equation}
This shows that $R^{-m}=R^{-m_1}$ and thus $m_1=m$.

From (\ref{i}) we also have
\begin{equation}
Tr  W(K) = Tr R^{-m} W(C)
\label{i2}
\end{equation}
for some integer $m$ and that $TrR^{-m} W(C)$ is a knot invariant.

Then let us show that the invariant $TrR^{-m} W(C)$  classifies
knots. Let $K_1$ and $K_2$ be two knots with the same invariant
$TrR^{-m} W(C)$. Then $K_1$ and $K_2$ are both with the same
invariant $R^{-m} W(C)$ where the trace is omitted.
Then by the above formula (\ref{i}) we have
\begin{equation}
 W(K_1) =  R^{-m}W(C) =  W(K_2)
\label{AA1}
\end{equation}
Thus $W(K_1)$ and $W(K_2)$ can be transformed to each other. Thus
$K_1$ and $K_2$ are equivalent. Thus the invariant $TrR^{-m}W(C)$  classifies knots.
It follows that the invariant $TrR^{-m}$ classifies knots and thus knots can be classified by the integer $m$, as was to be proved. $\diamond$

 \section{More Computations of Knot Invariant}

In this section let us give more computations of the knot invariant $Tr R^{-m}$.
We shall show by computation (with the chosen braiding formulae) that the fight-eight knot ${\bf 4_1}$ is assigned with the number $m=3$ and two composite knots composed by two trefoil knots (with the names reef knot and granny knot and denoted by ${\bf 3_1\star 3_1}$ and ${\bf 3_1\times 3_1}$ respectively) are assigned with the numbers $-m=4$ and $-m=9$ respectively.
The computation is quite tedious. In the next section we shall have a more efficient way to determine the integer $m$. Readers may skip this section
for the first reading.

 Let us first consider the figure-eight knot. From the figure of this knot in a above section we
have that the knot invariant of this knot is given by:
\begin{equation}
\begin{array}{rl}
 & Tr  W(z_6, w_1)W(w_1, z_2)W(z_1, w_1)W(w_1,z_7)\cdot \\
 &W(z_2, w_2)W(w_2, z_6)W(z_5, w_2)W(w_2,z_3)\cdot \\
 &  W(z_8, w_3)W(w_3, z_4)W(z_3, w_3)W(w_3,z_1)\cdot \\
 & W(z_4, w_4)W(w_4, z_8)W(z_7, w_4)W(w_4,z_5)
\end{array}
\label{more}
\end{equation}
In the above computation we have chosen $z_1$ as the
staring point (By the circling property  we may choose any point as the starting point). By repeatedly applying the braiding formulas
(\ref{m7a}),(\ref{m8a}) and (\ref{m9})
 we have that this invariant
is equal to:
\begin{equation}
 Tr R^{-3}W(w_2, z_3)W(z_8, w_2)W(z_3, z_8)
\label{more2}
\end{equation}
Then we have that (\ref{more2}) is equal to
\begin{equation}
 Tr W(z_3, z_8)R^{-3}W(w_2, z)W(z, z_3)W(z_8,z_1)W(z_1, w_2)
\label{more3}
\end{equation}
where $W(w_2, z_3)=W(w_2, z)W(z, z_3)$ with $z$ being a point on the
line represented by  $W(w_2, z_3)$ and
that $W(z_8, w_2)=W(z_8,z_1)W(z_1, w_2)$.
Since $z_1$ is  as the starting and ending point and is an intermediate point we have the
following braiding formula:
\begin{equation}
\begin{array}{rl}
&W(w_2, z_3)W(z_8, w_2) \\
 =&W(w_2,z)W(z, z_3)W(z_8,z_1)W(z_1, w_2)\\
=& R^{-1}W(z_8,z_1)W(z, z_3)W(w_2, z)W(z_1, w_2)\\
=& R^{-1}W(z_8,z_1)W(z_1, w_2)W(w_2, z)W(z, z_3)R^{-1} \\
= & R^{-1}W(z_8, w_2)W(w_2, z_3)R^{-1}
\end{array}
\label{more4}
\end{equation}
 Thus we have that (\ref{more2}) is equal to
\begin{equation}
\begin{array}{rl}
 & Tr W(z_3, z_8)R^{-3}R^{-1}W(z_8, w_2)W(w_2, z_3)R^{-1}\\
=& Tr W(z_3, z_8)R^{-4}W(z_8, z_3)R^{-1} \\
=:&  Tr W(z_3, z_8)R^{-4} \bar W(z_8, z_3)
\end{array}
\label{more5}
\end{equation}
Then  in (\ref{more5}) we have that
\begin{equation}
\begin{array}{rl}
& \bar W(z_8, z_3)\\
 = & W(z_8, z_3)R^{-1} \\
=& W(z_8,z_1)W(z_1,w_1)W(w_1,z_2)W(z_2, z_3)R^{-1}\\
=& W(z_8,z_1)W(z_2, z_3)W(w_1,z_2)W(z_1,w_1)R R^{-1}\\
=& W(z_8,z_1)W(z_2, z_3)W(w_1,z_2)W(z_1,w_1)
\end{array}
\label{more6}
\end{equation}
This shows that $\bar W(z_8, z_3)$ is a generalized Wilson line.
Then since generalized Wilson lines are with the same braiding formulas as Wilson lines
we have that by a braiding formula similar to (\ref{more4}) (for $z_1$ as the starting and ending point and as an intermediate point) the formula (\ref{more5}) is equal to:
\begin{equation}
\begin{array}{rl}
& Tr R^{-4}\bar W (z_8, z_3)W(z_3, z_8)\\
 = & Tr R^{-4} R W(z_3, z_8)\bar W (z_8, z_3)R \\
=& Tr R^{-3}  W(z_3, z_8) W (z_8, z_3)\\
=& Tr R^{-3}  W(z_3, z_3)
\end{array}
\label{more7}
\end{equation}
where the first equality is by a braiding formula which is
similar to the braiding formula (\ref{more4}).
This is the  knot invariant for the figure-eight knot
and we have that $m=3$ for this knot.

Let us then consider the composite knot ${\bf 3_1\star 3_1}$ in Fig.4.
The trace of the generalized loop of this knot is given by (In Fig.4 one of the two $w_3$ should be $w_1^{'}$):
\begin{equation}
\begin{array}{rl}
 & Tr  W(z_4,w_1)W(w_1,z_2)W(z_1,w_1)W(w_1,z_5)\cdot \\
 & W(z_2,w_2)W(w_2,z_6)W(z_5,w_2)W(w_2,z_3)\cdot \\

& W(z_3,w_1^{'})W(w_1^{'},z_5^{'})W(z_4^{'},w_1^{'})W(w_1^{'},z_2^{'})\cdot
\\
& W(z_5^{'},w_2^{'})W(w_2^{'},z_3^{'})W(z_2^{'},w_2^{'})W(w_2^{'},z_6^{'})
\cdot\\

& W(z_3^{'},w_3^{'})W(w_3^{'},z_1^{'})W(z_6^{'},w_3^{'})W(w_3^{'},z_4^{'})
\cdot \\
 & W(z_6,w_3)W(w_3,z_4)W(z_1^{'},w_3)W(w_3,z_1)
\end{array}
\label{more8}
\end{equation}

\begin{figure}[hbt]
\centering
\includegraphics[scale=0.5]{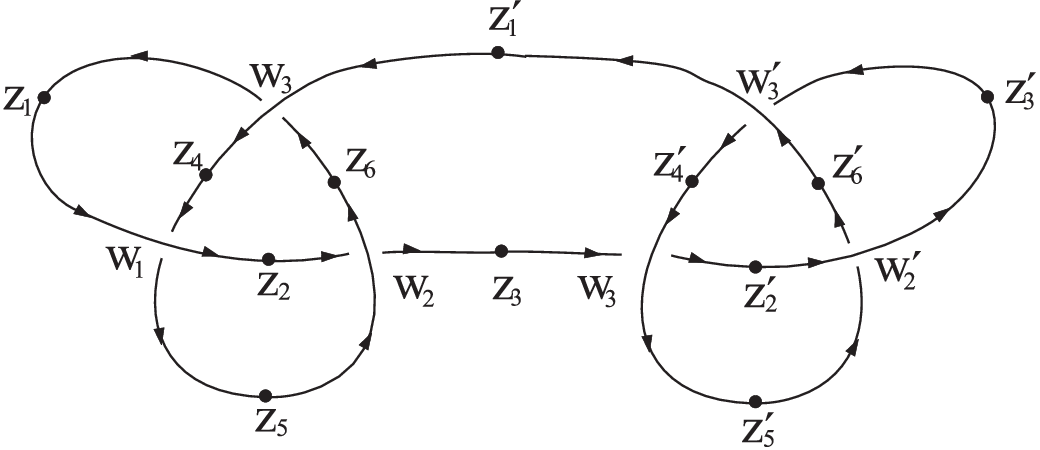}

Fig.4

\end{figure}

By repeatedly applying braiding formula (\ref{m9}) we have that
this invariant is equal to
\begin{equation}
\begin{array}{rl}
  & Tr  W(z_1,w_2^{'})W(z_5,z_1)W(w_2^{'},z_5)\\
=& Tr  W(z_1,w_2)W(w_2,w_2^{'})W(z_5,w_2)W(w_2,z_1)W(w_2^{'},z_5)\\
=& Tr W(z_1,w_2)W(w_2,z_1)W(z_5,w_2)W(w_2,w_2^{'})R^4W(w_2^{'},z_5)
\end{array}
\label{more8a}
\end{equation}
where the braiding of $W(w_2,z_1)$ and $W(w_2,w_2^{'})$ gives $R^4$.
This braiding formula comes from the fact that the Wilson line
$W(w_2,w_2^{'})$ represents a curve with end points $w_2$ and
$w_2^{'}$ such that one and a half loop
is formed which cannot be removed because the end point $w_2^{'}$
is attached to this curve itself to form the closed loop. This
closed loop gives a $3\pi$ phase angle  which is a topological
effect. Thus while the usual
braiding of two pieces of curves gives $R$ which is of a
$\pi$ phase angle  we have that the braiding of
$W(w_2,z_1)$ and $W(w_2,w_2^{'})$ gives $R$ and an additional
$3\pi$ phase angle  and thus gives $R^4$.

Then we have that (\ref{more8}) is equal to
\begin{equation}
\begin{array}{rl}
& Tr W(w_2^{'},z_5)W(z_1,w_2)W(w_2,z_1)W(z_5,w_2)W(w_2,w_2^{'})R^4\\
= & Tr W(w_2^{'},z_5)W(z_1,w_2)RW(z_5,w_2)W(w_2,z_1)
R^{-1}W(w_2,w_2^{'})R^4\\
= & Tr W(w_2^{'},z_5)W(z_1,w_2)RW(z_5,z_1)R^{-1}W(w_2,w_2^{'})R^4 \\
= & Tr RW(z_1,w_2)W(w_2^{'},z_5)W(z_5,z_1)R^{-1}W(w_2,w_2^{'})R^4 \\
= & Tr RW(z_1,w_2)W(w_2^{'},z_1)R^{-1}W(w_2,w_2^{'})R^4 \\
= & Tr W(w_2^{'},z_1)W(z_1,w_2)W(w_2,w_2^{'})R^4 \\
= & Tr W(w_2^{'},w_2)W(w_2,w_2^{'})R^4 \\
= & Tr W(w_2^{'},w_2^{'})R^4
\end{array}
\label{more9}
\end{equation}
This is the invariant of ${\bf 3_1\star 3_1}$.
 Thus we have that $-m=4$ for ${\bf 3_1\star 3_1}$.

Let us then consider the composite knot ${\bf 3_1\times 3_1}$ in Fig.5.
We have that the trace of the generalized Wilson loop of
${\bf 3_1\times 3_1}$ is given by (In Fig.5 one of the two $w_3$ should be $w_1^{'}$):
\begin{equation}
\begin{array}{rl}
 & Tr  W(z_4,w_1)W(w_1,z_2)W(z_1,w_1)W(w_1,z_5)\cdot \\
 & W(z_2,w_2)W(w_2,z_6)W(z_5,w_2)W(w_2,z_3)\cdot \\
& W(z_4^{'},w_1^{'})W(w_1^{'},z_2^{'})W(z_3,w_1^{'})W(w_1^{'},z_5^{'})
\cdot \\
 & W(z_2^{'},w_2^{'})W(w_2^{'},z_6^{'})W(z_5^{'},w_2^{'})W(w_2^{'},z_3^{'})
\cdot \\
& W(z_6^{'},w_3^{'})W(w_3^{'},z_4^{'})W(z_3^{'},w_3^{'})W(w_3^{'},z_1^{'})
\cdot \\
 & W(z_6,w_3)W(w_3,z_4)W(z_1^{'},w_3)W(w_3,z_1)
\end{array}
\label{more10}
\end{equation}

\begin{figure}[hbt]
\centering
\includegraphics[scale=0.5]{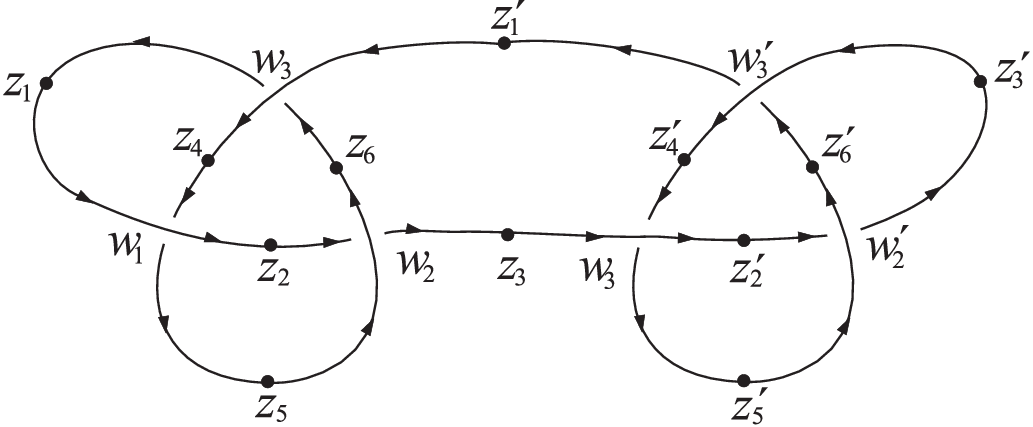}

Fig.5

\end{figure}

By repeatedly applying braiding formulas (\ref{m7a}), (\ref{m8a})
and (\ref{m9}) we have that this invariant is equal to
\begin{equation}
Tr RW(z_1^{'},w_1)R^2W(z_4^{'},z_6^{'})W(w_1,z_4^{'})W(z_6^{'},z_1^{'})
\label{more11}
\end{equation}
where the quantum Wilson line $W(w_1,z_4^{'})$ represents the piece of curve
which starts at $w_1$ and goes through $z_5$, $z_6$, $z_1$ and ends
at $z_4^{'}$. This curve includes a one and a half loop which cannot be
removed since $w_1$ is attached to this curve to form the loop.
This is of the same case as that in the knot ${\bf 3_1\star 3_1}$.
This is a topological property which gives a $3\pi$ phase angle.

We have that (\ref{more11}) is equal to
\begin{equation}
Tr RW(z_1^{'},w_1)R^2
W(z_4^{'},w_1^{'})W(w_1^{'},z_6^{'})
W(w_1,w_1^{'})W(w_1^{'},z_4^{'})W(z_6^{'},z_1^{'})
\label{more12}
\end{equation}
where the piece of curve represented by quantum Wilson line $W(w_1,w_1^{'})$ also
contains the closed loop. Now let this knot ${\bf 3_1\times 3_1}$ be
starting  and ending at $z_6^{'}$.
 Then by the braiding formula on $W(w_1,w_1^{'})$ and $W(z_4^{'},w_1^{'})$
as in the case of the knot ${\bf 3_1\star 3_1}$ we have that (\ref{more12})
is equal to
\begin{equation}
\begin{array}{rl}
& Tr RW(z_1^{'},w_1)R^2 \\
& R^4 W(w_1,w_1^{'})W(w_1^{'},z_6^{'})W(z_4^{'},w_1^{'})
W(w_1^{'},z_4^{'})W(z_6^{'},z_1^{'})\\
= & Tr RW(z_1^{'},w_1)R^6 \\
 & W(w_1,w_1^{'})RW(z_4^{'},w_1^{'})W(w_1^{'},z_6^{'})R^{-1}
 W(w_1^{'},z_4^{'})W(z_6^{'},z_1^{'})\\
= & Tr RW(z_1^{'},w_1)R^6 \\
 & W(w_1,w_1^{'})RW(z_4^{'},z_6^{'})
R^{-1}W(w_1^{'},z_4^{'})W(z_6^{'},z_1^{'})\\
= & Tr RW(z_1^{'},w_1)R^6 \\
 & W(w_1,w_1^{'})RW(z_4^{'},z_6^{'})
W(z_6^{'},z_1^{'})W(w_1^{'},z_4^{'})R^{-1}\\
= & Tr RW(z_1^{'},w_1)R^6 \\
& W(w_1,w_1^{'})RW(z_4^{'},z_1^{'})W(w_1^{'},z_4^{'})R^{-1}\\
= &Tr W(z_1^{'},w_1)R^6
 W(w_1,w_1^{'})RW(z_4^{'},z_1^{'})W(w_1^{'},z_4^{'})
\end{array}
\label{more13}
\end{equation}
where we have repeatedly applied the braiding formula (\ref{m9}).

Now let $z_4^{'}$ be the starting and ending point. Then we have that
(\ref{more13}) is equal to
\begin{equation}
\begin{array}{rl}
& Tr W(z_1^{'},w_1)R^6
 W(w_1,w_1^{'})W(w_1^{'},z_4^{'})W(z_4^{'},z_1^{'})R \\
= & Tr W(z_1^{'},w_1)R^6 W(w_1,z_1^{'})R \\
=& Tr W(z_1^{'},w_1)R^6
 W(w_1,w_1^{'})W(w_1^{'},z_4^{'})W(z_4^{'},w_1^{'})
W(w_1^{'},z_1^{'})R \\
=& Tr W(z_1^{'},w_1)R^6
 W(w_1,w_1^{'})W(w_1^{'},z_1^{'})W(z_4^{'},w_1^{'})W(w_1^{'},z_4^{'}) \\
=: & Tr W(z_1^{'},w_1)R^6 \bar W(w_1,z_1^{'}) \\
= & Tr R^6 \bar W (w_1,z_1^{'}) W(z_1^{'},w_1)
\end{array}
\label{more14}
\end{equation}
where $ \bar W(w_1,z_1^{'})$ denotes the following generalized Wilson line:
\begin{equation}
W(w_1,w_1^{'})W(w_1^{'},z_1^{'})W(z_4^{'},w_1^{'})W(w_1^{'},z_4^{'})
\label{more15}
\end{equation}
Then by the same braiding formula for generalized Wilson lines as that for Wilson lines (with $z_4^{'}$ as the starting and ending point and as an intermediate point)
we have that (\ref{more14}) is equal to:
\begin{equation}
\begin{array}{rl}
& Tr R^6 \bar W (w_1,z_1^{'}) W(z_1^{'},w_1) \\
= & Tr R^6 R W(z_1^{'},w_1) \bar W (w_1,z_1^{'})R \\
= & Tr R^6 R W(z_1^{'},w_1) W (w_1,z_1^{'})RR \\
= & Tr R^9 W(z_1^{'},z_1^{'})
\end{array}
\label{more15a}
\end{equation}
This is the knot invariant for the knot ${\bf 3_1\times 3_1}$.  Thus we have that $-m=9$ for the knot ${\bf 3_1\times 3_1}$. Then we have that the image of ${\bf 3_1\times 3_1}$ is with the knot invariant $Tr R^{-9} W(z_1^{'},z_1^{'})$.


\section{A Classification Table of Knots I}\label{sec1a}

In the above sections the computations of the knot invariant $Tr R^{-m}$ is tedious. In this section let us use another method to determine the integer $m$ without carrying out the tedious computations. We shall use only the connected sum operation on knots to find out the integer $m$. For simplicity we use the positive integer $|m|$ to form a classification table of knots where $m$ is assigned to a knot while $-m$ is assigned to
its mirror image if the knot is not equivalent to its mirror image.
Our main references on the connected sum operation on knots are \cite{Ada}-\cite{Rol}.

Let $\star$ denote the connected sum of two knots such that the
resulting total number of alternating crossings is equal to the
sum of alternating crossings of each of the two knots minus $2$.
As an example we have the reef knot (or the square knot) ${\bf
3_1\star 3_1}$ which is a composite knot composed with the knot
$\bf 3_1$ and its mirror image as in Fig.4.
This square knot has $6$ crossings
and $4$ alternating crossings. Then let $\times$ denote the
connected sum for two knots such that the resulting total number
of alternating crossings is equal to the sum of alternating
crossings of each of the two knots. As an example we have the granny
knot ${\bf 3_1\times 3_1}$ which is a composite knot composed with
two identical knots $\bf 3_1$ as in Fig.5
(For simplicity we use one notation $\bf 3_1$ to denote both the trefoil knot and its mirror image though these two knots are nonequivalent). This knot has $6$ alternating crossings which is equal to the total number of crossings. We have that the two operations $\star$ and $\times$ satisfy the commutative law and the associative law
 \cite{Ada}-\cite{Rol}. Further for each knot there is a unique factorization of this knot into a $\star$ and $\times$ operations of prime knots which is similar to the unique factorization of a number into a product of prime numbers \cite{Ada}-\cite{Rol}. We shall show that there is a deeper connection between these two factorizations.

We shall show that we can establish a classification table of knots where each knot is assigned with a number such that
prime knots are bijectively assigned with prime numbers such that the prime number $2$ corresponds to the trefoil knot
(The trefoil knot will be assigned with the number $1$ and is related to the prime number $2$). We have shown by computation  that the knot $\bf
3_1$ is with $m=1$, the knot $\bf 4_1$ is with $m=3$. Thus there are
no knots assigned with the number $2$ since other knots are with
crossings more than these two knots. We have
shown by computation that the knot ${\bf 3_1\star 3_1}$ is assigned with the
number $4$. Thus we have $1\star 1=4$ (Since knots are assigned with integers we may regard the $\star$ and $\times$ as operations on the set of numbers). This shows that the number
$1$ plays the role of the number $2$. Thus while the knot $\bf
3_1$ is with $m=1$ we may regard this $m=1$ is as the even prime
number $2$. We shall have more to say about this phenomenon
of $1$ and $2$. This phenomenon reflects that the operation $\star$ has partial properties  of addition and multiplication where $m=1$ is assigned to ${\bf 3_1}$ for addition while ${\bf 3_1}$ plays the role of $2$ is for multiplication. The aim of this section is to find out a table of the relation between knots and numbers by using only the operations $\star$ and $\times$ on knots and by using the following data as the initial step for induction:

{\bf Initial data for induction:} The prime knot ${\bf 3_1}$ is assigned with the number $1$ and it also plays the role of $2$. This means that the number $2$ is not assigned to other knots and is left for the prime knot ${\bf 3_1}$.
$\diamond$

{\bf Remark}. We shall say that the prime knot ${\bf 3_1}$ is assigned with the number $1$ and is related to the prime number $2$.
$\diamond$

We shall give an induction
on the number $n$ of $2^n$ for establishing the table. For each induction step on $n$ because of the special role of the trefoil knot ${\bf 3_1}$ we let the composite knot ${\bf 3_1}^n$ obtained by repeatedly taking $\star$ operation $n-1$ times on the trefoil knot ${\bf 3_1}$ be assigned with the number $2^n$ in this induction.

Let us first give the following table relating knots and numbers up to $2^5$ as a guide for the induction
for establishing the whole classification table of knots:

\begin{displaymath}
\begin{array}{|c|c|c|c|} \hline
\mbox{Type of Knot}& \mbox{Assigned number} \,\, |m|
 &\mbox{Type of Knot}& \mbox{Assigned number} \,\, |m|
\\ \hline
{\bf 3_1} & 1 & {\bf 6_3} &  17\\ \hline

 &  2
&  {\bf 3_1\times 4_1} &  18 \\ \hline

{\bf 4_1} &  3 & {\bf 7_1} &  19 \\ \hline

{\bf 3_1\star 3_1} &  4 &
{\bf 4_1\star 5_1} &  20
\\ \hline

{\bf 5_1} & 5 & {\bf 4_1\star(3_1\star 4_1) } &  21
\\ \hline

{\bf 3_1\star 4_1} & 6 & {\bf 4_1\star 5_2} & 22 \\ \hline

{\bf 5_2} &  7 & {\bf 7_2} & 23 \\ \hline

{\bf 3_1\star 3_1\star 3_1} &  8 &
{\bf 3_1\star (3_1\times 3_1)}& 24 \\ \hline

{\bf 3_1\times 3_1} & 9 &
{\bf 3_1\star (3_1\star  5_1)}& 25 \\ \hline

{\bf 3_1\star 5_1} &  10 & {\bf 3_1\star 6_1}
& 26 \\ \hline

{\bf 6_1} &  11 & {\bf 3_1\star (3_1\star  5_2)}& 27
\\ \hline

{\bf 3_1\star 5_2} &  12 & {\bf 3_1\star 6_2} & 28
\\ \hline

{\bf 6_2} &  13 & {\bf 7_3} & 29 \\ \hline

{\bf 4_1\star 4_1} &  14 & { \bf (3_1\star 3_1)\star (3_1\star
4_1)} & 30 \\ \hline

{\bf 4_1\star (3_1\star 3_1)} & 15 & {\bf 7_4} & 31
\\ \hline

{ \bf (3_1\star 3_1)\star (3_1\star 3_1)} & 16 & {\bf (3_1\star
3_1)\star (3_1\star 3_1\star 3_1)} &  32
\\ \hline
\end{array}
\end{displaymath}

From this table we see that the $\star$ operation is similar to the usual multiplication $\cdot$ on numbers. Without the $\times$ operation this $\star$ operation would be exactly the usual multipilcation on numbers if this $\star$ operation is regarded as an operation on numbers.
From this table we see that comparable composite knots (in a sense
from the table and we shall discuss this point later) are grouped
in each of the intervals between two prime numbers. It is
interesting that in each interval composite numbers are one-to-one
assigned to the comparable composite knots while prime numbers are one-to-one
assigned to prime knots. Here a main point is to introduce the $\times$ operation while keeping composite knots correspond to composite numbers and prime knots correspond to prime numbers. To this end we need to have rooms at the positions of composite numbers for the introduction of composite knots obtained by the $\times$ operation. We shall show that these rooms can be obtained by using the special property of the trefoil knot which is assigned with the number $1$ (for the addition property of the $\star$ and $\times$ operations) while this trefoil knot is similar to the number $2$ for the multiplication property of the $\star$ operation.

Let us then carry out the induction steps for obtaining the whole table. To this end let us investigate in more detail the above comparable properties of knots. We have the following definitions and theorems.

{\bf Definition}.
We write $K_1<K_2$ if
$K_1$ is before $K_2$ in the ordering of knots; i.e. the number assigned to $K_1$ is less than the number assigned to $K_2$.

{\bf Definition (Preordering)}. Let two knots be
written in the form $K_1\star K_2$ and $K_1\star K_3$ where we
have determined the ordering of $K_2$ and $K_3$. Then we say that
$K_1\star K_2$ and $K_1\star K_3$ are in a preordering in the sense that we put the ordering of these two knots to follow the ordering of $K_2$ and $K_3$. If this preordering is not changed by conditions from other preorderings on these two knots (which are from other factorization forms of these two knots) then this preordering becomes the ordering of these two knots.
We shall see that this preordering gives the comparable property in the above table. $ \diamond$

{\bf Remark}.
a) 
If $K_1$ is the unknot then we have
$K_1\star K_2$=$K_2$ and $K_1\star K_3$=$K_3$ and thus the ordering of $K_1\star K_2$ and $K_1\star K_3$ follows the ordering of $K_2$ and $K_3$.

b) We can also similarly define the preordering of two knots $K_1\times K_2$ and
$K_1\times K_3$ with the $\times$ operation. $\diamond$

We have the following theorem.
\begin{theorem}

Consider two knots of the form $K_1\star K_2$ and $K_1\star K_3$ where $K_1$, $K_2$ and $K_3$ are prime knots such that $K_2<K_3$.
Then we have $K_1\star K_2 < K_1\star K_3$.
\label{pre}
\end{theorem}
{\bf Proof}.
Since $K_1$, $K_2$ and $K_3$ are prime knots there are no other factorization forms of the two knots $K_1\star K_2$ and $K_1\star K_3$. Thus these two forms of the two knots are the only way to give preordering to the two knots and thus there are no other conditions to change the preordering given by this factorization form of the two knots. Thus we have that $K_2<K_3$ implies $K_1\star K_2 < K_1\star K_3$.
$\diamond$

\begin{theorem}

Suppose two knots are written in the form $K_1\star K_2$ and $K_1\star K_3$ for determining their ordering and that the other forms of these two knots are not for determining their ordering. Suppose that $K_2<K_3$. Then we have $K_1\star K_2 < K_1\star K_3$.
\end{theorem}
{\bf Proof}.
The proof of this theorem is similar to the proof of the above theorem. Since the other factorization forms are not for the determination of the ordering of the two knots in the factorization form $K_1\star K_2$ and $K_1\star K_3$ we have that the preordering of these two knots in this factorization form becomes the ordering of these two knots. Thus we have $K_1\star K_2 < K_1\star K_3$. $\diamond$

As a generalization of theorem \ref{pre} we have the following theorem.

\begin{theorem}

Let two knots be of the form $K_1\star K_2$ and $K_1\star K_3$ where $K_2$ and $K_3$ are prime knots. Suppose that $K_2<K_3$. Then we have $K_1\star K_2 < K_1\star K_3$.
\label{pre2}
\end{theorem}
{\bf Proof}.
We have the preordering that $K_1\star K_2$ is before $K_1\star K_3$. Then since $K_2$ and $K_3$ are prime knots we have that the other preordering of $K_1\star K_2$ and $K_1\star K_3$ can only from the factorization of $K_1$. Without loss of generality let us suppose that $K_1$ is of the form $K_1=K_4\star K_5$ where
$K_4< K_5$ and $K_4$ and $K_5$ are prime knots. Then we have the factorization
$K_1\star K_2=K_4\star (K_5\star K_2)$ and $K_1\star K_3=K_5\star (K_4\star K_3)$.
This factorization is the only factorization that might change the preordering that $K_1\star K_2$ is before $K_1\star K_3$. Then if $K_2\neq K_4$ or $K_3\neq K_5$
with this  factorization the two knots $K_1\star K_2$ and $K_1\star K_3$ are noncomparable in the sense that this factorization gives no preordering property and that the ordering of these two knots is determined by other conditions. Thus this factorization of $K_1\star K_2$ and $K_1\star K_3$ is not for the determination of the ordering of $K_1\star K_2$ and $K_1\star K_3$. Thus the preordering that $K_1\star K_2$ is before $K_1\star K_3$ is the ordering of $K_1\star K_2$ and $K_1\star K_3$. On the other hand if $K_2=K_4$ and $K_3=K_5$ then this factorization gives the same preordering that $K_1\star K_2$ is before $K_1\star K_3$. Thus for this case the preordering that $K_1\star K_2$ is before $K_1\star K_3$ is also the ordering of $K_1\star K_2$ and $K_1\star K_3$. Thus we have $K_1\star K_2 <K_1\star K_3$. $\diamond$

In addition to the above theorems we have the following theorems.
\begin{theorem}

Consider two knots of the form $K_1\times K_2$ and $K_1\times K_3$ where $K_1$, $K_2$ and $K_3$ are prime knots such that $K_2<K_3$.
Then we have $K_1\times K_2 < K_1\times K_3$.
\label{1}
\end{theorem}
{\bf Proof}. By using a preordering property for knots with $\times$ operation as similar to that for knots with $\star$ operation we have that the proof of this theorem is similar to the proof of the above theorems. $\diamond$

\begin{theorem}

Let two knots be of the form $K_1\times K_2$ and $K_1\times K_3$ where $K_2$ and $K_3$ are prime knots. Suppose that $K_2<K_3$. Then we have $K_1\times K_2 < K_1\times K_3$.
\end{theorem}
{\bf Proof}.
The proof of this theorem is also similar to the proof of the theorem \ref{pre2}. $\diamond$

{\bf Remark}. These two theorems will be used for introducing and ordering knots involved with a $\times$ operation which will have the effect of pushing out composite knots with the property of jumping over (to be defined) such that knots are assigned with a prime number if and only if the knot is a prime knot. $\diamond$

{\bf Remark}. The above theorems on preordering are used to find out  the ordering of knots.

Let us investigate more on the property of preordering. We consider the following

{\bf Definition (Preordering sequences)}.
 At the $n$th induction step let the prime knot ${\bf 3_1}$ take a $\star$ operation with the previous $(n-1)$th step. We call this obtained sequence of composite knots as a
preordering sequence. Thus from the ordering of the $(n-1)$th step we have a sequence of composite knots which will be for the construction of the $n$th step.

Then we let the prime knot ${\bf 4_1}$ (or the knot assigned with a prime number which is $3$ in the $2$nd step as can be seen from the above table) take a $\star$ operation with the previous $(n-2)$th step. From this we get a sequence of composite knots for constructing the $n$th step.
Then we let the prime knots ${\bf 5_1}$ and ${\bf 5_2}$ (which are prime knots in the same step assigned with a prime number which is $5$ or $7$ in the $3$rd step as can be seen from the above table) take a $\star$ operation with the previous $(n-3)$th step respectively. From this we get two sequences for constructing the $n$th step.

Continuing in this way until the sequences are obtained by a prime knot in the $(n-1)$th step taking a $\star$ operation with the step $n=1$ where the prime knot is assigned with a prime number in the $(n-1)$th step by induction (By induction each prime number greater than $2$ will be assigned to a prime knot).

We call these obtained sequences of composite knots as
the preordering sequences of composite knots for constructing the $n$th step. Also we call
the sequences truncated from these preordering sequences as preordering subsequences of composite knots for constructing the $n$th step. $\diamond$

We first have the follwing lemma on preordering sequence.
\begin{lemma}

Let $K$ be a knot in a preordering sequence of the $n$th step.
Then there exists a room for this $K$ in the $n$th step in the sense that this $K$ corresponds to a number in the $n$th step or in the $(n-1)$th step.
\end{lemma}
{\bf Proof}.
Let $K$ be of the form $K= {\bf3_1}\star K_1$ where $K_1$ is a knot in the previous $(n-1)$th step. By induction we have that $K_1$ is assigned with a number $a$ which is the position of $K_1$ in the previous $(n-1)$th step. Then since ${\bf3_1}$ corresponds to the number $2$ we have that $K$ corresponds to the number $2\cdot a$ in the $n$th step (We remark that
$K$ may not be assigned with the number $2\cdot a$). Thus there exists a room for this $K$ in the $n$th step.

Then let $K$ be of the form $K= {\bf4_1}\star K_2$ where $K_2$ is a knot in the previous $(n-2)$th step. By induction we have that $K_2$ is assigned with a number $b$ which is the position of $K_2$ in the previous $(n-2)$th step. Since ${\bf4_1}$ is by induction assigned with the prime number $3$ we have $3\cdot b >3\cdot2^{n-3}>2\cdot2^{n-3}=2^{n-2}$. Also we have
$3\cdot b <3\cdot2^{n-2}<2^2\cdot2^{n-2}=2^{n}$. Thus there exists a room for this $K$ in the $(n-1)$th step or the $n$th step.

Continuing in this way we have that this lemma holds. $\diamond$

{\bf Remark}. By using  this lemma we shall construct each $n$th step of the classification table by first filling the $n$th step with the preordering subsequences of the $n$th step. $\diamond$

{\bf Remark}.
When the number corresponding to the knot $K$ in the above proof is not in the $n$th step we have that the knot $K$ in the preordering sequences of the $n$th step has the function of pushing a knot $K^{\prime}$ out of the $n$th step where this knot $K^{\prime}$ is related to a number in the $n$th step in order for the knot $K$ to be filled into the $n$th step.

As an example in the above table the knot $K= {\bf4_1}\star{\bf5_1} $ (related to the number $3\cdot5$) in a preordering sequence of the $5$th step pushes the knot $K^{\prime}= {\bf5_1}\star {\bf5_1}$ related to the number $5\cdot5$ in the $5$th step out of the $5$th step. This relation of pushing out is by the chain $3\cdot5 \to 2\cdot2\cdot5 \to 5\cdot5$.

As another example in the above table the knot $K= {\bf3_1}\star({\bf3_1}\times{\bf3_1})$ (correspoded to the number $2\cdot9$) in a preordering sequence of the $5$th step pushes the knot $K^{\prime}= {\bf3_1}\star({\bf4_1}\star {\bf5_1})$ related to the number $2\cdot3\cdot5$ in the $5$th step out of the $5$th step. This relation of pushing out is by the chain $2\cdot9 \to 2\cdot2\cdot2\cdot3 \to 2\cdot3\cdot5$.
$\diamond$

\begin{lemma}

For $n\geq 2$ the preordering subsequences for the $n$th step can cover the whole $n$th step.
\end{lemma}
{\bf Proof}.
For $n=2$ we have one preordering sequence with number of knots $=2^0$ which is obtained by the prime knot ${\bf 3_1}$ taking $\star$ operation with the step $n=2-1=1$.
In addition we have the knot ${\bf 3_1}\star{\bf 3_1}$ which is assigned at the position of $2^n,n=2$ by the induction procedure.  Then since  the total rooms of this step $n=2$ is $2^1$ we have that
these two knots cover this step $n=2$.

For $n=3$ we have one preordering sequence with number of knots $=2^1$ which is obtained by the prime knot ${\bf 3_1}$ taking $\star$ operation with the step $3-1=2$. This sequence cover half of this step $n=3$ which is with $2^{3-1}=2^2$ rooms.
Then we have one more preordering sequence  which is obtained by the knot ${\bf 4_1}$ taking $\star$ operation with step $n=1$ giving the number $2^0=1$ of knots. This covers half of the remaining rooms of the step $n=3$ which is with $2^{2-1}=2^1$ rooms. Then
in addition we have the knot ${\bf 3_1}\star{\bf 3_1}$ which is assigned at the position of $2^n,n=2$ by the induction procedure. The total of these four knots thus cover the step $n=3$.

For the $n$th step we have one preordering sequence with the number of knots $=2^{n-2}$ which is obtained by the prime knot ${\bf 3_1}$ taking $\star$ operation with the $n-1$th step. This sequence cover half of this $n$th step which is with $2^{n-1}$ rooms. Then we have a preordering sequence  which is obtained by the knot ${\bf 4_1}$ taking $\star$ operation with the $(n-2)$th step giving the number $2^{n-3}$ of knots.  This covers half of the remaining rooms of the $n$th step which is with the remaining $2^{n-2}$ rooms. Then we have one preordering sequence obtained by picking a prime knot (e.g.${\bf 5_1}$) which by induction is assigned with a prime number (e.g. the number 5) taking $\star$ operation with the $(n-3)$th step.
Continue in this way until the knot ${\bf 3_1}^n$
is by induction assigned at the position of $2^n$. The total number of these knots is $2^{n-1}$ and thus cover this $n$th step. This proves the lemma. $\diamond$

{\bf Remark}.
Since there will have more than one prime number in the $k$th steps ($k>2$)
in the covering
of the $n$th step there will have knots from
the preordering sequences in repeat and in overlapping. These knots in repeat and in overlapping  may be deleted when  the ordering of the subsequences of the preordering sequences has been determinated for the covering of the $n$th step.

Also in the preordering sequences some knots which are in repeat and are not used for the covering of the $n$th step will be omitted when the ordering of the subsequences of the preordering sequences has been determinated for the covering of the $n$th step.
$\diamond$

Let us then introduce another definition for constructing the classification table of knots.

{\bf Definition (Jumping over of the first kind)}. At an induction $n$th step consider a knot
$K^{\prime}$ and the knot $K={\bf 3_1}^n$ which is a $\star$ product of $n$ knots
$3_1$. $K^{\prime}$ is said to jump over $K$, denoted by $K \prec
K^{\prime}$,
 if exist $K_2$ and $K_3$  such that
$K^{\prime}=K_2\star K_3$ and for any $K_0$, $K_1$
such that  $K = K_0\star K_1$ where $K_0$, $K_1$, $K_2$ and $K_3$
are not
equal to ${\bf 3_1}$ we have
\begin{equation}
2^{n_0}<p_1\cdot\cdot\cdot p_{n_2}, \quad
2^{n_1}>q_1\cdot\cdot\cdot q_{n_3}
\label{class1}
\end{equation}
or vice versa
\begin{equation}
2^{n_0}>p_1\cdot\cdot\cdot p_{n_2}, \quad
2^{n_1}<q_1\cdot\cdot\cdot q_{n_3}
\label{class2}
\end{equation}
where $2^{n_0}$, $2^{n_1}$ are the numbers assigned to $K_0$ and $K_1$ respectively ($n_0 +n_1=n$) and
\begin{equation}
K_2= K_{p_1}\star \cdot\cdot\cdot \star K_{p_{n_2}} \quad K_3=
K_{q_1}\star \cdot\cdot\cdot \star K_{q_{n_3}}
\label{class3}
\end{equation}
where $K_{p_i}$, $K_{q_j}$ are prime knots which have been
assigned with prime integers $p_i$, $q_j$ respectively; and the following inequality holds:
\begin{equation}
2^n=2^{n_0+n_1}>p_1\cdot\cdot\cdot p_{n_2}\cdot q_1\cdot\cdot\cdot
q_{n_3}
\label{class33}
\end{equation}
Let us call this definition as the property of jumping over of the first kind. $\diamond$

We remark that the definition of jumping over of the first kind is a generalization of the above
ordering of ${\bf 4_1}\star{\bf 5_1}$ and ${\bf 3_1\star 3_1\star
3_1\star 3_1}$ in the above table in the step $n=4$ of $2^4$.
Let us consider some examples of this
definition. Consider the knots $K^{\prime}=K_2\star K_3={\bf
4_1\star 5_1}$ and $K={\bf 3_1\star 3_1\star 3_1\star 3_1}$. For
any $K_0$, $K_1$ which are not equal to ${\bf 3_1}$ such that
$K=K_0\star K_1$ we have $2^{n_0}< 5$ and $2^{n_1}> 3$ (or vice
versa) where $3$, $5$ are the numbers of ${\bf 4_1}$ and
${\bf 5_1}$ respectively. Thus we have that ${\bf (3_1\star
3_1)\star(3_1\star 3_1)} \prec {\bf 4_1\star 5_1}$.

As another example we have that
$ {\bf 3_1\star(3_1\star 3_1)\star(3_1\star 3_1)} \prec
{\bf 5_1\star 5_1}$, ${\bf 4_1\star 4_1\star 4_1 }$,  and
${\bf 3_1\star (4_1\star 5_1)}$.

{\bf A Remark on Notation}. At the $n$th step let a composite knot of the form $K_1\star K_2\star\cdot\cdot\cdot\star K_q$ where each $K_i$ is a prime knot such that $K_i$ is assigned with a prime number $p_i$ in the previous $n-1$ steps. Then in general $K_1\star K_2\star\cdot\cdot\cdot\star K_q$ is not assigned with the number $p_1\cdot\cdot\cdot p_q$. However with a little confusion and for notation convenience we shall sometimes use the notation $p_1\cdot\cdot\cdot p_n$ to denote the knot $K_1\star K_2\star\cdot\cdot\cdot\star K_q$ and we say that this knot is related to the number $p_1\cdot\cdot\cdot p_n$ (as similar to the knot ${\bf3_1}$ which is related to the number $2$ but is assigned with the number $1$) and we keep in mind that the knot $K_1\star K_2\star\cdot\cdot\cdot\star K_q$ may not be assigned with the number $p_1\cdot\cdot\cdot p_n$. With this notation then we may say that the composite number $3\cdot5$ jumps over the number $2^4$ which means that the composite knot ${\bf4_1}\star{\bf5_1}$ jumps over the knot ${\bf3_1}\star{\bf3_1}\star{\bf3_1}\star{\bf3_1}$. $\diamond$

{\bf Definition (Jumping over of the general kind)}. At the $n$th step let a composite knot $K^{\prime}$ be related with a number $p_1\cdot p_2\cdot\cdot\cdot p_m$ where the number $p_1\cdot p_2\cdot\cdot\cdot p_m$ is in the $n$th step. Then we say that the knot $K^{\prime}$ (or the number $p_1\cdot p_2\cdot\cdot\cdot p_m$) is of jumping over of the general kind (with respect to the knot $K$ in the definition of the jumping over of the first kind and we also write $K\prec K^{\prime}$) if
$K$ satisfies one of the following conditions:

1) $K^{\prime}$ (or the number related to $K^{\prime}$) is of jumping over of the first kind; or

2) There exists a $p_i$ (for simplicity let it be $p_1$) and a prime number $q$ such that $p_1$ and $q$ are in the same step $k$ for some $k$ and $q$ is the largest prime number in this step such that the numbers $p_1\cdot p_2\cdot\cdot\cdot p_m$ and $q\cdot p_2\cdot\cdot\cdot p_m$ are also in the same step and that the knot $K_q^{\prime}$ related with $q\cdot p_2\cdot\cdot\cdot p_m$ is of jumping over of the first kind. $\diamond$

{\bf Remark}. The condition 2) is a natural generalization of 1) that if $K^{\prime}$ and the knot $K_q^{\prime}$ are as in 2) then they are both in the preordering sequences of an induction $n$th step or both not.
Then since $K_q^{\prime}$  is of jumping over into an $(n+1)$th induction step and thus is not in the preordering sequences of the induction $n$th step we have that $K^{\prime}$ is also of jumping over into this $(n+1)$th induction step (even if $K^{\prime}$ is not of jumping over of the first kind). This means that $K^{\prime}$ is of jumping over of the general kind. $\diamond$

{\bf Example of jumping over of the general kind}. At an induction step let $K^{\prime}$ be represented by $11\cdot5\cdot5$ (where we let $p_1=11$) and let $K_q^{\prime}$ be represented by $13\cdot5\cdot5$ (where we let $q=13$). Then $K_q^{\prime}$ is of jumping over of the first kind. Thus we have that $K^{\prime}$ is of jumping over (of the general kind). $\diamond$

We shall show that if $K={\bf 3_1}^n \prec K^{\prime}$ then we can set
$K={\bf 3_1}^n<K^{\prime}$. Thus we have, in the above first example, ${\bf
(3_1\star 3_1)\star(3_1\star 3_1)}< {\bf 4_1\star 5_1}$ while $2^4
> 3\cdot 5$.  From this
property we  shall have rooms for the introduction of the $\times$
operation such that composite numbers are assigned to composite knots and prime numbers are assigned to prime knots. We have the following theorem.

\begin{theorem}

If $K ={\bf 3_1}^n\prec K^{\prime}$ then it is consistent with the preordering property that $K={\bf 3_1}^n <K^{\prime}$ for setting up the table.
\end{theorem}

For proving this theorem let us first prove the following lemma.
\begin{lemma}

The preordering sequences for the construction of the $n$th step do not have knots of jumping over of the general kind.
\end{lemma}
{\bf Proof of the lemma}.
It is clear that the preordering sequence obtained by the ${\bf 3_1}$ taking a $\star$ operation with the previous $(n-1)$th step has no knots with the jumping over of the first kind property since ${\bf 3_1}$ is corresponded with the number $2$ and the  previous $(n-1)$th step has no knots with the jumping over of the first kind property for this $(n-1)$th step. Then preordering sequence obtained by the ${\bf 4_1}$ taking a $\star$ operation with the previous $(n-2)$th step has no knots with the jump over of the first kind property since ${\bf 4_1}$ is assigned with the number $3$ and $3<2^2$ and the  previous $(n-2)$th step has no knots with the jumping over  of the first kind property for this $(n-2)$th step. Continuing in this way we have that all the knots in these preordering sequences do not satisfy the property of jumping over of the first kind. Then let us show that these preordering sequences have no knots with the property of jumping over of the general kind. Suppose this is not true. Then there exists a knot with the property of jumping over of the general kind and let this knot be represented by a number of the form
$p_1\cdot p_2\cdot\cdot\cdot p_m$ as in the definition of jumping over of the general kind such that there exists a prime number $q$ and that $p_1$ and $q$ are in the same step $k$ for some $k$ and $q$ is the largest prime number in this step such that the numbers $p_1\cdot p_2\cdot\cdot\cdot p_m$ and $q\cdot p_2\cdot\cdot\cdot p_m$ are also in the same step and the knot $K_q$ represented by $q\cdot p_2\cdot\cdot\cdot p_m$ is of jumping over of the first kind. Then since $p_1$ and $q$ are in the same step $k$ we have that the two knots represented by $p_1\cdot p_2\cdot\cdot\cdot p_m$ and $q\cdot p_2\cdot\cdot\cdot p_m$ are elements of two preordering sequences for the construction of the same $n$th step. Now since we have shown that the preordering sequences for the construction of the $n$th step do not have knots of jumping over of the first kind we have that this is a contradiction. This proves the lemma.
$\diamond$

{\bf Proof of the theorem}.
By the above lemma if $K={\bf 3_1}^n\prec K^{\prime}$ then $ K^{\prime}$ is not in the preordering sequences for the $n$th step and thus is pushed out from the $n$th step by the preordering sequences for the $n$th step and thus we have $K={\bf 3_1}^n< K^{\prime}$, as was to be proved.
$\diamond$

{\bf Remark}. We remark that there may exist knots (or numbers related to the knots) which are not in the preordering sequences and are not of jumping over. An example of such special knot is the knot ${\bf4_1}\star{\bf5_1}\star{\bf5_1}$ related with $3\cdot 5\cdot 5$ (but is not assigned with this number).
$\diamond$

{\bf Definition}.
When there exists a knot which is not in the preordering sequences of the $n$th step and is not of jumping over we put this knot back into the $n$th step to join the  preordering sequences for the filling and covering of the $n$th step if there are rooms to be filled in the $n$th step after the filling with the preordering sequences. Let us call the preordering sequences together with the knots which are not in the preordering sequences of the $n$th step and are not of jumping over as the generalized preordering sequences (for the filling and covering of the $n$th step).
$\diamond$

{\bf Remark}. By using the generalized preordering sequences for the covering of the $n$th step we have that the knots (or the number related to the knots) in the $n$th step pushed out of the $n$th step by the generalized preordering sequences  are just the knots of jumping over (of the general kind). $\diamond$

Then we also have the following theorem.
\begin{theorem}

At each $n$th step ($n>3$) in the covering of the $n$th step ($n>3$) with the generalized preordering sequences there are rooms for introducing new knots with the $\times$ operations.
\label{times}
\end{theorem}
{\bf Proof}.
We want to show that at each $n$th step ($n>3$) there are rooms for
introducing new knots with the $\times$ operations.
At $n=4$ we have shown that there is the room at the position $9$ for introducing the knot ${\bf 3_1}\times{\bf 3_1}$ with the $\times$ operation.
Let us suppose that this property holds at an
induction step $n-1$. Let us then consider the induction step $n$.
For each $n$ because of the relation between $1$ and $2$ for ${\bf
3_1}$  as a part of the induction step $n$ the number $2^n$ is assigned to the knot
${\bf 3_1}^n$ which is a $\star$
product of $n$ ${\bf 3_1}$. Then we want to show that for this
induction step $n$ by using the $\prec$ property we have rooms for
introducing the $\times$ operation. Let $K^{\prime}$ be a knot
such that ${\bf 3_1}^{n-1}\prec K^{\prime}$
and $K^{\prime}=K_2\star K_3$
is as in the definition of $\prec$ of jumping over of the first kind
such that $p_1\cdot\cdot\cdot p_{n_2}\cdot q_1\cdot\cdot\cdot q_{n_3}<2^{n-1}$
(e.g. for $n-1=4$
we have $K^4={\bf 3_1}\star{\bf 3_1}\star{\bf 3_1}\star{\bf 3_1}$
and $K^{\prime}=K_2\star K_3={\bf 4_1}\star{\bf 5_1}$). Then let
us consider $K^{\prime\prime}=({\bf 3_1}\star K_2)\star K_3$.
Clearly we have ${\bf 3_1}^n\prec K^{\prime\prime}$. Thus for each
$K^{\prime}$  we have a $K^{\prime\prime}$ such that
${\bf 3_1}^n\prec K^{\prime\prime}$. Clearly all these $ K^{\prime\prime}$ are different.

Then from $K^{\prime}$ let us
construct more $K^{\prime\prime}$, as follows. Let $K^{\prime}$ be a knot of jumping over
of the first kind. Let
$p_1\cdot\cdot\cdot p_{n_2}$ and $q_1\cdot\cdot\cdot q_{n_3}$ be
as in the definition of jumping over
of the first kind. Then as in the definition of jumping over
of the first kind (w.l.o.g)
we let
\begin{equation}
2^{n_0}<p_1\cdot\cdot\cdot p_{n_2} \quad \mbox {and}\quad
2^{n_1}>q_1\cdot\cdot\cdot q_{n_3}
\label{forjumpingover}
\end{equation}

Then we have
\begin{equation}
2^{n_0+1}<(2\cdot p_1\cdot\cdot\cdot p_{n_2})-1\quad \mbox {and}\quad
2^{n_1}>q_1\cdot\cdot\cdot q_{n_3}
\label{jumpover12}
\end{equation}
Also it is trivial that we have
$2^{n_0}<(2\cdot p_1\cdot\cdot\cdot p_{n_2})-1$ and
$2^{n_1+1}>q_1\cdot\cdot\cdot q_{n_3}$. This shows that ${\bf 3_1}^n\prec
K^{\prime\prime}:= K_{2a}\star K_{3}$ where $K_{2a}$ denotes the
knot with the number $(2\cdot p_1\cdot\cdot\cdot p_{n_2})-1$ as in
the definition of jumping over of the first kind (We remark that this $K^{\prime\prime}$
corresponds to the knot ${\bf 4_1}\star({\bf 4_1}\star{\bf 4_1})$
in the above induction step where $K_{2a}={\bf 4_1}\star{\bf 4_1}$
is with the number $2\cdot 5-1=3\cdot 3$).

It is clear that all these more $K^{\prime\prime}$ are different from the above $K^{\prime\prime}$ constructed by the above method of taking a $\star$ operation with ${\bf 3_1}$.
Thus there are more $ K^{\prime\prime}$ than $K^{\prime}$. Thus at this $n$th step there are rooms for introducing new knots with the $\times$ operations.
This proves the theorem. $\diamond$

{\bf Remark}.
In the proof of the above theorem we have a way to construct the knots $ K^{\prime\prime}$ by replacing a number $a$ with the number $2a-1$. There is another way of constructing
the knots $ K^{\prime\prime}$ by replacing a number $b$ with the number $2b+1$. For this way we need to check that the number related to $ K^{\prime\prime}$ is in the $(n-1)$th step for $ K^{\prime\prime}$ of jumping over into the $n$th step.

 As an example let us consider the knot $ K^{\prime}={\bf 4_1}\star{\bf 4_1}\star{\bf 4_1}$ of jumping over into the $6$th step with the following data:
\begin{equation}
2^{3}<3\cdot3 \quad \mbox {and}\quad
2^{2}>3
\label{jumpingover9}
\end{equation}

From this data we have:
\begin{equation}
2^{3+1}<2\cdot3\cdot3-1 =17\quad \mbox {and}\quad
2^{2}>3
\label{jumpingover10a}
\end{equation}
This data gives a knot $ K^{\prime\prime}$ with the related number $3\cdot17 $.

On the other hand from the data (\ref{jumpingover9}) we have:
\begin{equation}
2^{3}<3\cdot3 \quad \mbox {and}\quad
2^{2+1}>2\cdot3+1
\label{jumpingover10}
\end{equation}
Since $(3\cdot3)(2\cdot3+1)=(2\cdot5-1)(2\cdot3+1)=2\cdot5\cdot2\cdot3+2\cdot2-1<2\cdot2\cdot2^4-1<2^6$ we have that the knot $ K^{\prime\prime}={\bf 4_1}\star{\bf 4_1}\star{\bf 5_2}$ related with the number $3\cdot3\cdot7$ is of jumping over into the $7$th step (We shall show that ${\bf 5_2}$
is assigned with the number $7$). $\diamond$

{\bf Remark}. The above theorem shows that at each $n$th step there are rooms for introducing new knots with the $\times$ operations and thus we may establish a one-to-one correspondence of knots and numbers such that prime knots are bijectively assigned with prime numbers. Further to this theorem we have the following main theorem:

\begin{theorem}
A classification table of knots can be formed (as partly described by the above table up to $2^n$ with $n=5$) by induction on the number $2^n$
such that knots are one-to-one assigned with an integer and prime knots are bijectively assigned with prime numbers such that the prime number $2$ corresponds to the trefoil knot. This assignment is onto the set of positive integers except $2$ where the trefoil knot is assigned with 1 and is related to $2$ and at each $n$th induction step of the number $2^n$ there are rooms for introducing new knots with the $\times$ operations only.

Further this assignment of knots to numbers for the $n$th induction step 
effectively includes the determination of the distribution of prime numbers in the $n$th induction step and is by induction determined by this assignment for the previous $n-1$ induction steps.

Furthermore  this assignment of knots to numbers for the $n$th induction step is by induction  
 formed by the assignments
(or structures) 
of the previous steps which are as the preordering sequences in the forming of the  assignment (or structure) of the $n$th induction step such that the preordering sequences are arranged consistently with the ordering  in the forming of the assignment (or structure) of the $n$th induction step. Thus  all the properties of the structures of the previous steps  
are extended as properties of the structure of the $n$th induction step.

\label{maintheorem}
\end{theorem}

{\bf Remark}. Let us also call this assignment of knots to numbers as the structure of numbers obtained by assigning numbers to knots. This structure of numbers is the original number system together with the one-to-one assignment of numbers to knots.
 
{\bf Proof}.
By the above lemmas and theorems we have that
the generalized preordering sequences 
have the function of pushing out those composite knots 
of jumping over from the $n$th step.
It follows that for step $n>3$ there must exist
 chains of transitions whose initial states are composite knots in repeat (to be replaced by the new composite knots with the $\times$ operations only); or the knots  
 jumping over into this $n$th step from the previous $(n-1)$th step or the knots in the preordering sequences with the $\times$ operations; such that the composite knots 
 of jumping over are pushed out from the $n$th step by these chains. These chains are obtained by ordering the subsequences of the preordering sequences
 such that the preordering property holds in the 
 $n$th step. 
 Further 
the intermediate states of the chains must be positions of composite numbers. This is because that if a chain is transited to an intermediate state which is a position of prime number then there are no composite knots related with  
this prime number and thus this chain can not be transited to the next state and  is stayed at the intermediate state forever and thus the chain can not push out the composite knot of jumping over. Then when a composite knot is transited to the position of an intermediate state (which is a position of composite number as has just been proved) this knot is definitely assigned with this composite number.
Then when a composite knot which is in repeat is transited to the position of an intermediate state this knot is also definitely assigned with this composite number. 
It follows that when the chains are completed we have that the ordering of the subsequences of preordering sequences is determined.

Then the remaining knots (which are not at the transition states of the chains) which are not in repeat are definitely assigned with the number of the position of these knots in the $n$th step. For these knots the numbers of positions assigned to them are just the number related to them respectively.

Then the remaining knots (which are not at the transition states of the chains) which are in repeat must be replaced by new prime knots because of the repeat and that no other composite knots related with numbers in this $n$th step in the generalized preordering sequences can be used to replace the remaining knots.  
This means that the numbers of the positions of these remaining knots in repeat are prime numbers in this $n$th step. This is because that if the number of the position assigned to the new prime knot is a composite number then the composite knot related  
with this composite number is either in a transition state or is not in transition.
If the composite knot is not in transition then the composite number related to 
 this composite knot is just the number assigning to this composite knot and since this number is  also assigned to the new prime knot that this is a contradiction. Then if this composite knot is in transition state then this means that the remaining knot is also in transition state and this is a contradiction since by definition the remaining knot is not at the transition states of the chains. 
 
 Thus prime numbers in the $n$th step are assigned and  are only assigned to prime knots which replace the remaining knots in repeat 
 in the $n$th step. 
 Thus from the preordering sequences we have determined the positions (i.e. the distribution) of prime numbers in the $n$th step. Now since the preordering sequences are constructed by  the previous steps  we have shown that the basic structure (in the sense of above proof) of this assignment of knots with numbers for the $n$th step (including the determination of the distribution of prime numbers in the $n$th step) is determined by this assignment of knots with numbers for the previous $n-1$ steps. In other words we have that the basic structure of the $n$th induction step  
 is determined by the structure of the previous $n-1$ steps.

In summary, we have shown that  the basic structure of the $n$th induction step is formed by the structures of the previous steps which are as the preordering sequences in the forming of the structure of the $n$th induction step such that the preordering sequences are arranged consistently with the ordering  in the forming of the structure of the $n$th induction step.
 
 To complete the proof of this theorem 
 let us  show that at each $n$th induction step ($n>3$) there are rooms for introducing new composite knots with the $\times$ operations only and we can determine the ordering of these composite knots with the $\times$ operations only in each $n$th induction step. 
 
In the above proof we have shown that the basic structure of the $n$th induction step is determined by the structure of the previous steps such that the positions of the composite knots with the $\times$ operations only in the $n$th induction step are determined by the structures of the previous steps. Then by the preordering of these composite knots with the $\times$ operations only (with the $\star$ operation replaced by the $\times$ operation),
these positions 
can be fitted for the correct composite knots with the $\times$ operations only constructed (by the $\times$ operations) by knots in the previous steps, just as the preordering of those composite knots with the $\star$ operations  in the previous steps. 
Thus for this $n$th induction step the introducing and ordering of composite knots with the $\times$ operations only is also determined by the structures of the previous $n-1$ steps.
 
Since  the basic structure of the $n$th induction step is formed by the structures of the previous steps which are as the preordering sequences in the forming of the structure of the $n$th induction step such that the preordering sequences are arranged consistently with the ordering  in the forming of the structure of the $n$th induction step, and that the structure of an induction step consists of the properties of the structure of this induction step,
 we have that all the properties of the structures of the previous steps  
are extended as properties of the structure of the $n$th induction step.


 Thus the  ordering properties  
of composite knots with the $\times$ operations only  in the previous steps 
are extended as  ordering properties of  the composite knots with the $\times$ operations only in the $n$th induction step.
(These ordering properties of the composite knots with the $\times$ operations only can be used to find out the correct composite knots with the $\times$ operations only to be assigned at the correct positions in the $n$th step).

With this fact let  
us then show that at each $n$th induction step ($n>3$) there are rooms for introducing new composite knots with the $\times$ operations only.
 As in the proof of the theorem \ref{times} we first construct more $ K^{\prime\prime}$ by the method following (\ref{forjumpingover}). Let us start at the step $n=4$. For this step we have the knot $K^{\prime}={\bf 4_1}\star{\bf 5_1}$ jumps over into the step $n=5$. For this $K^{\prime}$ we have the following data as in (\ref{forjumpingover}):
\begin{equation}
2^{2}<5 \quad \mbox {and}\quad 
2^{2}>3
\label{jumpingover3}
\end{equation}
From (\ref{jumpingover3}) we construct a $ K^{\prime\prime}$ for the step $n=5$ by the following data:
\begin{equation}
2^{2+1}<2\cdot5-1= 3\cdot3 \quad \mbox {and}\quad 
2^{2}>3
\label{jumpingover4}
\end{equation}
This data gives one more $ K^{\prime\prime}={\bf 4_1}\star{\bf 4_1}\star{\bf 4_1}$. Then from (\ref{jumpingover3}) we construct one more $ K^{\prime\prime}$ for the step $n=5$ by the following data:
\begin{equation}
2^{3}>5 \quad \mbox {and}\quad 
2^{1+1}<2\cdot3-1=5
\label{jumpingover5}
\end{equation}
This data gives one more $K^{\prime\prime}={\bf 5_1}\star{\bf 5_1}$. Thus in this step $n=5$ there are two rooms  for the two knots $ K^{\prime}={\bf 4_1}\star{\bf 5_1}$ and ${\bf 3_1}\star({\bf 3_1}\times{\bf 3_1})$ coming from the preordering sequences  and there exists exactly one room for introducing a  new composite knot with the $\times$ operations only (Recall that we also have a $K^{\prime\prime}={\bf 3_1}\star{\bf 4_1}\star{\bf 5_1}$). 
From the ordering of knots in the previous steps we determine that ${\bf 3_1}\times{\bf 4_1}$ is the composite knot with the $\times$ operations only for this step.
Thus at the $4$th and $5$th steps we can and only can introduce exactly one composite knot with the $\times$ operations only and they are the knots ${\bf 3_1}\times {\bf 3_1}$ and ${\bf 3_1}\times {\bf 4_1}$ respectively.
This shows that at the $4$th and the $5$th steps we can determine the number of prime knots with the minimal number of crossings $=3$ and $=4$ respectively (These two prime knots are denoted by ${\bf 3_1}$ and ${\bf 4_1}$ respectively and we do not distinguish knots with their mirror images for this determination of the ordering of knots with the $\times$ operations only. This also shows that there are rooms for introducing new composite knots with the $\times$ operations only in the $4$th and $5$th steps).

Then since this property is 
extended as a property
in the $6$th step  
we can thus determine that the $6$th step is a step for introducing new composite knots with the $\times$ operations only of the form ${\bf 3_1}\times {\bf 5_{(\cdot)}}$ where ${\bf 5_{(\cdot)}}$ denotes a prime knot with the minimal number of crossings $=5$ (and thus there are rooms for introducing new composite knots with the $\times$ operations only in this $6$th step). Also since the properties in the $ 4$th and $5$th steps are
extended as properties
in the $6$th step we can determine
the number of prime knots with the minimal number of crossings $=5$ by the knots of the form ${\bf 3_1}\times {\bf 5_{(\cdot)}}$ as
this is a property of knots with the $\times$ operations only in the $4$th and $5$th steps (In the classification table  in the next section we show that there are exactly two composite knots of the form ${\bf 3_1}\times {\bf 5_1}$ and ${\bf 3_1}\times {\bf 5_2}$ in the $6$th step whose ordering are determined by the preordering property of knots and the structure of the $6$th step. This thus shows that there are exactly two prime knots with
the minimal number of crossings $=5$ and they are denoted by ${\bf 5_1}$ and ${\bf 5_2}$ respectively).

Then since the properties of the $4$th, $5$th and $6$th steps are  
extended as properties
in the  $7$th step we can determine that  the $7$th step is a step for introducing new composite knots with the $\times$ operations only of the form ${\bf 3_1}\times {\bf 6_{(\cdot)}}$ where ${\bf 6_{(\cdot)}}$ denotes a prime knot with the minimal number of crossings $=6$ (and thus there are rooms for introducing new composite knots with the $\times$ operations only in this $7$th step). Also since the properties in the $4$th, $5$th and $6$th steps are  
extended as properties
in the $7$th step we can determine
the number of prime knots with the minimal number of crossings $=6$ by the knots of the form ${\bf 3_1}\times {\bf 6_{(\cdot)}}$ as 
this is a property of knots with the $\times$ operations only in the $4$th, $5$th and $6$th steps (In the classification table  in the next section we show that there are exactly three composite knots of the form ${\bf 3_1}\times {\bf 6_1}$, ${\bf 3_1}\times {\bf 6_2}$ and ${\bf 3_1}\times {\bf 6_3}$ in the $7$th step whose ordering are determined by the preordering property of knots and the structure of the $7$th step. This thus shows that there are exactly three prime knots with
the minimal number of crossings $=6$ and they are denoted by ${\bf 6_1}$, ${\bf 6_2}$ and ${\bf 6_3}$ respectively).

Continuing in this way we thus have that at each $n$th induction step $(n>3)$ we can determine the number of prime knots with the minimal number of crossings $=n-1$ and there are rooms for introducing new composite knots with the $\times$ operations only. This proves the theorem. $\diamond$

{\bf Example}. Let us consider the above classification table up to $2^5$ (with $n$ up to $5$) as an example.
For the induction step at $n=2$ (or at $2^2$) we have one preordering sequence obtained by letting ${\bf 3_1}$ to take a $\star$ operation with the step $n=1$ (For the step $n=1$ the number $2^1$ is related to the trefoil knot ${\bf 3_1}$): ${\bf 3_1}\star {\bf 3_1}$. Then we fill the step $n=2$ with this preordering sequence and we have the following ordering of knots for this step $n=2$:
\begin{equation}
{\bf 3_1}\star {\bf 3_1},
{\bf 3_1}\star  {\bf 3_1}
\label{step2}
\end{equation}
where the first ${\bf 3_1}\star  {\bf 3_1}$ placed at the position $3$ is the preordering sequence while the second ${\bf 3_1}\star  {\bf 3_1}$ placed at the position $2^2$ is required by the induction procedure. For this step there is no numbers of jumping over.  Then we have that the first ${\bf 3_1}\star  {\bf 3_1}$ is a repeat of the second ${\bf 3_1}\star  {\bf 3_1}$.
Thus this repeat one must be replaced by a new prime knot. Let us choose the prime knot ${\bf 4_1}$ to be this new prime knot since ${\bf 4_1}$ is the smallest of prime knots other than the trefoil knot. Then this new prime knot must be at the position of a prime number, as we have proved in the above theorem. Thus we have determined that $3$ is a prime number in this step $n=2$ by using the structure of numbers of step $n=1$ which is only with the prime number $2$.

Then for the induction step at $n=3$ (or at $2^3$) we have two preordering sequence obtained by letting ${\bf 4_1}$ to take a $\star$ operation with the step $n=1$ and by letting ${\bf 3_1}$ to take a $\star$ operation with the step $n=2$:
\begin{equation}
{\bf 4_1}\star{\bf 3_1}; {\bf 3_1}\star{\bf 4_1},
{\bf 3_1}\star ({\bf 3_1}\star {\bf 3_1})
\label{step3a}
\end{equation}
where the first knot is the preordering sequence obtained by letting ${\bf 4_1}$ to take a $\star$ operation with the step $n=1$ and the second and third knots is the preordering sequence obtained by letting ${\bf 3_1}$ to take a $\star$ operation with the step $n=2$.

For this step there is no numbers of jumping over and thus there are no chains of transition. Thus the ordering of the above three knots in this step follow the usual ordering of numbers.
Thus the number assigned to the knot ${\bf 4_1}\star {\bf 3_1}={\bf 3_1}\star{\bf 4_1}$ must be assigned with a number less than that of ${\bf 3_1}\star {\bf 3_1}\star {\bf 3_1}$ by the ordering of ${\bf 3_1}\star{\bf 4_1}$ and ${\bf 3_1}\star {\bf 3_1}\star {\bf 3_1}$ in the second preordering sequence.
By this ordering of the two preordering sequences  we have that the step $n=3$ is of the following form:
\begin{equation}
{\bf 4_1}\star{\bf 3_1};{\bf 3_1}\star{\bf 4_1},
{\bf 3_1}\star ({\bf 3_1}\star {\bf 3_1});
{\bf 3_1}\star {\bf 3_1}\star {\bf 3_1}
\label{cover}
\end{equation}
where the fourth knot ${\bf 3_1}\star {\bf 3_1}\star {\bf 3_1}$ is put at the position of $2^3$ and is assigned with the number $2^3$ as required by the induction procedure.
 Thus the third knot ${\bf 3_1}\star ({\bf 3_1}\star {\bf 3_1})$  is a repeated one and thus must be replaced by a prime knot and the position of this prime knot is determined to be a prime number.  Thus we have determined that the number $7$ is a prime number.
 Then since there are no chains of transition we have that the composite knot ${\bf 3_1}\star{\bf 4_1}$ must be assigned with the number related to this knot and this number is $2\cdot 3=6$. Thus the composite knot ${\bf 3_1}\star{\bf 4_1}$ is at the position of $6$ and that
the first knot ${\bf 4_1}\star{\bf 3_1}$ is a repeat of the second knot and thus must be replaced by a prime knot. Then since this prime knot is at the position of $5$ we have that $5$ is determined to be a prime number. Now the two prime knots at $5$ and $7$ must be the prime knots ${\bf 5_1}$ and ${\bf 5_2}$ respectively since these two knots are  the smallest prime knots other than ${\bf 3_1}$ and ${\bf 4_1}$ (We may just put in two prime knots and then later determine what these two knots will be. If we put in other prime knots then this will not change the distribution of prime numbers determined by the structure of numbers of the previous steps
and it is only that the prime knots are assigned with incorrect prime numbers. Further as shown in the above proof by using knots of the form ${\bf 3_1}\times{\bf 5_{(\cdot)}}$  we can determine that there are exactly two prime knots with minimal number of crossings $=5$ and they are denoted by ${\bf 5_1}$ and ${\bf 5_2}$ respectively. From this we can then determine that these two prime knots are ${\bf 5_1}$ and ${\bf 5_2}$).
 Thus we have
the following ordering for $n=3$:
\begin{equation}
{\bf 5_1}<{\bf 3_1}\star{\bf 4_1}<
{\bf 5_2}<{\bf 3_1}\star {\bf 3_1}\star {\bf 3_1}
\label{order1}
\end{equation}
where ${\bf 5_1}$ is assigned with the prime number $5$ and ${\bf 5_2}$ is assigned with the prime number $7$. This gives the induction step $n=3$. For this step there is no knot with $\times$ operation since there is no knots of jumping over.

Let us then consider the step $n=4$ (or $2^4$). For this step we have the following three preordering sequences obtained from the steps $n=1,2,3$:
\begin{equation}
\begin{array}{rl}
&{\bf 5_1}\star {\bf3_1};\\
&{\bf 4_1}\star {\bf 4_1},{\bf 4_1}\star {\bf 3_1}\star {\bf 3_1};\\
&{\bf 3_1}\star {\bf5_1},{\bf 3_1}\star {\bf 3_1}\star {\bf 4_1},{\bf 3_1}\star {\bf 5_2},{\bf 3_1}\star {\bf 3_1}\star {\bf 3_1}\star{\bf 3_1};\\
 \end{array}
\label{order2}
\end{equation}
where the third sequence is obtained by taking $\star$ operation of the knot ${\bf 3_1}$ with step $n=3$ while the third sequence is obtained by taking $\star$ operation of the knot ${\bf 4_1}$ with the step $n=2$ and the first sequence is obtained by taking $\star$ operation of the knot ${\bf 5_1}$ with step $n=1$. Then as required by the induction procedure the knot ${\bf 3_1}\star {\bf 3_1}\star {\bf 3_1}\star {\bf 3_1}$ is assigned at the position of $2^4$. The total number of knots in (\ref{order2}) plus this knot is exactly $2^3$ which is the total number of this step $n=4$.

{\bf Remark}. We  have one more preordering sequence  obtained by taking $\star$ operation of the knot ${\bf 5_2}$ with step $n=1$. This preordering sequence gives the knot
${\bf 5_1}\star {\bf3_1}$. However since the knots in (\ref{order2}) and the knot ${\bf 3_1}\star {\bf 3_1}\star {\bf 3_1}\star {\bf 3_1}$ assigned at the position of $2^4$ are enough for covering this step $n=4$ and that the knot ${\bf 5_1}\star {\bf3_1}$ of this preordering sequence is a repeat of the knot ${\bf 5_1}\star {\bf3_1}$ in (\ref{order2}) that this preordering sequence obtained by taking $\star$ operation of the knot ${\bf 5_2}$ with step $n=1$ can be omitted. $\diamond$

Then to find the chains of transition for this step let us order the three preordering sequences with the following ordering where we rewrite the preordering sequences in column form and the knot ${\bf 3_1}\star {\bf 3_1}\star {\bf 3_1}\star {\bf 3_1}$ assigned at the position of $2^4$ is put to follow the three sequences:
\begin{equation}
\begin{array}{rl}
&{\bf 5_1}\star {\bf3_1};\\
&{\bf 3_1}\star {\bf5_1},\\
&{\bf 3_1}\star {\bf 3_1}\star {\bf 4_1},\\
&{\bf 3_1}\star {\bf 5_2},\\
&{\bf 3_1}\star {\bf 3_1}\star {\bf 3_1}\star{\bf 3_1};\\
&{\bf 4_1}\star {\bf 4_1},\\
&{\bf 4_1}\star {\bf 3_1}\star {\bf 3_1};\\
&{\bf 3_1}\star {\bf 3_1}\star {\bf 3_1}\star{\bf 3_1}
 \end{array}
\label{step4}
\end{equation}
We notice that this column exactly fills the step $n=4$.

For this step we have that the number $3\cdot5$ (or the knot ${\bf 4_1}\star {\bf 5_1}$ related with $3\cdot5$ ) is of jumping over.
From (\ref{step4}) we have the following chain of transition for pushing out ${\bf 4_1}\star {\bf 5_1}$ at $3\cdot5$ by a knot with the $\times$ operation replacing the repeated knot ${\bf 5_1}\star {\bf 3_1}$ at the position of $9=3\cdot3$:
\begin{equation}
{\bf 3_1}\times {\bf 3_1} (\mbox {at} 3\cdot3)\to {\bf 4_1}\star {\bf 4_1}(\mbox {at} 2\cdot7) \to{\bf 3_1}\star {\bf 5_2} (\mbox {at} 2\cdot2\cdot3)\to
{\bf 3_1}\star {\bf 3_1}\star {\bf 4_1} (\mbox {at} 3\cdot5)\to {\bf 4_1}\star {\bf 5_1} (\mbox {pushed out})
\label{Chain2}
\end{equation}
where we choose the knot ${\bf 3_1}\times{\bf 3_1}$ as the knot with the $\times$ operation since ${\bf 3_1}\times{\bf 3_1}$ is the smallest one of such knots.
For this chain the intermediate states are at positions of composite numbers $2\cdot7$, $2\cdot2\cdot3$ and $3\cdot5$. Thus the knots in this chain at the positions of these composite numbers are assigned with these composite numbers respectively.

Then once this chain of pushing out ${\bf 4_1}\star {\bf 5_1}$ at $3\cdot5$ is set up we have that the other knots in repeat must by replaced by prime knots and that their positions must be prime numbers.
These positions are at $11$ and $13$ and thus $11$ and $13$ are determined to be prime numbers (The knot ${\bf 3_1}\star{\bf 3_1}\star{\bf 3_1}\star{\bf 3_1}$ at the end of this step must be assigned with $2^4=16$ by the induction procedure and thus the knot at $13$ is a repeat).
Then  the new prime knots ${\bf 6_1}$ and ${\bf 6_2}$ are suitable knots corresponding to the prime numbers $11$ and $13$ respectively since they are the smallest  prime knots other than ${\bf 3_1}$, ${\bf 4_1}$, ${\bf 5_1}$ and ${\bf 5_2}$ (As the above induction step we may just put in two prime knots and then later determine what these two prime knots will be. As shown in the above proof by using knots of the form ${\bf 3_1}\times{\bf 6_{(\cdot)}}$  we can determine that there are exactly three prime knots with minimal number of crossings $=6$ and they are denoted by ${\bf 6_1}$, ${\bf 6_2}$ and ${\bf 6_3}$ respectively. From this we can then determine that these two prime knots are ${\bf 6_1}$ and ${\bf 6_2}$).

 This completes the step $n=4$.
Thus the structure of numbers
of this step (including distribution of prime numbers in this step) is  determined by the structure of numbers of the previous induction steps.

Let us then consider the step $n=5$. For this step we have the following four preordering sequences from the previous steps $n=1,2,3,4$:
\begin{equation}
{\bf 6_1}\star{\bf 3_1}
\label{step5a}
\end{equation}
and
\begin{equation}
\begin{array}{rl}
 & {\bf 5_2}\star{\bf 4_1},\\
 & {\bf 5_2}\star ({\bf 3_1}\star {\bf 3_1})
\end{array}
\label{step5b}
\end{equation}
and
\begin{equation}
\begin{array}{rl}
 & {\bf 4_1}\star{\bf 5_1},\\
 & {\bf 4_1}\star ({\bf 3_1}\star{\bf 4_1}),\\
 & {\bf 4_1}\star {\bf 5_2},\\
 & {\bf 4_1}\star ({\bf 3_1}\star {\bf 3_1}\star {\bf 3_1})
 \end{array}
\label{CC2}
\end{equation}
and
\begin{equation}
\begin{array}{rl}
 & {\bf 3_1}\star ({\bf 3_1}\times {\bf 3_1}),\\
 & {\bf 3_1}\star({\bf 3_1}\star {\bf 5_1}),\\
 & {\bf 3_1}\star {\bf 6_1},\\
 & {\bf 3_1}\star ({\bf 3_1}\star{\bf 5_2}),\\
 & {\bf 3_1}\star{\bf 6_2}, \\
 & {\bf 3_1}\star ({\bf 4_1}\star {\bf 4_1}),\\
 & {\bf 3_1}\star ({\bf 3_1}\star {\bf 3_1}\star {\bf 4_1}),\\
 & {\bf 3_1}\star ({\bf 3_1}\star {\bf 3_1}\star {\bf 3_1}\star {\bf 3_1})
\end{array}
\label{CC4}
\end{equation}

The total number of knots (including repeat) in the above sequences plus the knot ${\bf 3_1}\star {\bf 3_1}\star {\bf 3_1}\star {\bf 3_1}\star {\bf 3_1}$ to be assigned at the position of $2^5$ exactly cover this $n=5$ step.

{\bf Remark}. As similar to the step $n=4$ two preordering sequences ${\bf 5_1}\star {\bf 4_1}, {\bf 5_1}\star {\bf 3_1}\star {\bf 3_1}$ and ${\bf 6_2}\star{\bf 3_1}$ are omitted since these sequences are with knots which are repeats of the knots in the above preordering sequences. $\diamond$

Then to find the chains of transition for this step let us order these four preordering sequences with the following ordering where the knot ${\bf 3_1}\star {\bf 3_1}\star {\bf 3_1}\star {\bf 3_1}\star {\bf 3_1}$ assigned at the position of $2^5$ is put to follow the four sequences:
\begin{equation}
\begin{array}{rl}
& {\bf 6_1}\star {\bf 3_1}; \\
& {\bf 5_2}\star{\bf 4_1},\\
& {\bf 5_2}\star {\bf 3_1}\star {\bf 3_1};\\
& {\bf 4_1}\star{\bf 5_1}, \\
& {\bf 4_1}\star ({\bf 3_1}\star{\bf 4_1}),\\
& {\bf 4_1}\star {\bf 5_2}, \\
& {\bf 4_1}\star ({\bf 3_1}\star {\bf 3_1}\star {\bf 3_1});\\
& {\bf 3_1}\star ({\bf 3_1}\times {\bf 3_1}), \\
& {\bf 3_1}\star({\bf 3_1}\star {\bf 5_1}),\\
& {\bf 3_1}\star {\bf 6_1}, \\
& {\bf 3_1}\star ({\bf 3_1}\star{\bf 5_2}),\\
& {\bf 3_1}\star{\bf 6_2}, \\
& {\bf 3_1}\star ({\bf 4_1}\star {\bf 4_1}), \\
& {\bf 3_1}\star ({\bf 3_1}\star {\bf 3_1}\star {\bf 4_1}),\\
& {\bf 3_1}\star ({\bf 3_1}\star {\bf 3_1}\star {\bf 3_1}\star {\bf 3_1});\\
& ({\bf 3_1}\star {\bf 3_1})\star {\bf 3_1}\star {\bf 3_1}\star {\bf 3_1}
\end{array}
\label{CC5}
\end{equation}

For this step we have three composite knots ${\bf 3_1}\star ({\bf 4_1}\star{\bf 5_1})$, ${\bf 5_1}\star {\bf 5_1}$ and ${\bf 4_1}\star ({\bf 4_1}\star{\bf 4_1})$ (related with $2\cdot3\cdot5$,$5\cdot5$ and $3\cdot3\cdot3$ respectively) of jumping over and there are two new knots ${\bf 4_1}\star {\bf 5_1}$ and ${\bf 3_1}\star ({\bf 3_1}\times{\bf 3_1})$ coming from the previous step. Thus there is a room for the introduction of new knot obtained only by the $\times$ operation. Then this new knot must be the composite knot ${\bf 3_1}\times {\bf 4_1}$ since besides the composite knot ${\bf 3_1}\times {\bf 3_1}$ it is the smallest of composite knots of this kind.

From (\ref{CC5})
there is a chain of transition given by $18\to 21\to 22\to 26\to 28 \to 27$ and the composite knot ${\bf 4_1}\star ({\bf 4_1}\star{\bf 4_1})$ related with $27=3\cdot3\cdot3$ is pushed out into the next step by the composite knot ${\bf 5_2}\star{\bf 4_1}$ at the starting position $18$. Then this repeated knot must be replaced by a new composite knot obtained by the $\times$ operation only and this new composite knot must be the knot ${\bf 3_1}\times {\bf 4_1}$.

Then the composite knots at the intermediate states are assigned with the numbers of these states respectively.

In addition to the above chain there are two more chains: $24\to 30$ and $20\to 25$.
The
chain $24\to 30$ starts from ${\bf 3_1}\star ({\bf 3_1}\times{\bf 3_1})$ at $24$ and the composite knot ${\bf 3_1}\star ({\bf 4_1}\star{\bf 5_1})$ at $30$ is pushed out by the composite knot ${\bf 3_1}\star ({\bf 3_1}\star{\bf 3_1}\star{\bf 4_1})$.
Then the chain $20\to 25$ starts from ${\bf 4_1}\star {\bf 5_1}$ at $20$ and  the composite knot
${\bf 5_1}\star {\bf 5_1}$ at $25$ is pushed out by the composite knot ${\bf 3_1}\star ({\bf 3_1}\star{\bf 5_1})$.

Then the knots ${\bf 3_1}\star ({\bf 3_1}\star{\bf 3_1}\star{\bf 4_1})$ and ${\bf 3_1}\star ({\bf 3_1}\star{\bf 5_1})$ at the intemediate states of these two chains are assigned with the numbers $30=2\cdot3\cdot5$ and $25=5\cdot5$ respectively.

Now the remaining repeated composite knots at the positions $17,19, 23,29,31$ must be replaced by new prime knots  and thus $17,19, 23,29,31$ are determined to be prime numbers and they are determined by the prime numbers in the previous induction steps. Then we may follow the usual table of knots to determine that the new prime knots for the prime numbers $17,19, 23,29,31$ are ${\bf 6_3}$, ${\bf 7_1}$, ${\bf 7_2}$, ${\bf 7_3}$ and ${\bf 7_4}$ respectively (As the above induction steps we may just put in five prime knots and then later determine what these five prime knots will be. As shown in the above proof by using knots of the form ${\bf 3_1}\times{\bf 7_{(\cdot)}}$  we can determine the number of prime knots with minimal number of crossings $=7$. From this we can then determine these five prime knots).

In summary we have the following form of the step $n=5$:
\begin{equation}
\begin{array}{rl}
& {\bf 6_3} \\
& {\bf 3_1}\times{\bf 4_1}\\
& {\bf 7_1}\\
& {\bf 4_1}\star{\bf 5_1} \\
& {\bf 4_1}\star ({\bf 3_1}\star{\bf 4_1})\\
& {\bf 4_1}\star {\bf 5_2} \\
& {\bf 7_2}\\
& {\bf 3_1}\star ({\bf 3_1}\times {\bf 3_1}) \\
& {\bf 3_1}\star({\bf 3_1}\star {\bf 5_1})\\
& {\bf 3_1}\star {\bf 6_1} \\
& {\bf 3_1}\star ({\bf 3_1}\star{\bf 5_2})\\
& {\bf 3_1}\star{\bf 6_2} \\
& {\bf 7_3} \\
& {\bf 3_1}\star ({\bf 3_1}\star {\bf 3_1}\star {\bf 4_1})\\
& {\bf 7_4}\\
& ({\bf 3_1}\star {\bf 3_1})\star {\bf 3_1}\star {\bf 3_1}\star {\bf 3_1}
\end{array}
\label{CC6}
\end{equation}

This completes the induction step at $n=5$. We have that the structure of numbers of this step (including distribution of prime numbers in this step) is determined by the structure of numbers of the previous induction steps.
$\diamond$

\section{A Classification Table of Knots II}\label{sec2a}

Following the above classification table up to $2^5$ let us in this section give the table up to $2^7$.
Again we shall see from the table that the preordering property is clear. At the $7$th step  there is a special composite knot ${\bf 4_1}\star{\bf 5_1}\star{\bf 5_1}$ which is not of jumping over and is not in the preordering sequences (On the other hand the knot ${\bf 5_1}\star{\bf 5_1}\star{\bf 5_1}$ is of jumping over).

We remark again that it is interesting that (by the ordering of composite knots with the $\times$ operation only) at the $6$th step we require exactly two prime knots with minimal number of crossings $=5$ to form the two composite knots obtained by the $\times$ operation only. From this we can determine the number of prime knots with minimal number of crossings $=5$ without using the actual contruction of these prime knots. We then denote these two prime knots by ${\bf 5_1}$ and ${\bf 5_2}$ respectively and the two composite knots obtained by the $\times$ operation only by ${\bf 3_1}\times{\bf 5_1}$ and ${\bf 3_1}\times{\bf 5_2}$ respectively. Similarly at the $7$th step we can determine that there are exactly three prime knots with minimal number of crossings $=6$ and we denote these three prime knots by ${\bf 6_1}$ and ${\bf 6_2}$ and ${\bf 6_3}$ respectively. These three prime knots give the composite knots ${\bf 3_1}\times{\bf 6_1}$, ${\bf 3_1}\times{\bf 6_2}$ and ${\bf 3_1}\times{\bf 6_3}$ respectively. We can then expect that at the next $8$th step we may determine that the number of prime knots with minimal number of crossings $=7$ is $7$ and then at the next $9$th step the number of prime knots with minimal number of crossings $=8$ is $21$, and so on; as we know from the well known table of prime knots \cite{Rol}.
Here the point is that we can determine the number of prime knots with the same minimal number of crossings without using the actual construction of these prime knots (and by using only the classification table of knots).
\begin{displaymath}
\begin{array}{|c|c|c|} \hline
\mbox{Type of Knot}& \mbox{Assigned number} \,\, |m|
 &\mbox{Repeated Knots being replaced}
\\ \hline

{\bf 3_1\star 6_3} & 33 & {\bf } \\ \hline

{\bf 3_1\star(3_1\times 4_1)} & 34 & {\bf } \\ \hline

{\bf 3_1\star 7_1} & 35 & {\bf } \\ \hline

{\bf 3_1\times 5_1} & 36 & {\bf 3_1\star(4_1\star 5_1)} \\ \hline

{\bf 7_5} & 37 & {\bf 3_1\star(4_1\star 3_1\star 4_1)} \\ \hline

{\bf 3_1\times 5_2} & 38 & {\bf 3_1\star(4_1\star 5_2)} \\ \hline

{\bf 3_1\star 7_2} & 39 & {\bf } \\ \hline

{\bf 3_1\star(3_1\star 3_1\times 3_1)} &40  & {\bf } \\ \hline

{\bf 7_6} & 41 & {\bf 3_1\star(3_1\star 3_1\star 3_1\star 5_1)} \\ \hline

{\bf 5_1\star 5_1} & 42 & {\bf } \\ \hline

{\bf 7_7} & 43 & {\bf 5_1\star (3_1\star4_1)} \\ \hline

{\bf 5_1\star 5_2} & 44 & {\bf } \\ \hline

{\bf 4_1\times 4_1} & 45 & {\bf 5_1\star(3_1\star3_1\star3_1), 5_2\star(3_1\star4_1)} \\ \hline

{\bf 5_2\star5_2} & 46 & {\bf } \\ \hline

{\bf 8_1} & 47 & {\bf 5_2\star(3_1\star3_1\star3_1)},{\bf 4_1\star(3_1\times3_1)}\\ \hline

{\bf 4_1\star(3_1\star5_1)} & 48 & {\bf } \\ \hline

{\bf 4_1\star6_1} & 49 & {\bf } \\ \hline

{\bf 4_1\star(3_1\star5_2)} & 50 & {\bf } \\ \hline

{\bf 4_1\star6_2} & 51 & {\bf } \\ \hline

{\bf 4_1\star(4_1\star4_1)} & 52 & {\bf } \\ \hline

{\bf 8_2} & 53 & {\bf 4_1\star(4_1\star3_1\star3_1)} \\ \hline

{\bf 3_1\star(3_1\star3_1\star5_1)} & 54 & {\bf } \\ \hline

{\bf 3_1\star(3_1\star6_1)} & 55 & {\bf } \\ \hline

{\bf 3_1\star(3_1\star3_1\star5_2)} & 56 & {\bf } \\ \hline

{\bf 3_1\star(3_1\star6_2)} & 57 & {\bf } \\ \hline

{\bf 3_1\star7_3} & 58 & {\bf } \\ \hline

{\bf 8_3} & 59 & {\bf 3_1\star(3_1\star3_1\star3_1\star4_1)} \\ \hline

{\bf 3_1\star7_4} &60  & {\bf } \\ \hline

{\bf 8_4} &61  & {\bf 3_1\star(3_1\star3_1\star3_1\star 3_1\star3_1) } \\ \hline

{\bf 4_1\star(4_1\star3_1\star3_1)} & 62 & {\bf } \\ \hline

{\bf 4_1\star(3_1\star3_1\star3_1\star3_1)} & 63 & {\bf } \\ \hline

{\bf 3_1\star(3_1\star3_1\star3_1\star3_1\star3_1)} &64  & {\bf } \\ \hline
\end{array}
\end{displaymath}

\begin{displaymath}
\begin{array}{|c|c|c|} \hline
\mbox{Type of Knot}& \mbox{Assigned number} \,\, |m|
 &\mbox{Repeated Knots being replaced}
\\ \hline
{\bf 3_1\star(3_1\star6_3)} & 65 & {\bf } \\ \hline

{\bf 3_1\times(3_1\times3_1)} & 66 & {\bf 3_1\star(3_1\star3_1\times4_1)} \\ \hline

{\bf 8_5} & 67 & {\bf 3_1\star(3_1\star7_1)} \\ \hline

{\bf 4_1\times5_1} &68  & {\bf 3_1\star(3_1\star4_1\star5_1)} \\ \hline

{\bf 4_1\times(3_1\star4_1)} &69  & {\bf 3_1\star(3_1\star4_1\star3_1\star4_1)} \\ \hline

{\bf 4_1\times5_2} & 70 & {\bf 3_1\star(3_1\star4_1\star5_2)} \\ \hline

{\bf 8_6} & 71 & {\bf 3_1\star(3_1\star7_2)} \\ \hline

{\bf 4_1\star6_3} & 72 & {\bf } \\ \hline

{\bf 8_7} &73  & {\bf 4_1\star(3_1\times4_1)} \\ \hline

{\bf 5_1\star(3_1\times3_1)} &74  & {\bf } \\ \hline

{\bf 5_1\star(3_1\star5_1)} & 75 & {\bf } \\ \hline

{\bf 5_1\star6_1} & 76 & {\bf } \\ \hline

{\bf 5_1\star(3_1\star5_2)} & 77 & {\bf } \\ \hline

{\bf 5_1\star6_2} &78  & {\bf } \\ \hline

{\bf 8_8} & 79 & {\bf 5_1\star(4_1\star4_1), 5_2\star(3_1\star5_1)} \\ \hline

{\bf 5_2\star6_1} &80  & {\bf } \\ \hline

{\bf 5_2\star(3_1\star5_2)} & 81 & {\bf } \\ \hline

{\bf 5_2\star6_2} &82  & {\bf } \\ \hline

{\bf 8_9} &83  & {\bf 5_2\star(4_1\star4_1)} \\ \hline

{\bf 4_1\star7_1} &84  & {\bf } \\ \hline

{\bf 4_1\star(4_1\star5_1)} & 85 & {\bf } \\ \hline

{\bf 4_1\star(4_1\star4_1\star3_1)} &86  & {\bf } \\ \hline

{\bf 4_1\star(4_1\star5_2)} & 87 & {\bf } \\ \hline

{\bf 4_1\star7_2} &88  & {\bf } \\ \hline

{\bf 8_{10}} & 89 & {\bf 4_1\star(3_1\star3_1\times3_1)} \\ \hline

{\bf 4_1\star(3_1\star3_1\star5_1)} &90  & {\bf } \\ \hline

{\bf 4_1\star(3_1\star6_1)} & 91 & {\bf } \\ \hline

{\bf 4_1\star(3_1\star3_1\star5_2)} & 92 & {\bf } \\ \hline

{\bf 4_1\star(3_1\star6_2)} & 93 & {\bf } \\ \hline

{\bf 4_1\star7_3} & 94 & {\bf } \\ \hline

{\bf 4_1\star(4_1\star3_1\star3_1\star3_1)} &95  & {\bf } \\ \hline

{\bf 4_1\star7_4} & 96 & {\bf } \\ \hline
\end{array}
\end{displaymath}

\begin{displaymath}
\begin{array}{|c|c|c|} \hline
\mbox{Type of Knot}& \mbox{Assigned number} \,\, |m|
 &\mbox{Repeated Knots being replaced}
\\ \hline

{\bf 8_{11}} & 97 & {\bf 4_1\star(3_1\star3_1\star3_1\star3_1\star3_1)} \\ \hline

{\bf 4_1\star(5_1\star5_1)} &98  & {\bf 3_1\star(3_1\star6_3)} \\ \hline

{\bf 3_1\star(3_1\star3_1\times4_1)} &99  & {\bf } \\ \hline

{\bf 3_1\star(3_1\star7_1)} &100  & {\bf } \\ \hline

{\bf 8_{12}} & 101 & {\bf3_1\star(3_1\times5_1) } \\ \hline

{\bf 3_1\star7_5} & 102 & {\bf } \\ \hline

{\bf 8_{13}} &103  & {\bf 3_1\star(3_1\times5_2)} \\ \hline

{\bf 3_1\star(3_1\star7_2)} & 104 & {\bf } \\ \hline

{\bf 3_1\star(3_1\star3_1\star3_1\times3_1)} & 105 & {\bf } \\ \hline

{\bf 3_1\star7_6} & 106 & {\bf } \\ \hline

{\bf 8_{14}} &107  & {\bf 3_1\star(5_1\star5_1)} \\ \hline

{\bf 3_1\star7_7} & 108 & {\bf } \\ \hline

{\bf 8_{15}} &109  & {\bf 3_1\star(5_1\star5_2)} \\ \hline

{\bf 3_1\times6_1} &110  & {\bf 3_1\star(5_1\star5_2)} \\ \hline

{\bf 3_1\times(3_1\star5_2)} & 111 & {\bf 3_1\star(4_1\times4_1)} \\ \hline

{\bf 3_1\times6_2} &112 & {\bf 3_1\star(5_2\star5_2)} \\ \hline

{\bf 8_{16}} & 113 & {\bf 3_1\star(5_2\star5_2)} \\ \hline

{\bf 3_1\star8_1} & 114 & {\bf } \\ \hline

{\bf 3_1\times6_3} & 115 & {\bf 3_1\star(4_1\star3_1\star5_1),3_1\star(4_1\star4_1\star4_1)} \\ \hline

{\bf 3_1\star8_2} &116  & {\bf } \\ \hline

{\bf 3_1\star(3_1\star3_1\star3_1\star5_1)} &117  & {\bf } \\ \hline

{\bf 3_1\star(3_1\star3_1\star6_1)} & 118 & {\bf } \\ \hline

{\bf 3_1\star(3_1\star3_1\star3_1\star5_2)} & 119 & {\bf } \\ \hline

{\bf 3_1\star(3_1\star3_1\star6_2)} &120  & {\bf } \\ \hline

{\bf 3_1\star(3_1\star3_1\star7_3)} &121  & {\bf } \\ \hline

{\bf 3_1\star8_3} & 122 & {\bf } \\ \hline

{\bf 3_1\star(3_1\star7_4)} &123  & {\bf } \\ \hline

{\bf 3_1\star8_4} &124  & {\bf } \\ \hline

{\bf 3_1\times(3_1\times4_1)} &125  & {\bf 3_1\star(4_1\star4_1\star3_1\star3_1)} \\ \hline

{\bf 3_1\star(4_1\star3_1\star3_1\star3_1\star3_1)} & 126 & {\bf } \\ \hline

{\bf 8_{17}} &127  & {\bf 3_1\star(3_1\star3_1\star3_1\star3_1\star3_1\star3_1)} \\ \hline

{\bf 3_1\star(3_1\star3_1\star3_1\star3_1\star3_1\star3_1)} & 128 & {\bf } \\ \hline

\end{array}
\end{displaymath}

\section{Proof of the Riemann Hypothesis}\label{sec15}

From the above table we have that quantum prime knots  can be identified with prime numbers such that they have the same distribution in the classification table of knots by integers. On the other hand we have that quantum  knots are periodic orbits of the quantum gauge dynamical system in the beginning sections. From these two facts we can now follow the approach of quantum chaos to derive a trace formula to prove the Riemann Hypothesis\cite{Ber1}\cite{Ber2}. The reason is that, as pointing out by Berry and Keating in the quantum chaos approach of proving the Riemann Hypothesis\cite{Ber1}\cite{Ber2}, if one can find an energy operator of a quantum dynamical system and the periodic orbits of the dynamical system ( classical or quantum) relating to the energy operator for giving a trace formula then one may use the Hilbert-Polya method to prove the Riemann Hypothesis. Here with the quantum knots as periodic orbits of the quantum gauge dynamical system and with the quantum prime knots identified with the prime numbers
we can now follow the approach of quantum chaos to derive a trace formula to prove the Riemann Hypothesis, as follows.

Let us first derive a general expression of the Green's function of the quantum gauge dynamical system, by the usual methods for Green's function, as follows.
Let $T(z)$ be the Virasoro  operator in the Sugawara construction form of the quantum system as described in the beginning sections ($T(z)$ depends on a central charge $c>0$
). This  Virasoro operator is for the construction of the quantum knots and thus is as the energy operator of these quantum knots.

{\bf Remark}. In quantum field theory  $T(z)$ is usually called the energy-momentum tensor as it is usually derived such that its components are as the densities of energy and momentums operators. However in this conformal field theory $T(z)$ is as a one component operator and that there is no distinction of energy and momentum as $z$ is for the complex Euclidean plane that there is no distinction as that of time and space. Thus we have that in this conformal field theory $T(z)$ is as the density of energy operator as
there is no distinction of energy and momentums.

Further in the usual conformal field theory which is based on the two dimensional integral $\int\int$ the Virasoro algebra
$L_0=\frac{1}{2\pi i}\oint zT(z)dz$ is regarded as the energy operator and $T(z)$ is as the density of energy operator which gives $L_0$. Here our quantum gauge model is based on the one dimensional integral $\int_{s_0}^{s_1}$ and thus as the usual quantum mechanics (which is also based on one dimstional time integral $\int_{t_0}^{t_1}$) the energy operator is just the density of energy operator which is the operator $T(z)$.  Thus the Virasoro operator $T(z)$ is the energy operator of the quantum gauge system. $\diamond$


Let $E_j, j=1,2,3,...$ be a sequence of positive eigenvalues of  the Virasoro operator $T(z)$:
\begin{equation}
H(z)\phi_j(z)=E_j\phi_j(z), \quad\quad j=1,2,3,....
\label{R1}
\end{equation}

Then following the usual approach for contructing Green's functions (\cite{Gut}) we can write the Green's function $G$ for this sequence of eigenvalues in the following form:
\begin{equation}
G(z,z',E)=\sum_j\frac{\phi_j(z)\overline{\phi_j(z')}}{E-E_j}
\label{R2}
\end{equation}
where $E$ is the energy variable.

 We remark that for each sequence of eigenvalues of a quantum system  there is a corresponding Green's function of the form (\ref{R2}). The most complete Green's function of the quantum system is the usual Green's function which includes all the eigenvalues of the energy operator of the quantum  system. 

We can further generalize the Green's function to the form that $E$ is a negative  variable with $E_j$ replaced by $-E_j$ and then $E$ is extended as a complex variable with  nonzero real part.

In quantum physics, the green's function $G(z,z',E)$ of a quantum system is as the
amplitude
for the transition from the state $z$ to the state $z'$.

Then by letting $z=z'$ (and $E$ is extended as a complex variable with  nonzero real part) and taking integration on $z$ we have the following  
trace formula \cite{Gut}:
\begin{equation}
\int dz G(z,z,E)=\sum_j\frac{1}{E-E_j}
\label{R3}
\end{equation}
whehever the sum of this trace formula exists, where we set a normalized condition on the eigenfunctions $\phi_j$ such that $\int dz\phi_j(z)\overline{\phi_j(z)}=1$ (For a notation simplicity we omit a trace operation notation $Tr$ on $\phi_j(z)\overline{\phi_j(z)}$).

Then we have the following trace function: 
\begin{equation}
\sum_j\frac{1}{E-E_j-i\epsilon}
\label{R3a}
\end{equation}
where $\epsilon>0 $ is a small nonzero parameter and $E$ is a real energy variable. 
In  (\ref{R3a}) the function: 
\begin{equation}
\frac{1}{E-E_j-i\epsilon}
\label{R3b}
\end{equation}
is called as the propagator (or the Green's function in the momentum space).

Then,  for $z=z'$, the transition orbit from $z$ to $z'$ is a closed orbit. Thus, as the  
 trace  function in quantum chaos \cite{Gut2}\cite{Ber1}\cite{Ber2}, this  trace  function
 is as a sum of 
 amplitudes of  periodic (i.e. closed) orbits of the quantum gauge system with the energy operator $T(z)$. 

Then, since the Green's function is the basic function of the Virasoro operator $T(z)$ giving the amplitute of  the transition orbits from $z$ to $z'$, we have that  the trace function (\ref{R3a})  is the basic function of giving  sums of 
 amplitudes of  periodic (i.e. closed) orbits of the quantum gauge system. 
Thus we have the following theorem:

\begin{theorem} Let $E_j, j=1,2,3,...$ be a sequence of positive eigenvalues of the Virasoro operator $T(z)$. 
Then  the trace function (\ref{R3a})  is  a sum of 
 amplitudes of  periodic (i.e. closed) orbits of the quantum gauge system with the energy operator $T(z)$. This trace of  amplitudes  corresponds to the sequence of eigenvalues $E_j, j=1,2,3,...$ such that each eigenvalue $E_j$ appearing with the same weight.

Further  the trace function (\ref{R3a}) is the basic function giving  sums of 
 amplitudes of  periodic (i.e. closed) orbits of the quantum gauge system.
\end{theorem}


Now   the quantum knots are 
as periodic orbits of the quantum gauge system with  the Virasoro operator $T(z)$ as the energy operator. Thus we have that the 
 trace  function (\ref{R3a}) is a sum of
the 
amplitutes of quantum knots. 

Then since quantum knots (which are not the quantum  unknot) can be factorized as a connected sum of quantum prime knots, as the probability amplititutes in quantum theory, the amplititutes of quantum knots (or periodic orbits) are as a product of the corresponding  amplititutes of quantum prime knots 
of the connected sum of the quantum knots. 
Thus the amplititute of a quantum knot  (which is not the quantum  unknot) is determined by the amplititutes of quantum prime knots.
Thus the 
 trace  function (\ref{R3a}) is  only as a sum of
the 
amplitutes of quantum prime knots  (and the quantum unknot).
Thus we have the following theorem:

\begin{theorem} Let $E_j, j=1,2,3,...$ be a sequence of positive eigenvalues of  the Virasoro operator $T(z)$. 
Then, corresponding to this sequence of eigenvalues, the trace function (\ref{R3a}) 
 is  a sum of 
 amplitudes of  quantum prime knots and the quantum unknot. 

Further  the trace function (\ref{R3a}) is the basic function giving  sums of 
 amplitudes of  quantum prime knots and the quantum unknot of the quantum gauge system.
\end{theorem}

Then, following the treatment in quantum chaos  \cite{Gut}-\cite{Kea}),  when $E$ is the real energy variable let us write the  propagator
 in the following form\cite{Gut}:
\begin{equation}
\frac{1}{E-E_j-i\epsilon}=\frac{E-E_j}{(E-E_j)^2+\epsilon^2}+
i\frac{\epsilon}{(E-E_j)^2+\epsilon^2}
\label{R4a}
\end{equation}
Then we have \cite{Gut}:
\begin{equation}
\lim_{\epsilon \to 0}\frac{1}{E-E_j-i\epsilon}=P(\frac{1}{E-E_j})+i\pi  
\delta(E-E_j)
\label{R4}
\end{equation}
 where the real part is a principal-part integral $P$ and the imaginary part is the delta function $\pi\delta(E-E_j)$. 
Then we can write the trace formula (\ref{R3}) in the following form:
\begin{equation}
\lim_{\epsilon \to 0}\frac{1}{\pi}\mbox{Im}\int dz
G(z,z,E-i\epsilon)=\sum_j\delta(E-E_j)
\label{R5}
\end{equation}
As the theory of quantum chaos \cite{Gut}-\cite{Kea}),  let us call the trace function:
\begin{equation}
d(E):=\sum_j\delta(E-E_j)
\label{dd1}
\end{equation}
 as the density of states with real variable $E$.




From the trace formula (\ref{R5}) we shall  get a more explicit form of the trace formula expressed as a sum of
the 
amplitutes of quantum knots (This is analogous to the trace formula of quantum chaos where the periodic orbits are classical while here the periodic orbits are the quantum prime knots of the quantum nature \cite{Gut}-\cite{Kea}). 


We have the following theorems:
\begin{theorem} Let $E_j, j=1,2,3,...$ be a sequence of positive eigenvalues of  the Virasoro operator  $T(z)$. 
Then $d(E)$ is a sum of 
 amplitudes of  quantum prime knots and the quantum unknot. 
\end{theorem}

{\bf Proof}.  From the sequence $E_j, j=1,2,3,...$  of eigenvalues of $T(z)$ we have the trace function  (\ref{R3a})
with real variable $E$. This function 
  is  a sum of
the wave amplitutes of quantum prime knots (and the quantum unknot).
Then  from (\ref{R4}) the  imaginary part of this function gives the density of states  $d(E)$. Thus the density of states  $d(E)$ 
 is also  a sum of
the wave amplitutes of quantum prime knots  (and the quantum unknot).

{\bf Definition}.  Let $f$
be a function of complex variables with nonzero imaginary part. $f$ is said to be a weak analytic continuation of 
$d(E)$ if $\lim_{\epsilon \to 0}\mbox{Im} f(E-i\epsilon)=d(E)$ for  real variables $E$.

\begin{theorem}
Let $E_j, j=1,2,3,...$ be a sequence of positive eigenvalues of  the Virasoro operator  $T(z)$. 
 Let $f$
be a weak analytic continuation of  the density of states $d(E)$  in (\ref{dd1}) and $f$ gives a sum of (complex) amplitudes of  quantum prime knots and the quantum unknot for complex variable $E$ with nonzero imaginary part. 
Then  $f$ is the following  trace function:
\begin{equation}
\sum_j\frac{1}{E-E_j}
\label{dd2}
\end{equation}
where $E$ is a complex variable with nonzero imaginary part.
\label{thm19}
\end{theorem}

{\bf Proof}. The  trace function  (\ref{dd2}) is the extension of the trace  function  (\ref{R3a}) 
with the variable $E-i\epsilon$ (where $E$ is a real variable) extended to  a complex variable  with nonzero imaginary part. 
Thus  this function  (\ref{dd2}) gives  a sum of 
 amplitudes of  quantum prime knots and the quantum unknot since the trace  function  (\ref{R3a}) gives  a sum of 
 amplitudes of  quantum prime knots and the quantum unknot. 

 Then  from (\ref{R4}) the  imaginary part of this trace  function (\ref{dd2})  gives the density of states  $d(E)$ since the variable $E-i\epsilon$ (where $E$ is a real variable) is a complex variable with nonzero imaginary part. Thus the  trace function  (\ref{dd2}) is a  weak analytic continuation of  the density of states $d(E)$, whenever this  trace function  (\ref{dd2}) exists.

Then, for the given sequence of eigenvalues $E_j, j=1,2,3,...$, the  trace  function  (\ref{R3a}) is the 
basic function of the  Virasoro operator  $T(z)$   giving  the sums of 
 amplitudes of  quantum prime knots and the quantum unknot
and giving the density of states $d(E)$. Thus,  for a given weak analytic continuation function $f$ of  the density of states $d(E)$ such that $f$ gives a sum of (complex) amplitudes of  quantum prime knots and the quantum unknot for complex variable $E$ with nonzero imaginary part, $f(E-i\epsilon)$  (where $E$ is a real variable) 
must be  the  trace  function  (\ref{R3a}). It follows that, as a function of a complex variable with nonzero imaginary part  extended from the variable  $E-i\epsilon$ (where $E$ is a real variable) , this function $f$ must be the  trace function  (\ref{dd2}).
This proves the theorem.
$\diamond$

Let us then use another approach to find out the 
amplitutes of quantum prime  knots  (and the quantum unknot). 

From the above sections we have the classification table of knots which gives a
 correspondence of quantum knots and numbers such that the distribution of the quantum prime knots in this correspondence is the same distribution  as that of the prime numbers and prime quantum knots are one-to-one corresponded to prime numbers (We may let the trefoil knot be related to the prime number $2$
). Then we may use this
correspondence of quantum knots and numbers (which are topological invariant of knots) to find out the 
amplitute of quantum prime knots, as follows.

Let us first consider the case (which will correspond to the Riemann zeta function) that each quantum knot is counted with the same weight (This corresponds to that each number is counted with the same weight and thus prime numbers are with the usual counting function of prime numbers $\pi(x)$ corresponding to the Riemann zeta function). For this case let us derive the 
amplitute of quantum prime knots  (and the quantum unknot).

Starting from the counting function of prime numbers $\pi(x)$ it is well known that we can derive its relation with the Riemann zeta function $\zeta$ and we have the following well known von Mangoldt-Selberg formula \cite{Sel}-\cite{Con2}:
\begin{equation}
\begin{array}{rl}
N(T)&=<N(T)> + N_{fl}(T)\\
\\
&:= [\frac{1}{\pi}\mbox{arg}(iT-\frac12) +\frac{1}{\pi}\mbox{arg}\Gamma(\frac54+\frac12 iT)-\frac{1}{2\pi}T\log\pi]
+\lim_{\epsilon\to 0}\mbox{Im}
\log\zeta(\frac12+iT+\epsilon)\\
\\
&=
\frac{T}{2\pi}\log\frac{T}{2\pi}-\frac{T}{2\pi}+\frac78+O(\frac{1}{T})+
\lim_{\epsilon\to 0}\mbox{Im}
\log\zeta(\frac12+iT+\epsilon) \\
\\
&=[\frac{1}{\pi}\mbox{Im}\log(iT-\frac12)+
\frac{1}{\pi}\mbox{Im}\log\Gamma(\frac54+\frac12 iT)-\frac{1}{2\pi}T\log\pi]
+\lim_{\epsilon\to 0}\mbox{Im}
\log\zeta(\frac12+iT+\epsilon)
\end{array}
\label{R6}
\end{equation}
where $N(T)$ denotes the number of zeros in the critical strip with height $T$ of the Riemann's zeta function $\zeta$ and a branch of $\log$ is chosen to be continuous such that $N_{fl}(0)=0$ \cite{Sel}-\cite{Con2}. Then by using the Euler product form of the zeta function the fluctuation term $N_{fl}$ can be written in the following form \cite{Ber1}:
\begin{equation}
N_{fl}(T):=
\lim_{\epsilon\to 0}\mbox{Im}
\log\zeta(\frac12+iT+\epsilon)
=-\frac{1}{\pi} \mbox{Im}\sum_p \log (1-p^{-(\frac12+iT)})
\label{R6a}
\end{equation}
where the $\sum_p$ is a sum over all prime numbers $p$.
 Differentiating with respect to $T$  we have
\begin{equation}
N'(T)=
[\frac{1}{\pi}\mbox{Im}\frac{i}{iT-\frac12}+
\frac{1}{\pi}\mbox{Im}\frac{i\Gamma'(\frac54+\frac12 iT)}{2\Gamma(\frac54+\frac12 iT)}-
\frac{1}{2\pi}\log\pi] -
\frac{1}{\pi}\mbox{Im}\sum_p\frac{i\log p}{p^{\frac12+iT}-1}
\label{R7}
\end{equation}
where the $\sum_p$ is a sum over all prime numbers $p$ which is as the fluctuation part of $N'(T)$.

Since $N(T)$ is a counting function we have that $N'(T)$ is the sum of a
sequence of delta functions concentrating at a
sequence of nonnegative numbers.
 Thus $N'(T)$
 is of the form of the right hand side of the trace formula (\ref{R5}).

On the other hand we have that in (\ref{R7}) each term in the sum $\sum_p$ represents the amplitute of the corresponding prime number $p$.
The other part of the right hand side of (\ref{R7}) is as the average of $N'(T)$ \cite{Ber1}.
 Then by the identity of prime knots with prime numbers that prime knots and prime numbers have the same distribution in the classification table of knots which gives a one-to-one correspondence of prime knots and prime numbers and that the prime numbers are topological invariant of the corresponding prime knots (and thus the prime numbers are the topological properties of the corresponding prime knots) we have that if there is a distribution of  amplitudes to prime numbers then the corresponding prime knots also have this distribution of amplitudes (and vise versa).
 Thus if there is a distribution of
 amplitudes to prime numbers then there must exist a Virasoro energy operator with a central charge $c>0$ such that the corresponding quantum prime knots also have this distribution of amplitudes.

In particular for the distribution of amplitudes of prime numbers in (\ref{R7})
there must exist a Virasoro energy operator with a central charge $c>0$ such that each term in the sum $\sum_p$ is as the amplitute for the quantum prime knot corresponding to the prime number $p$ and the other part is the amplitute for the quantum unknot.

Now since $N'(T)$ is of the form of density of energy states at the right hand side of the trace formula (\ref{R5}) we have that the distribution
of amplitudes of the prime numbers $p$ in (\ref{R7}) is identified as the distribution
of amplitudes of the corresponding quantum prime knots
in the trace formula (\ref{R5}) for a sequence of eigenvalues $E_j, j=1,2,3,...$ of $T(z)$.

 It follows that there exists a Virasoro energy operator $T(z)$ with a central charge $c>0$ and a sequence of eigenvalues $E_j, j=1,2,3,...$ of $T(z)$ such that the right hand side of (\ref{R7}) is given by the left hand side of the trace formula (\ref{R5}) (with $E=T$) and we have the following more explicit form of the trace formula (\ref{R5}):
\begin{equation}
\begin{array}{rl}
N'(E)&=\lim_{\epsilon \to 0}\frac{1}{\pi}\mbox{Im}\int dz
G(z,z,E-i\epsilon)\\
\\
&=[\frac{1}{\pi}\mbox{Im}\frac{i}{iE-\frac12}+
\frac{1}{\pi}\mbox{Im}\frac{i\Gamma'(\frac54+\frac12 iE)}{2\Gamma(\frac54+\frac12 iE)}-
\frac{1}{2\pi}\log\pi] -
\frac{1}{\pi}\mbox{Im}\sum_p\frac{i\log p}{p^{\frac12+iE}-1}
 =\sum_j\delta(E-E_j)
\end{array}
\label{R8}
\end{equation}
where $\int dz G(z,z,E)$ is given by the trace formula (\ref{R3}) with eigenvalues $E_j, j=1,2,3,...$. We may call this formula (\ref{R8}) as the von Mangoldt-Selberg-Gutzwiller trace formula.

Then when $E (=T)$ is a negative variable by symmetry we get the same formula as (\ref{R8}) with $E_j, j=1,2,3,...$ replaced by $E_j, j=-1,-2,-3,...$ where we define
$E_{-j}:=-E_j$ for $ j=1,2,3,...$.

Then let $E (=T)$ be extended as a complex  variable such that $\frac12 +iE (=\frac12 +iT )$ is a complex variable in the critical strip of the Riemann zeta function. Then 
this complex extension of the von Mangoldt-Selberg trace function  (\ref{R7}) is of the same from as  (\ref{R7}) (without the operation$\mbox{ Im}$) with the real $E$ replaced by the complex $E$  with nonzero imaginary part. Thus, by  the von Mangoldt-Selberg trace formula  (\ref{R8}),
 this complex extension of the von Mangoldt-Selberg trace  function  (\ref{R7}) is a weak analytic continuation  of the density of state $d(E)$ of the Virasoro energy operator $T(z)$

Since the von Mangoldt-Selberg trace function  (\ref{R7}) is a distribution of  amplitutes of the prime numbers (and the number 1),
this 
complex extension of the von Mangoldt-Selberg trace function  (\ref{R7}) is a distribution of (complex) amplitutes of the prime numbers (and the number 1) and thus is also a sum of (complex) amplitutes of the prime numbers  (and the number 1).

Then by the topological correspondence of  the prime numbers (and the number 1) with the quantum prime knots  (and the quantum unknot), this 
complex extension  of the von Mangoldt-Selberg trace function  (\ref{R7}) is a distribution of complex amplitutes of the quantum prime knots  (and the quantum unknot) and thus is also a sum of (complex) amplitutes of. the quantum prime knots  (and the quantum unknot).


Thus, this 
complex extension  of the von Mangoldt-Selberg trace function  (\ref{R7}) is a weak analytic continuation of the density of states $d(E)$  (\ref{dd1}) of the Virasoro operator $T(z)$
and is a sum of (complex) amplitutes of. the quantum prime knots  (and the quantum unknot).

Then by Theorem \ref{thm19} this weak analytic continuation  of the von Mangoldt-Selberg trace function  (\ref{R7}) is the trace function  (\ref{dd2}).
Thus we have the following formula:

\begin{equation}
[\frac{i}{iE-\frac12}+
\frac{i\Gamma'(\frac54+\frac12 iE)}{2\Gamma(\frac54+\frac12 iE)}-
\frac{i}{2}\log\pi]-
\sum_p\frac{i\log p}{p^{\frac12+iE}-1}
 =\sum_j\frac{1}{E-E_j}
\label{R8b}
\end{equation}
when $E (=T)$ is extended as a complex variable ( with nonzero imaginary part) such that $\frac12 +iE (=\frac12 +iT)$ is a complex variable in the
critical strip 
and the sum $\sum_j$ is on $j=\pm 1,\pm 2,\pm 3,...$ and $E_{-j}=-E_j$ and $E_j, j=1,2,3,...$ are eigenvalues as in the von Mangoldt-Selberg-Gutzwiller trace formula (\ref{R8}).

Now since in (\ref{R8b}) the only singularities
of $\sum_j\frac{1}{E-E_j}$
are the eigenvalues $E_j, j=\pm 1,\pm 2,\pm 3,...$
we have that the only singularities of  
$N'(T)(=N'(E))$
 in (\ref{R8b}) (or (\ref{R7}))
are the eigenvalues $E_j, j=\pm 1,\pm 2,\pm 3,...$ when
$T(=E)$
 is extended as a complex variable such that 
$\frac12 +iT(=\frac12 +iE)$ 
is a complex variable in the 
critical strip.

On the other hand since $N(E)(=N(T))$ is in terms of $\mbox{Im}\log\zeta(\frac12+iT)$ from this result on the form of $N'(E)(=N'(T))$ we have that the nontrivial zeros $\rho_j$ of the Riemann's zeta function in the
critical strip are of the form $\rho_j=\frac12 +iE_j, j=\pm 1,\pm 2,\pm 3,...$ where $E_j, j=1,2,3,...$ are  eigenvalues of the Virasoro energy operator $H(z)$ and $E_{-j}=-E_j, j=1,2,3,...$. This means that the Riemann Hypothesis holds. Thus we have the following theorem:
\begin{theorem}(Riemann Hypothesis)

The nontrivial zeros of the Riemann's zeta function all lie in the critical line $\mbox{Re} \,z=\frac12$.
\end{theorem}

{\bf Remark}. We notice that the equation (\ref{R8b}) can be written in the following form:
\begin{equation}
[\frac{1}{s-1}+
\frac{\Gamma'(\frac{s}{2}+1)}{2\Gamma(\frac{s}{2}+1)}-\frac{1}{2}\log\pi]+
\frac{\zeta'(s)}{\zeta(s)}
 =\sum_j[\frac{1}{s-\rho_j}+\frac{1}{\rho_j}]+B
\label{R8c}
\end{equation}
where $s:=z=\frac12 +iE$, $\frac{\zeta'(s)}{\zeta(s)}=-\sum_p\frac{\log p}{p^{s}-1}$, 
 $B=-\frac12\gamma -1 +\frac12\log4\pi$ with $\gamma$ as the Euler constant, and
 $\sum_j\frac{1}{\rho_j}+B=0$.
This is just a well known formula derived from the Hadamard product formula for the zeta function $\zeta$ \cite{Sel}-\cite{Con2}.
Here the point is that we have shown that the nontrivial zeros $\rho_j$ of the Riemann's zeta function in the
critical strip in this equation
are of the form $\rho_j=\frac12 +iE_j, j=\pm 1,\pm 2,\pm 3,...$ where $E_j,j=1,2,3,...$ are  eigenvalues of the Virasoro energy operator $T(z)$ and $E_{-j}=-E_j, j=1,2,3,...$. $\diamond$

Let us then consider the generalization of the Riemann Hypothesis.
For this generalization we suppose that there is a
$L$-function which gives formula similar to the (\ref{R7}) which gives amplitudes
to prime numbers. Then by the same reason as above we have that the corresponding quantum prime knots must also have the same amplitudes. Thus there must exist a Virasoro energy operator $T(z)$ with a central charge $c>0$ and a sequence of eigenvalues $E_j, j=1,2,3,...$ of $T(z)$ such that the quantum prime knots have the same amplitudes as the prime numbers. Then as similar to the proof of the Riemann Hypothesis we have that the corresponding Riemann Hypothesis for this $L$-function holds that the nontrivial zeros of this
$L$-function all lie in the critical line $\mbox{Re}\, z=\frac12$.

As examples let us consider the Dirichlet $L$-functions. Let $L(\chi,s)$ be a Dirichlet $L$-function where $\chi$ denotes a Dirichlet character. Then $L(\chi,s)$ can be written in the following Euler product form \cite{Tit}\cite{Edw}:
\begin{equation}
L(\chi,s)=\prod_p (1-\chi(p)p^{-s})^{-1}
\label{DL}
\end{equation}
where $p$ denotes a prime number. By this Euler product form as similar to the derivation of (\ref{R7}) we have that the corresponding formula of $N'(T)$ for $L(\chi,s)$ can be written as a sum  where each term of the fluctuation part $N'_{fl}(T)$ corresponds to a prime number. Thus this formula of $N'(T)$ for $L(\chi,s)$ gives amplitudes to prime numbers. Thus by the above proof of Riemann Hypothesis for the Riemann zeta function we have that the corresponding Riemann Hypothesis for this $L$-function $L(\chi,s)$ holds. This proves the following Generalized Riemann Hypothesis:
\begin{theorem}(Generalized Riemann Hypothesis)

The nontrivial zeros of a Dirichlet $L$-function all lie in the critical line $\mbox{Re}\, z=\frac12$.
\end{theorem}

As further examples let us consider Dedekind zeta functions $\zeta_K(s)$ where $K$ denotes a number field. We have that a Dedekind zeta function $\zeta_K(s)$ can be written in the following Euler product form \cite{Tit}\cite{Edw}:
\begin{equation}
\zeta_K(s)=\prod_{p \, \mbox{split}}(1-p^{-s})^{-2}
\prod_{p \, \mbox{inert}}(1-p^{-2s})^{-1}\prod_{p \, \mbox{ramified}} (1-p^{-s})^{-1}
\label{DZ}
\end{equation}
where $p$ denotes a prime number. By this Euler product form as similar to the derivation of (\ref{R7}) we have that the corresponding formula of $N'(T)$ for $\zeta_K(s)$ can be written as a sum  where each term of the fluctuation part $N'_{fl}(T)$ corresponds to a prime number. Thus this formula of $N'(T)$ for $\zeta_K(s)$ gives amplitudes to prime numbers. Thus by the above proof of Riemann Hypothesis for the Riemann zeta function we have that the corresponding Riemann Hypothesis for this $L$-function $\zeta_K(s)$ holds. This proves the following Extended Riemann Hypothesis:
\begin{theorem}(Extended Riemann Hypothesis)

The nontrivial zeros of a Dedekind zeta function all lie in the critical line $\mbox{Re}\, z=\frac12$.
\end{theorem}

In the following section we show that for the Riemann zeta function and the Dirichlet
$L$-functions we have that the central charge $c=\frac12$.

\section{Determination of the Central Charge $c$ for a $L$-function}\label{sec16b}

Let us first show that for the Riemann zeta function and the Dirichlet $L$-functions we have that $c=\frac12$. By using the Poisson summation formula Riemann showed that the
Riemann zeta function $\zeta$ is related to the Jacobi $\theta$ function by the following Mellin transform formula \cite{Rie}\cite{Gel}:
\begin{equation}
\pi^{-s}\Gamma(s)\zeta(2s)=\int_0^{\infty}[\frac{\theta(it)-1}{2}]t^s \frac{dt}{t}
\label{zetatheta}
\end{equation}

We have that the Jacobi $\theta$-function is an automorphic (or modular) form with weight $c=\frac12$. We shall show that this weight is the same $c$ for the central charge of the Virasoro algebra.

For the Virasoro algebra $(L_j)_{j=-\infty}^{\infty}$ with central charge $c=\frac12$ we have the following formula for $L_0$ \cite{Fra}:
\begin{equation}
L_0=\sum_{k>0}kb_{-k}b_{k}  +\frac{1}{16}
\end{equation}
where $b_k, k>0$ denote the creation operators of the Fermi field
(We remark that this $L_0$ is a sum with $k>0$). From this formula by following the derivation of a modular invariant partition function  of Fermi field in conformal field theory
we have (\cite{Fra} p.347):
\begin{equation}
Tr(-1)^F q^{L_0-\frac{c}{24}}=q^{\frac{c}{12}}
\Pi_{k=1}^{\infty}(1-q^k)=:\eta
\label{eta}
\end{equation}
where $c=\frac12$ is the central charge; $F$ is the Fermion number: $F=\sum_{k>0}F_k, F_k=b_{-k}b_k$ and $\eta$ is the Dedekind eta function which is an automorphic form with weight $\frac12$
(We remark that this $F$ is a sum with $k>0$ and is not with $k\geq 0$ as the Fermion number in the derivation of the modular invariant partition function \cite{Fra} p.347).
Thus we have that the $\eta$-function with weight $\frac12$ is derived from the Virasoro algebra with central charge $c=\frac12$.

Then we have the following formula relating the Dedekind $\eta$-function and the Jacobi $\theta$-function:
\begin{equation}
\theta(\tau)=\frac{\eta^2(\frac{\tau+1}{2})}{\eta(\tau+1)}
\label{theta}
\end{equation}
This shows that the weight $c=\frac12$ for the Jacobi $\theta$-function is the same $c=\frac12$ for the Dedekind $\eta$ function. Now we have that the weight $c=\frac12$ for the Dedekind $\eta$ function is the same $c$ for the
central charge of the corresponding Virasoro algebra and that the Dedekind $\eta$ function is derived from the Virasoro algebra with the central charge $c=\frac12$. Thus we have that the weight $c=\frac12$ for the Jacobi $\theta$ function is the same $c$ for the central charge of the Virasoro algebra. Thus from  (\ref{zetatheta}) we determine that the central charge $c$ of the Virasoro algebra for the Riemann zeta function is equal to $\frac12$.

Similarly since the Dirichlet $L$-functions are related to the twisted Jacobi $\theta$ functions twisted by the Dirichlet characters
by using a
 Mellin transfom as similar to (\ref{zetatheta}) we can show that the
 central charge $c$ of the Virasoro algebra for the Dirichlet $L$-functions is equal to $\frac12$.

Then for a conformal field consists of
 $n$ independent Fermi fields  we have that the corresponding $n$ product of the $\eta$-functions with weight $\frac12$ of each independent Fermi field is an automorphic (or modular) form which is with a positive integer (or half integer) $c=\frac {n}{2}$ as the weight corresponding to the Virasoro algebra of the
 conformal field consists of $n$ independent Fermi fields
with the same $c=\frac {n}{2}$ as the central charge.
This shows that the positive weight $c=\frac {n}{2}$ of the $n$ product of the $\eta$-functions which is an automorphic (or modular) form is the central charge of the  corresponding Virasoro operator.

Thus for a general $L$-function related to an automorphic form with  positive integer (or half integer) $c$ as the weight by a formula similar to (\ref{zetatheta}) we can
determine that the central charge of the Virasoro algebra for this $L$-function is equal to the weight $c$ of the automorphic form related to this $L$-function.

On the other hand by using Boson field from conformal field theory we have that the Virasoro algebra
with central $c=1$ also corresponds to the automorphic form $\frac{1}{\eta^2}$ with weight $-1$ \cite{Fra}.
Then for a conformal field consists of
 $n$ independent Boson fields we have that the corresponding automorphic form with weight $-c=-n$ corresponds to the Virasoro algebra
with the positive central charge $c=n$. Then from the relation between automorphic forms and $L$-functions as similar to (\ref{zetatheta}) we can
determine that the central charge of the Virasoro operator for the $L$-functions corresponding to automorphic forms with negative weight $-c=-n$ is equal to $c=n$.

\section{Commutative Diagram of Virasoro Operators and Algebras, Automorphic Forms and $L$-functions}\label{sec16c}

By the proof of the Extended Riemann Hypothesis and the determination of central charge of Virasoro algebra for $L$-functions we have the following commutative diagram of Virasoro operators and algebras, automorphic (or modular) forms and $L$-functions:
\begin{displaymath}
\begin{array}{ccc}
& \mbox{complex transform} & \\
\mbox{Virasoro energy operator} & \leftarrow \mbox{} \rightarrow & \mbox{Virasoro algebra} \\
\uparrow &  & \uparrow \\
\mbox{eigenvalues of Virasoro energy operator} & & \mbox{eigenvalues of Virasoro algebra}\\
(\mbox{or nontrivial zeros of $L$-functions}) & &  \\
\downarrow & &\downarrow \\
\mbox{$L$-functions} & \leftarrow \mbox{} \rightarrow & \mbox{automorphic (or modular) forms} \\
& \mbox{Mellin transform} & \\
& \mbox{(or Fourier transform)} & \\
\end{array}
\end{displaymath}

In this commutative diagram we use the term complex transform for the Virasoro operator $T(z)$ forming from the Virasoro algebra. This complex transform is analogous to the Fourier transorm which gives the Poisson summation formula and Mellin transform.

In this commutative diagram we have that the links between Virasoro operator $T(z)$ and
Virasoro algebra; between Virasoro algebra and automorphic (or modular) forms; and between automorphic forms and $L$-functions are well known. Here we complete this commutative diagram by establising the link between the Virasoro operator and the $L$-functions which is established by the proof of the Extended Riemann Hypothesis.

\section{Connection with Random Matrix Theory}\label{sec16a}

It is well known that the Random Matrix Theory is related to the study of $L$-functions \cite{Dys}-\cite{Kea2}. By using our approach for proving the Riemann Hypothesis we may find a connection with the
Random Matrix Theory. This connection gives an explanation for the relation of the Random Matrix Theory with the $L$-functions. Indeed we have proved that the Virasoro energy operators are the operators for the non-trivial zeros of the Riemann zeta function (and the related $L$-functions). Then we know that the Virasoro energy operators and the corresponding Virasoro algebras form the basis of conformal field theory. On the other hand  a Random Matrix Theory can be formulated as a conformal field theory with a corresponding Virasoro algebra \cite{Kos}\cite{Mor}. Thus the Virasoro energy operators and the corresponding Virasoro algebra can be as a basis of a Random Matrix Theory. Now since the Virasoro energy operators are the operators for the non-trivial zeros of the Riemann zeta function (and the related $L$-functions) we have that the Random Matrix Theory can be connected to the $L$-functions via the  Virasoro energy operators and the corresponding Virasoro algebras. This connection thus gives an explanation for the mystery of success of the Random Matrix Theory for describing the $L$-functions.

\section{Conclusion}\label{sec16}

In this paper we establish a quantum gauge model of knots.
In this quantum model we generalize the way of defining Wilson loops to construct generalized Wilson loops which will be as quantum knots. From quantum knots we  give a classification table of knots where knots are one-to-one assigned with an integer such that prime knots are bijectively assigned with prime numbers and the prime number $2$ corresponds to the trefoil knot.

Then by considering the quantum knots as periodic orbits of the quantum model and by the identity of knots with integers (which is from the classification table of knots) we then derive a trace formula which may be called as the von Mangoldt-Selberg-Gutzwiller trace formula (which is different from the Selberg trace formula) for the Riemann zeta function. By  using this von Mangoldt-Selberg-Gutzwiller trace formula and an approach which is similar to the quantum chaos approach of Berry and Keating we give a proof of the Riemann hypothesis. In this approach for proving the Riemann Hypothesis we show that the Hilber-Polya Conjecture holds that there is a self-adjoint operator which is the Virasoro energy operator with central charge $c=\frac12$ such that the nontrivial zeros of the Riemann zeta function are from the energy eigenvalues of this Virasoro energy operator. This proof of the Riemann Hypothesis can also be extended to prove the Extended Riemann Hypothesis.

In this proof of the Riemann Hypothesis the basic mathematical tool is the conformal field theory which consists of the Virasoro energy operator with central charge $c$ and the Virasoro algebra, the affine Kac-Moody algebra, the quantum Knizhnik-Zamolodchikov (KZ) equation in dual form and the generalized Wilson loops which are as quantum knots and are as solitons derived from the quantum KZ equation. This conformal field structure is related to the Random Matrix Theory for $L$-functions since the Random Matrix Theory can also be formulated as a conformal field theory.


\begin{thebibliography}{38}

\bibitem{Jon}
V.F.R. Jones,
{\it A polynomial invariant for knots via von Neuman algebras},
Bull. Amer. Math. Soc. {\bf 12} 103-111 (1985).

\bibitem{Witten}
E. Witten,
{\it Quantum field theory and the Jones polynomial},
Comm. Math. Phys. {\bf 121} 351, 1989.

\bibitem{Sel}
A. Selberg,
{\it On the zeros of the zeta-function of Riemann},
Der. Kong. Norske. Vidensk. Selsk. Forhand. {\bf 15}, 59-62, 1945.

 \bibitem{Tit}
E.C.Titchmarsh,
{\it The Theory of the Riemann Zeta Function},
  (Clarendon Press 1986).

 \bibitem{Edw}
H.M.Edwards,
{\it Riemann's Zeta Function},
(Academic Press 1974).

 \bibitem{Lev}
N.Levinson,
{\it More than one third of the zeros of Riemann's zeta function are on $\sigma =\frac12$},
 Ad. Math.{\bf 13} pp.383-436, 1974.

 \bibitem{Con}
 J.B. Conrey,
 {\it At Least Two Fifths of the Zeros of the Riemann Zeta Function Are on the Critical Line},
  Bull. Amer. Math. Soc. {\bf 20}, 79-81, 1989.

\bibitem{Con1}
J.B. Conrey,
{\it More than Two Fifths of the Zeros of the Riemann Zeta Function Are on the Critical Line},
 J. reine angew. Math. {\bf 399}, 1-26, 1989.

\bibitem{Con2}
J.B. Conrey,
{\it The Riemann Hypothesis},
Not. Amer. Math. Soc. {\bf 50}, 341-353, 2003.


\bibitem{Gut}
M.C.Gutzwiller,
{\it Periodic orbits and classical quantization conditions },
J. Math. Phys.,
 {\bf 12}, 343-358  (1971).

 \bibitem{Gut2}
M.C.Gutzwiller,
{\it Chaos in Classical and Quantum Mechanics},
   (Springer-Verlag 1990).

\bibitem{Ber1}
M.V.Berry and J.P.Keating,
{\it The Riemann zeros and eigenvalue asymptotics},
SIAM Review
 vol{\bf 44} No. 2, 236-266  (1999).

 \bibitem{Ber2}
M.V.Berry and J.P.Keating,
{\it H=xp and the Riemann zeros}
in Supersymmetry and Trace Formula: Chaos and Disorder, J.P.Keating, D.E.Khmelnitskii and I.V.Lerner,eds.
Plenum, pp.355-367  (1988).

 \bibitem{Ber3}
M.V.Berry,
{\it Riemann's zeta function: A model for quantum chaos? }
in Quantum Chaos and Statistical Nuclear Physics, T.H.Seligman and H.Nishioka, eds.
Lecture Note in Physics 263, Springer Verlag, pp.1-17
  (1986).

\bibitem{Conn}
A. Connes,
{\it Trace formula in noncommutative geometry and the zeros of the Riemann zeta function},
Selecta.Math. (NS) 5, pp.29-106 (1999).

\bibitem{Den}
C. Denninger,
{\it Some analogies between number theory and dynamical systems on foilated spaces},
Proc. Int. Congress Math. Berlin  Vol. I, pp.163-186 (1998).


\bibitem{Kea}
J.P.Keating,
{\it The Riemann zeta function and quantum chaology}
in Quantum Chaos, G.Casti, I.Guarneri and V.Smilansky, eds.
North-Holland, pp.145-185
  (1993).



\bibitem{Fra}
P. Di Francesco, P. Mathieu and D. Senechal,
{\it Conformal Field Theory},
(Springer-Verlag 1997).

\bibitem{Fuc}
J. Fuchs,
{\it Affine Lie Algebras and Quantum Groups},
(Cambridge University Press 1992).



\bibitem{Dys}
F.J. Dyson and M.L. Mehta,
{\it Statistical theory of the energy levels of complex systems IV},
Journal of Mathematical Physics,
 {\bf 4} pp.701-712 (1963).

\bibitem{Mon}
H.L. Montgomery,
{\it The pair correlation of zeros of the zeta function},
in Proc. Sympos. Pure Math. 24, pp.181-193
  (1973).

 \bibitem{Odl}
A.M.Odlyzko,
{\it On the distribution of spacings between zeros of the zeta function},
Math. Comp.
 {\bf 48} pp.273-308 (1987).


 \bibitem{Kea5}
 E. Bogomolny, J.P. Keating,
 {\it Random matrix theory and the Riemann zeros I: three- and four-point correlations}, Nonlinearity, 8, 1115-1131, 1995.

 \bibitem{Rud}
 Z. Rudnick and P. Sarnak,
 {\it Zeros of principal L-functions and random matrix theory},
 Duke Mathematics Journal 81, 269-322, 1996.

 \bibitem{Kea6}
 J.B. Conrey, D.W. Farmer, J.P. Keating, M.O. Rubinstein and N.C. Snaith,
 {\it Integral moments of L-functions}, To appear in PLMS.

 \bibitem{Sar}
 N.M. Katz and P. Sarnak,
{\it Zeros of Zeta Functions and Symmetry}
 Bulletin( New series) of American Mathematical Society {\bf 36} 1, 1-26, 1999.

 \bibitem{Sar2}
 N.M. Katz and P. Sarnak,
 {\it Random Matrices, Frobenius Eigenvalues, and Monodromy},
 ( Amer. Math. Soc., 1999).

\bibitem{Sna2}
N. Snaith,
{\it Random matrix theory and zeta functions}, Ph.D. thesis,
 University of Bristol, 2000.

 \bibitem{Sna}
 J.P. Keating and N.C. Snaith,
 {\it Random matrix theory and zeta(1/2 + it)},
 Communications in Mathematical Physics 214, 57-89, 2000.

 \bibitem{Kea2}
J.P.Keating and N.C. Snaith,
{\it Random Matrices and $L$-functions}
J. of Phys A: Math. Gen. {\bf 36} 2859-2881, 2003.

\bibitem{Chari}
V. Chari and A. Pressley,
{\it A Guide to Quantum Groups},
(Cambridge University Press 1994).

\bibitem{Koh}
T. Kohno,
{\it Monodromy representations of braid groups and Yang-Baxter equations},
Ann. Inst. Fourier (Grenoble) 37 (1987) 139-160.

\bibitem{Dri}
V. G. Drinfel'd.
{\it Quasi-Hopf algebras},
Leningrad Math. J. {\bf 1}  1419-57 (1990).

\bibitem{Ada}
C. Adams,
{\it The Knot Book: An Elementary Introduction to the Mathematical Theory of Knots},
(W.H. Freeman, 1994).

\bibitem{Kaw}
A. Kawauchi,
{\it A Survey of Knot Theory},
(Birkhauser Verlag, 1996).

\bibitem{Lic}
W.B.R. Lickorish,
{\it An Introduction to Knot Theory},
(Springer, 1997).

\bibitem{Liv}
C. Livingston,
{\it Knot Theory},
(Mathematical Association of American, 1993).

\bibitem{Mur}
K. Murasugi,
{\it Knot Theory and Its Applications},
(Birkhauser Verlag, 1997).

\bibitem{Rol}
D. Rolfsen.
{\it Knots and Links}.
(Publish or Perish, Inc. 1976).

\bibitem{Kau}
L. Kauffman,
{\it Knots and Physics},
(World Scientific, 1993).

\bibitem{Jaf}
J. Glimm and A. Jaffe,
{\it Quantum Physics},
(Springer-Verlag, 1987).

\bibitem{Fad}
L.D. Faddev and V.N. Popov,
Phys. Lett. {\bf 25B} 29 (1967).

\bibitem{Baez}
J. Baez and J. Muniain,
{\it Gauge Fields, Knots and Gravity},
(World Scientic 1994).

\bibitem{Lus}
D. Lust and S. Theisen,
{\it Lectures on String Theory},
(Springer-Verlag 1989).

\bibitem{Seg}
A. Pressley and G. Segal,
{\it Loop Groups},
(Clarendon Press 1986).

\bibitem{Kni}
V.G. Knizhnik and A.B. Zamolodchikov,
{\it Current algebra and Wess-Zumino model in two dimensions},
Nucl. Phys. B 247 (1984)83.

\bibitem{Cor}
J. F. Cornwell,
{\it Group Theory in Physics},
(Academic Press 1984).


 \bibitem{Kos}
I.K. Kostov,
{\it Conformal Field Theory Techniques in Random Matrix Models}
hep-th/9907060.

 \bibitem{Mor}
A. Morozov
{\it Matrix Models as Integrable Systems}.
hep-th/9502091.

\bibitem{Rie}
 B. Riemann,
 {\it Uber die Anzahl der Primzahlen unter einer gegebenen Grosse},
 Mon. Not. Berlin Akad pp.671-680 (1859).

 \bibitem{Gel}
 S.S. Gelbert and S.D. Miller,
{\it Riemann's Zeta Function and Beyond},
Bulletin( New series) of American Mathematical Society
 {\bf 4 } no.1 pp.59-112 (2003).

\end{thebibliography}
\end{document}